\documentclass[a4paper]{amsart}
\usepackage[all]{xy}
\usepackage{amsmath}
\usepackage{amsfonts,amssymb}

\usepackage{booktabs}
\usepackage{fixmath}

\numberwithin{equation}{section}
\newtheorem{satz}{Satz}
\newtheorem{proposition}[satz]{Proposition}
\newtheorem{theorem}{Theorem}
\newtheorem*{theorem*}{Theorem}

\newtheorem{lemma}[satz]{Lemma}

\newtheorem{corollary}[satz]{Corollary}

\theoremstyle{definition}
\newtheorem*{remark*}{Remark}
\newtheorem{remark}[satz]{Remark}
\newtheorem*{notation}{Notation}

\newtheorem{example}[satz]{Example}
\newtheorem*{example*}{Example}

\numberwithin{satz}{section}
\newcommand{\tensor}{\otimes}

\newcommand{\map}[1]{\stackrel{#1}{\longrightarrow}}

\newcommand{\un}[1]{\ensuremath{\protect\underline{#1}}}

\def\GL{\textsf{GL}}

\def\SL{\textsf{SL}}

\DeclareMathOperator{\Coh}{Coh}
\DeclareMathOperator{\Pic}{Pic}
\DeclareMathOperator{\Hecke}{Hecke}
\DeclareMathOperator{\Bun}{Bun}

\DeclareMathOperator{\uPic}{\underline{\Pic}}

\DeclareMathOperator{\Hom}{Hom}
\def\cHom{\mathcal{H}\mathit{om}}

\DeclareMathOperator{\Spec}{Spec}
\DeclareMathOperator{\Ext}{Ext}
\def\uExt{\underline{\Ext}}
\DeclareMathOperator{\Ker}{ker}

\DeclareMathOperator{\End}{End}
\DeclareMathOperator{\Sym}{Sym}

\DeclareMathOperator{\rank}{rank}

\DeclareMathOperator{\supp}{supp}

\DeclareMathOperator{\tr}{tr}

\DeclareMathOperator{\rk}{rk}

\DeclareMathOperator{\Var}{Var}

\DeclareMathOperator{\md}{mod}
\def\sat{\textrm{sat}}
\DeclareMathOperator{\sst}{ss}
\DeclareMathOperator{\out}{out}
\DeclareMathOperator{\gensurj}{gen-surj}
\DeclareMathOperator{\fin}{fin}

\DeclareMathOperator{\im}{Im}
\DeclareMathOperator{\coeff}{Coeff}

\def\gr{\mathrm{gr}}

\def\univ{\textrm{\tiny univ}}

\def\stable{{\text{stable}}}

\def\1halb{\frac{1}{2}}
\def\tto{\twoheadrightarrow}
\def\pprime{{\prime\prime}}

\def\sxymat{\xymatrix@C=1.5ex@R=0.8ex}
\def\grp{$\xymatrix{ R\times_{X}R  \ar[r]^-{\mu} & R \ar@<1ex>[r]^-{s}\ar@<-1ex>[r]_-{t} & X}$}
\def\dar{\ar@<-0.5ex>[r]\ar@<0.5ex>[r]}
\def\tar{\ar[r]\ar@<1ex>[r]\ar@<-1ex>[r]}
\newcommand{\dmap}[2]{\ar@<-0.5ex>[r]_-{#2}\ar@<0.5ex>[r]^-{#1}}
\newcommand{\dotarrow}[2]{\xymatrix{{#1}\ar@{..>}[r]&{#2}}}
\def\cart{\ar@{}[dr]|{\square}}

\def\cA{\mathcal{A}}

\def\cE{\mathcal{E}}
\def\cF{\mathcal{F}}

\def\cK{\mathcal{K}}
\def\cL{\ensuremath{\mathcal{L}}}
\def\cM{\mathcal{M}}

\def\cO{\mathcal{O}}

\def\cQ{\mathcal{Q}}

\def\cT{\mathcal{T}}

\def\cX{\mathcal{X}}

\def\bA{{\mathbb A}}

\def\bC{{\mathbb C}}

\def\bG{{\mathbb G}}
\def\bH{{\mathbb H}}

\def\bL{{\mathbb L}}

\def\bN{{\mathbb N}}

\def\bP{{\mathbb P}}
\def\bQ{{\mathbb Q}}
\def\bR{{\mathbf R}}

\def\bZ{{\mathbb Z}}

\def\K0hat{\widehat{K}_0(\Var)}
\def\ocM{\cM}
\def\rcM{\mathring{\cM}}

\begin{document}
%\input{deckbl}
%\title[Chains and Higgs bundles]{On the motives of moduli of chains and Higgs~bundles}
\title[Moduli of chains and Higgs bundles]{On the motives of moduli of chains and Higgs~bundles}
%\author[L. Alvarez-Consul]{Luis Alvarez-Consul}
\author[O. Garcia-Prada]{Oscar Garcia-Prada}
	\address{Instituto de Ciencias Matem\'aticas CSIC-UAM-UC3M-UCM,
	calle Nicol\'as Cabrera, 15, Campus de Cantoblanco, 28049 Madrid, Spain}
	\email{oscar.garcia-prada@icmat.es}
\author[J. Heinloth]{Jochen Heinloth}
	\address{University of Amsterdam, Korteweg-de Vries Institute for Mathematics, Science Park 904, 1098 XH Amsterdam,
The Netherlands}
	\email{J.Heinloth@uva.nl}
\author[A. Schmitt]{Alexander Schmitt}
	\address{Freie Universit\"at Berlin,
	Institut f\"ur Mathematik,
	Arnimallee 3,
	D-14195 Berlin,
	Germany}
	\email{alexander.schmitt@fu-berlin.de}
%\dedicatory{Version  \today}
%\keywords{}
%\subjclass{}
\begin{abstract}
We take another approach to Hitchin's strategy of computing the cohomology of moduli spaces of Higgs bundles by localization with respect to the circle-action. Our computation is done in the dimensional completion of the Grothendieck ring of varieties and starts by describing the classes of moduli stacks of chains rather than their coarse moduli spaces. 

As an application we show that the $n$-torsion of the Jacobian acts trivially on the middle dimensional cohomology of the moduli space of twisted $\SL_n$-Higgs-bundles of degree coprime to $n$ and we give an explicit formula for the motive of the moduli space of Higgs bundles of rank 4 and odd degree. This provides new evidence for a conjecture of Hausel and Rodr\'{\i}guez-Villegas. Along the way we find explicit recursion formulas for the motives of several types of moduli spaces of stable chains.
\end{abstract}
\maketitle

In this article we take another approach to implement Hitchin's strategy \cite[\S 7]{Hitchin} of computing the cohomology of moduli space $M_n^d$ of stable Higgs bundles of rank $n$ and degree $d$ on a curve $C$ by localization with respect to the circle-action, assuming that $n$ and $d$ are coprime.

In order to handle different possible choices of cohomology theories
in a uniform way, we want to compute the classes of these spaces in the (dimensional completion) of the Grothendiek ring of varieties $\K0hat$ (see Section \ref{K0intro}) and give the result in a form that allows us --- in case that $C$ is defined over the complex numbers --- to read off the Hodge polynomials of the spaces without further effort. 

Before explaining our strategy let us give an overview of the main results we obtain as an application of our approach.
We show (Theorem \ref{TrivialAction} in \S \ref{Application}) that in the case of coprime rank and degree the $n$-torsion points of the Jacobian of $C$ act trivially on the middle-dimensional cohomology of the space of twisted $\SL_n$-Higgs bundles. This answers a question of T.\ Hausel which was motivated by the conjecture of \cite{HT1}. 

Also we give an explicit formula (Theorem \ref{Higgs4} in \S \ref{HiggsRk4}) for the motive of the space of semistable Higgs bundles of rank 4 and odd degree. We implemented this formula in {\tt Maple} and checked that for genus $\leq 21$ the result confirms the conjecture of Hausel and Rodr\'{\i}guez-Villegas (\cite{HRV}) on the Poincar\'e polynomial of this space.

Along the way we obtain recursive formulas for the motives of spaces of $\alpha$-semistable chains of several types. %, e.g., chains of rank $(n,1,\dots,1)$ and in case of large stability parameters also for chains of rank $(n,\dots,n)$.
In particular, one could apply our method to give a recursive description of the motive of the space of stable $U(n,1)$-Higgs bundles. This application will be given elsewhere.

Let us now explain our strategy to obtain these results. It is known --- we will recall this in Section \ref{recollection} --- that the class $[M_n^d]\in \K0hat$ can be expressed in terms of classes of moduli spaces of stable chains of vector bundles on $C$. Here a chain is simply a collection $(\cE_0,\dots,\cE_r)$ of vector bundles together with morphisms $\phi_i\colon \cE_i\to \cE_{i-1}$ for $i=1,\dots r$. For chains one usually considers a stability criterion depending on parameters $\alpha=(\alpha_0,\dots,\alpha_r)$ and the notion of stability for Higgs bundles corresponds to a particular choice of $\alpha$.

Our strategy to compute the classes of the moduli spaces of chains of vector bundles is somewhat similar to the computation of the number of points of the moduli space of stable vector bundles given by Harder and Narasimhan \cite{HN} and the computation of the cohomology of this space over the complex numbers given by Atiyah and Bott \cite{AB}: In a first step we want to compute the motive of the whole stack of chains. The second step is then to study the Harder--Narasimhan stratification of the space of chains, i.e., the stratification according to the different types of canonical destabilizing subchains. As in the case of vector bundles, the Harder--Narasimhan strata are fibered over spaces of semistable chains of lower rank, for which we know the motive by induction.

In order to deduce the motive of the strata we need to describe these fibrations. Quite surprisingly the fibrations turn out to be smooth whenever the stability parameter $\alpha$ is larger or equal to the stability parameter needed for the application to moduli of Higgs-bundles (Lemma \ref{HN-Ext-Rechnung}). For other stability parameters this property would fail in general and this may indicate why the strategy to compute the cohomology by variation of the stability parameter turned out to be so difficult.

For the first step of our strategy, we make use of a result of Behrend and Dhillon \cite{BehrendDhillon}. They observed that the calculation of the cohomology of the stack of vector bundles on a curves given in \cite{BGL} can be interpreted in $\K0hat$. For stacks of arbitrary chains of vector bundles we define a stratification into pieces that we can compute explicitly in terms of the classes of moduli stacks of vector bundles (Proposition \ref{ChainStrata} and Corollary \ref{ChainStrata2}). The formulas turn out to be very simple. 

Let us immediately remark that to apply our programme in general, one has to overcome the problem, that although we can compute the pieces of the stratification, the summation over all strata does not always converge in $\K0hat$. 
If the summation does converge the above approach immediately gives a recursion formula for the class of the moduli stack of semistable chains and this is how we find our recursive formulas.

Further, in the cases needed to compute the cohomology of the space of Higgs-bundles of rank $4$ we show how to overcome the convergence problem by a truncation procedure. In particular we find a formula for the class of the stack of semistable chains of rank $(2,2)$ (Proposition \ref{m22}). For this example previous methods failed to compute the cohomology.

Let us briefly review the structure of the article. In Section \ref{K0intro} we recall the definition of the Grothendieck ring of varieties and a variant that contains the classes of algebraic stacks with affine stabilizer groups. We recall basic results on motivic zeta functions of varieties and the class of the moduli space of bundles on a curve that was calculated by Behrend and Dhillon. We end the section by explaining how to read off results on mixed Hodge polynomials from our formulas.

In Section \ref{HiggsIntro} we collect some known results on Higgs bundles and in particular explain how the computation of the class of the moduli space of stable Higgs bundles reduces to computations for moduli spaces of (holomorphic) chains.

After these preliminary sections we introduce (in Section \ref{BasicTools}) two general ingredients that will be the key to our computations. First we recall how to compute the class of stacks classifying of extensions of various types of objects and second, we give the class of the space of modifications of a family of vector bundles.

The core of the article is then contained in Sections \ref{ChainsGeneral} to \ref{ChainsExamples}. In Section \ref{ChainsGeneral}  we describe our general strategy to compute the class of moduli spaces of chains. In Section \ref{Application} we give the application to the middle-dimensional cohomology of the space of twisted $\SL_n$-Higgs bundles. 

In Section \ref{ChainsExamples} we apply our general strategy to do
explicit calculations. We first consider cases where our general
strategy carries through without additional effort to produce
recursive formulas for moduli spaces of chains.
We then show how in all chain types needed in order to deduce the class of the moduli space of Higgs bundles of rank $4$ one can solve the convergence problem mentioned above. 
The computation of the class of $[M_4^1]$ is then given in Section \ref{HiggsRk4}. Here we restrict ourselves to Higgs bundles in order to reduce the number of parameters involved, but the same arguments can also be used to obtain similar formulas for moduli spaces of bundles with a Higgs-field taking values in a line bundle $\cL$ with $\deg(\cL)>2g-2$.

For completeness we have included an appendix containing a quick
computation of the classes of $[M_2^1]$ and $[M_3^1]$. The class of
$[M_2^1]$ is implicitly already contained in Hitchin's original
article \cite{Hitchin}. The Poincar\'e polynomial of $[M_3^1]$ has
been computed by Gothen \cite{Gothen} and probably his argument could
also be refined to compute the motive, however the corresponding
Hodge-polynomial does not  seem to be available in the literature.

\noindent{\bf Acknowledgments:} J.H.\ would like to thank T.\ Hausel for sharing his insights and stimulating discussions and S. Mozgovoy for helpful remarks and corrections. O.G.P.\ was partially supported
by Ministerio de Ciencia e Innovaci\'on (Spain) Grant MTM2007-67623. J.H.\ was partially supported by the GQT-cluster of NWO and NFS grant No. DMS-0635607. A.S.\ was supported by SFB 647 ``Space-Time-Matter''. We thank the Newton Institute for ideal working conditions that allowed us to finish this work.

\noindent{\bf Notations:}
Throughout the article we will fix a smooth projective, geometrically connected curve $C$, defined over a field $k$. We will furthermore assume that there exists a line bundle of degree $1$ on $C$.

For a vector $\un{n}=(n_i)\in \bN^k$ we  denote $|n|:=\sum_{i=1}^k n_i$. 

%Given real numbers $a,b\in \bR$ we use the short hand $\sum_{k>a}^{<b}$ for the summation over all integers $k$ in the open interval $(a,b)$. For the summation over integers in the half-open interval $(a,b]$ we will write $\sum_{k>a}^{b}$. 

\tableofcontents

\section{Recollection on $\widehat{K_0(\Var)}$ and classes of algebraic stacks}\label{K0intro}

To explain the basic setup in which we do our calculations we need to
recall the definition of the motive of a local quotient stack from
the article \cite{BehrendDhillon} by Behrend and
Dhillon.\footnote{Other authors have considered similar definitions, see Joyce \cite{Joyce}, To\"en \cite{Toen} and Ekedahl \cite{Ekedahl}.} At the end of this section we collect the formulas needed to read off the mixed Hodge polynomials from our formulas. 

Our main reason for using motives rather than Hodge polynomials is that our computations only make use of geometric decompositions of the moduli spaces and thus it is natural to write our formulas in terms that reflect the underlying geometry. 

\subsection{The ring $\widehat{K_0}(\Var)$}
We denote by $K_0(\Var_k)$ the Grothendieck ring of varieties over $k$, i.e., it is the free abelian group generated by isomorphism classes of (quasi-projective) varieties $[X]$ subject to the relation $[X]=[X\setminus Z] +[Z]$ whenever $Z\subset X$ is a closed subvariety. (See \cite[Section 2]{LooijengaMotivic} for background on this ring.)

One denotes the class of the affine line by $\bL:=[\bA^1]\in K_0(\Var_k)$. Also for any quasi-projective variety $X$ we denote its symmetric powers by $\Sym^i(X)=X^{(i)}$. As observed in \cite{LarsenLunts} this extend this to classes in $K_0(\Var_k)$, by 
$$\Sym^n([X]+[Y]):= \sum_{i+j=n} [\Sym^i(X)][\Sym^j(Y)].$$

In $K_0(\Var_k)[\frac{1}{\bL}]$ we have the filtration defined by the subgroups generated by classes $[X]\bL^{-m}$ with $\dim(X)-m\leq n$ for $n\in \bN$ fixed. The completion of $K_0(\Var_K)$ according to this filtration is called the dimensional completion of $K_0(\Var_k)[\frac{1}{\bL}]$, we denote it by $\K0hat$.

Before defining  the motive of an algebraic stack with affine stabilizer groups, 
observe that $[\GL_n]=\prod_{k=0}^{n-1} (\bL^n-\bL^{k})$, by the usual argument that the first column of an invertible matrix is an arbitrary element of $\bA^n-\{0\}$, the second then gives a factor $\bA^n-\bA^1$ and so on. 
Since in $\K0hat$ we have $(\bL^{n}-1)^{-1}=\bL^{-n} \sum_{k=0}^\infty \bL^{-kn}$, we see that $[\GL_n]$ is an invertible element in $\K0hat$. 

Now suppose a stack $\cM$ is a quotient stack defined by an action of
$\GL_n$ on a scheme $X$, i.e., $\cM\cong [X/\GL_n]$. (Unfortunately
the standard notation for quotient stacks uses the same type of
brackets $[\;]$ that are used for classes in $\K0hat$.)  Then Behrend and Dhillon  define its class as
$$[\cM]:= \frac{[X]}{[\GL_n]}\in \K0hat.$$

It turns out that this definition does not depend on the choice of a presentation of $\cM=[X/\GL_n]$, because $\GL_n$-bundles are locally trivial in the Zariski topology.

In particular, for a quotient by an affine group $G$ one  can choose a faithful representation $G\to \GL_n$ and write any quotient as $[X/G]=[X\times^G \GL_n/\GL_n]$. Similarly one can then define the class of a stack which is stratified by locally closed substacks which are quotient stacks. This also makes sense for stacks  which are only locally of finite type, as long as they possess a stratification $\cM=\bigcup \cM_\alpha$, where for any $n\in \bZ$ only finitely many $\cM_\alpha$ are of dimension $\geq n$.

All stacks occurring in this article will admit a stratification into locally closed substacks that are of the form $[X/\GL_n]$. This follows for example from a theorem of Kresch \cite[Proposition 3.5.9]{Kresch}, by which it suffices to check that the stabilizer groups of all objects are affine.

\begin{example}[Behrend-Dhillon \cite{BehrendDhillon}]
Using the argument of \cite{BGL}, Behrend and Dhillon calculate the motive of the stack $\Bun_n^d$ classifying vector bundles of rank $n$ and degree $d$ on a smooth projective curve $C$: Denote by $$Z(C,t):=\sum_{k\geq 0} [C^{(k)}]t^k$$ the zeta function of $C$ and denote by $\Pic^0$ the Jacobian of $C$. Then Behrend and Dhillon \cite{BehrendDhillon} show:
$$ [\Bun_n^d] = \bL^{(n^2-1)(g-1)} \frac{[\Pic^0]}{\bL -1}\prod_{k=2}^{n} Z(C,\bL^{-k}).$$
\end{example}

\begin{notation}
We will often drop the degree and denote by $\Pic$ the Jacobian of $C$ and by $\uPic$ the stack of line bundles of degree $0$, so that we have $\uPic\cong \Pic \times B\bG_m$ and therefore $[\uPic]=\frac{[\Pic]}{\bL-1}$. 
\end{notation}

\subsection{Zeta functions and their relation with Hodge polynomials}\label{sym}
To compare our formulas with more classical formulas, we need to recall several facts on zeta functions from \cite{Kapranov} and \cite{FranziskaZeta}.
\begin{enumerate}
\item For any variety $X$ its zeta function is the formal power series $Z(X,t)=\sum [X^{(n)}] t^n$. The relation defining $\K0hat$ implies that the zeta function is multiplicative: for $Y\subset X$ closed we have  $Z(X,t)=Z(Y,t)Z(X\setminus Y,t)$. 
\item For the curve $C$ define $h^1(C):= [C]-1-\bL$ and set\footnote{In terms of cohomology $\Sym^i$ becomes the graded symmetric power, i.e., $\Sym^i h^1(C)$ corresponds to the exterior power of the first cohomology group.} $$P(t):= \sum_{i=0}^{2g} \Sym^i h^1(C) t^i.$$ Then 
	\begin{align*}
	 Z(C,t) &= \frac{P(t)}{(1-t)(1-\bL t)},\\
	 [\Pic] &= \sum_{i=0}^{2g} \Sym^i h^1(C) = P(1).\\
	\end{align*}
Thus we get $[\Bun_2^d] = \bL^{3g} \frac{P(1)P(\bL^{-2})}{(\bL-1)^2(\bL^2-1)}.$ Using the functional equation \cite[Section 3]{FranziskaZeta} $$P\left(\frac{1}{t\bL}\right)=\bL^{-g}t^{-2g} P(t)$$ for $t=\bL$ this simplifies to:
$$[\Bun_2^d] = \frac{P(1)P(\bL)}{(\bL-1)^2(\bL^2-1)}.$$

\item For any variety $X$ and $N,M\in \bZ$ we have
  \begin{align*}
   \sum_{k=0}^M X^{(k)} \bL^{Nk} &= \coeff_{t^0} \frac{Z(X,t)\bL^{NM}}{(1-t\bL^{-N})t^{M}}.
  \end{align*}
  This follows simply by expanding $\frac{1}{1-t\bL^{-N}}$ as a geometric series.
\end{enumerate}
Since we will need it later, let us give a simple application of these facts:
\begin{example}\label{SymCPn} 
The class $[(C\times \bP^{n-1})^{(l)}]\in \K0hat$ is given by:
$$[(C\times \bP^{n-1})^{(l)}] = \coeff_{t^l} \frac{\prod_{i=0}^{n-1} P(\bL^i t)}{\prod_{i=0}^{n-1} ((1-\bL^i t)(1-\bL^{i+1} t))}.$$
\end{example}
\begin{proof}
Since $[\bP^{n-1}]=\sum_{i=0}^{n-1} [\bA^i]$ we know from (1) that $Z(C\times \bP^{n-1},t)= \prod_i Z(C\times \bA^i,t)=\prod_i Z(C,\bL^i t)$. Together with (2) and (3) this implies the claimed formula.
\end{proof}

The preceding formulas will allow us to read off the compactly supported Hodge-polynomial of the moduli spaces we study from the their classes in $\K0hat$. 
% This is usually called the E-polynomial. The examples of moduli spaces of semi-stable chains will be projective, so that cohomology is equal to cohomology with compact supports. For the moduli space of Higgs-bundles this is not true. However in the situation we consider the cohomology is known to be pure and to satisfy Poincar\'e duality, so that the $E$-polynomial still determines the cohomology of the space.
First note that the $E$-polynomial (see, e.g.,\cite[\S 2]{HRV}) can be viewed as a map 
$E\colon K_0(\Var) \to \bZ[u,v],$
because for any closed subvariety $Y \subset X$ the long exact sequence for cohomology with compact supports implies that $E(X)=E(X\setminus Y)+ E(Y)$. 

Since $E(\bL)=uv$, this map extends to a map $$E\colon \K0hat \to \bZ[u,v]\bigg[\bigg[\frac{1}{uv}\bigg]\bigg],$$ taking values in Laurent-series in $(uv)^{-1}$.

\begin{example}\label{Hodge}
For the polynomial $P(t)=\sum_{i=0}^{2g} \Sym^i h^1(C) t^i$ defined above the description of the cohomology of symmetric products due to Macdonald \cite{MacdonaldSym} shows:
$$E(P(t))=(1-tu)^g(1-tv)^g \textrm{ and } E(Z(C,t))=\frac{(1-tu)^g(1-tv)^g}{(1-t)(1-tuv)}.$$
\end{example}

Our formulas will be given in terms of $P(t), [C^{(i)}]$ and $\bL$, so the above formulas will suffice to read off the E-polynomial of the moduli spaces from their class in $\K0hat$. If the cohomology of a variety $X$ is pure, e.g., if $X$ smooth and projective, then the E-polynomial determines the Hodge-polynomial by the formula
$$ H(X,u,v,t)=(uvt^2)^{\dim(X)}E\bigg(-\frac{1}{ut},-\frac{1}{vt}\bigg).$$ 
In this case the Poincar\'e polynomial of $X$ is given by
$$P(X,t)=t^{2\dim(X)}E\bigg(-\frac{1}{t},-\frac{1}{t}\bigg).$$

\section{Recollection on Higgs bundles}\label{recollection}\label{HiggsIntro}

In this section we collect the basic definitions on moduli spaces of Higgs bundles as well as how their topology is determined by the topology of moduli spaces of chains. For the convenience of the reader we briefly sketch the main arguments, more details can be found in \cite[Section 2 and 9]{HT1}.

A {\em Higgs bundle} is a pair $(\cE,\theta\colon \cE \to \cE \tensor
\Omega_C)$, where $\cE$ is a vector bundle on $C$, $\theta$ is an
$\cO_C$-linear map and $\Omega_C$ is the sheaf of differentials on
$C$. We will denote by $\cM_{n}^{d}$ the moduli stack of Higgs bundles
of rank $n$ and degree $d$ on $C$. In this stack we can consider the
open substack $\cM_{n}^{d,\sst}$ of semistable Higgs bundles, i.e.,
those $(\cE,\theta)$ such that for all subsheaves $\cF\subset \cE$
with $\theta(\cF)\subset \cF\tensor \Omega_C$ we have $\mu(\cF)\leq
\mu(\cE)$, where, as usual $\mu(\cF)$ is the ratio between the degree
and the rank of $\cF$. A Higgs bundle is called {\em stable} if this last inequality is a strict inequality  for all proper $(\cF,\theta_{|\cF})\subsetneq (\cE,\theta)$. Stability is also an open condition, so stable Higgs bundles define an open substack $\cM_{n}^{d,\stable} \subset \cM_{n}^{d,\sst}$. 

The stack of stable Higgs bundles turns out to be smooth. This follows from deformation theory (see Nitsure \cite{Nitsure} or Biswas and Ramanan \cite{BiswasRamanan}), showing that the first order infinitesimal deformations of a Higgs bundle $(\cE,\theta)$ are given by the cohomology of the complex $(\End(\cE) \to \End(\cE)\tensor \Omega_C)$. For a stable Higgs bundle, the only automorphisms are scalar automorphisms, so that $H^0$ of this complex is $1$-dimensional, and by Serre duality the same holds for $H^2$. Also the Euler-charactersitic of $H^*(C,\End(\cE) \to \End(\cE)\tensor \Omega_C)$ is $-2n^2(g-1).$

Therefore $\cM_{n}^{d,\stable}$ is a smooth stack of dimension $2n^2(g-1)+1$ and it is a $\bG_m$-gerbe over its coarse moduli space $M_{n}^{d}$, which is therefore smooth of dimension $2n^2(g-1)+2$.

Finally we have to recall the Hitchin map $$f\colon \cM_{n}^{d} \to \cA:=\bigoplus_{i=1}^n H^0(C,\Omega_C^i),$$ given by $f(\cE,\theta)=(\tr \wedge^i \theta) \in \bigoplus_{i=1}^n H^0(C,\Omega_C^i)$, i.e., $f$ maps a Higgs bundle $(\cE,\theta)$ to the coefficients of the characteristic polynomial of $\theta$. If $(n,d)=1$ Nitsure also proved \cite{Nitsure} that the induced map on the coarse moduli space $M_{n}^{d} \to \cA$ is proper. 

The moduli space of Higgs bundles has an action of $\bG_m$, given by multiplication of scalars on the Higgs field $\theta$.  The Hitchin map $f$ becomes equivariant with respect to this action, if we let $\bG_m$ act by the character $\lambda \mapsto \lambda ^i$ on the subspace $H^0(C,\Omega_C^i)\subset \cA$. 

We collect the known properties of the  Bia\l ynicki-Birula decomposition with respect to this action in the following proposition (which has been observed in \cite[Section 9]{HT1}):
\begin{proposition}\label{Stratification}
Let $n,d$ be a fixed pair of coprime, positive integers.
\begin{enumerate}
\item The fixed point scheme $(M_{n}^{d})^{\bG_m}$ of the $\bG_m$ action on $M_{n}^{d}$ is a disjoint union of connected, smooth schemes $F_i$ contained in the special fiber $f^{-1}(0)$ of the Hitchin map.

\item There are $\bG_m$-subvarieties $F_i^{+},F_i^-\subset M_{n}^{d}$  such that $F_i$ is a closed subscheme of $F_i^{\pm}$ and $F_i^{\pm}$ is a Zariski-locally trivial fibration over $F_i$, with fibres isomorphic to affine spaces.
For any $x\in F^i$ we have $T_x(F_i^\pm)=T_x(M)^0\oplus T_x(M)^{\pm}$ where $T_x(M)^{0,+,-}$ are the weight spaces of the tangent space at $x$ with weights repectively 0, positive or negative.

\item We have $M_{n}^{d}=\bigcup F_i^+$ and $(f^{-1}(0))=\bigcup F_i^-$ and for all points $x\in (M_{n}^{d})^{\bG_m}$ we have $\dim (T_x(M)^+)=\1halb \dim(M_{n}^{d})$. In particular the closure of the $F_i^-$ in $M_{n}^{d}$ are the irreducible components of $f^{-1}(0)$.
\end{enumerate}
\end{proposition}
\begin{proof}
The first part of the lemma is a direct consequence of the Bia\l ynicki-Birula  decomposition theorem (\cite[Theorem 4.1]{BialBirula} for algebraically closed $k$ and \cite[Theorem 5.8]{Hesselink} in general). The varieties $F_i^+$ consist of those points such that $\lim_{t\to 0} t.x\in F_i$ and $F_i^-$ consists of the points such that $\lim_{t\to \infty} t.x\in F_i$. 

Next we use that according to Nitsure \cite{Nitsure} the Hitchin map is proper. Since we know that $\bG_m$ acts with positive weights on $\cA$ this implies that the $F_i$ have to be contained in the special fibre of the fibration. Moreover, every point of $M_{n}^{d}$ is contained in exactly one $F_i^+$. Similarly, every point in the special fibre has to be contained in some $F_i^-$ and outside of the special fiber the limit $\lim_{t\to 0} t^{-1}.x$ does not exist.
 
As observed by Hausel and Thaddeus \cite{HT1} the claim on the dimension of $T_x(M)^+$ is a consequence of Laumon's result \cite{LaumonCone} that the downward flow is Lagrangian.
In particular, this implies that all the $F_i^-$ are smooth of dimension $\1halb \dim M_{n}^{d}$.
\end{proof}
This implies that the class of $[M_{n}^{d}]\in \K0hat$ can be computed in a very simple way from the classes of the $F_i$. This was observed in \cite[Proposition 9.1]{HT1} in terms of E-polynomials:
\begin{corollary}[Hausel-Thaddeus \cite{HT1}]\label{HT}
Write $N:= \1halb \dim M_{n}^{d}=n^2(g-1)+1$ then we have  $$ [M_{n}^{d}] = \bL^N \sum_i [F_i]\in \K0hat.$$
\end{corollary}

In order to make use of this result we need to recall the modular description of the fixed point strata of Lemma \ref{Stratification}, due to Hitchin and Simpson (see \cite[Lemma 9.2]{HT1}): If $(\cE,\theta)$ is a fixed point of the $\bG_m$ action, $(\cE,\theta)\cong (\cE,\lambda \theta)$ for all $\lambda\in \bG_m$. Then either $\theta=0$ or if $\theta\neq 0$ then this implies that the automorphism group of $\cE$ contains a copy of $\bG_m$. This implies that $\cE = \bigoplus \cE_i$ decomposes as a direct sum of weight spaces for this action and $\theta: \cE_i \to \cE_{i-1}\tensor \Omega_C$. This implies that each $F_i$ is a moduli space of bundles $\cE_i$ together with maps $\phi_i\colon \cE_i \to \cE_{i-1} \tensor \Omega$, such that the corresponding Higgs-bundle $(\bigoplus \cE_i, \bigoplus \phi_i)$ is stable. These are called {\em moduli spaces of stable chains} and we will recall their properties in more detail in Section \ref{ChainsGeneral}. 
The main part of this paper  will be devoted to the computation of the classes of these moduli spaces. 

\section{The basic ingredients of our calculations: extensions and modifications}\label{BasicTools}

Our calculations rely on two basic results, which we would like to explain in this section. Firstly, an observation already contained in SGA 4 \cite[Expos\'e XVIII 1.4]{SGA4} allows to compute the class of spaces of extensions of bundles, or more generally of extensions of chains whenever the extension problem is unobstructed. Secondly, we can compute the class of the stack of modifications of vector bundles. As any morphism of bundles $\cE_1\map{\phi} \cE_0$ can be viewed as an extension of vector bundles $\ker(\phi) \to \cE_1 \to \im(\phi)$, followed by a modification $\im(\phi) \to \im(\phi)^\sat$ (here we denoted by $\im(\phi)^\sat$ the saturation of the image) and another extension of vector bundles $\im(\phi)^\sat \to \cE_0 \to \cE_0/\im(\phi)^\sat$ the above ingredients will allow us to describe moduli spaces of chains.

\subsection{Stacks classifying extensions of objects}

In the following we will often consider stacks parameterizing
extensions of bundles or chains. In order to compute their motives we will apply the following result, which appeared in SGA4 \cite[Expos\'e XVIII, Proposition 1.4.15]{SGA4}, the statement as below can also be found in \cite{Diss}:

\begin{proposition}\label{VectorBundleStack}
Let $\cX$ be an algebraic stack, $\cE_0,\cE_1$ vector bundles on $\cX$ and $\cE_0\map{d} \cE_1$ a morphism. Viewing  $\cE_0$ as an affine group scheme over $\cX$ acting on $\cE_1$ via $d$ we form the quotient stack $[\cE_1/\cE_0]$. 
Then for any affine scheme $T\map{t}\cX$ over $\cX$ the category $[\cE_1/\cE_2](T)$ is equivalent to:
$$[\cE_1/\cE_0](T) \cong \left\langle \textrm{Objects}= H^1\big(T,t^*(\cE_0\to\cE_1)\big), \text{ Morphisms } H^0\big(T,t^*(\cE_0\to\cE_1)\big)\right\rangle.$$
\end{proposition}

Stacks of the form $[\cE_1/\cE_0]$ as occurring in the above proposition are called {\em vector bundle stacks}. The above proposition shows that quasi-isomorphic complexes $\cE_\bullet$ define equivalent stacks. 

As an illustration of how we will apply the above proposition let us recall a well-known application.
Denote by $\Coh_{n}^d$ the stack of coherent sheaves of rank $n$ and degree $d$ on $C$. Denote by $\uExt((n^\pprime,d^\pprime),(n^\prime,d^\prime))$ the stack classifying extensions $\cF^\prime \to \cF \to \cF^\pprime$ of coherent sheaves with $(\rk(\cF^\prime),\deg(\cF^\prime))=(n^\prime,d^\prime)$ and  $(\rk(\cF^\pprime),\deg(\cF^\pprime))=(n^\pprime,d^\pprime)$. 

\begin{corollary}\label{ext-stack}
The forgetful map $p\colon \uExt((n^\pprime,d^\pprime),(n^\prime,d^\prime))\to \Coh_{n^\prime}^{d^\prime}\times \Coh_{n^\pprime}^{d^\pprime}$ defines a vector bundle stack of relative dimension $n^\prime n^\pprime ((g-1)+\frac{d^\pprime}{n^\pprime} - \frac{d^\prime}{n^\prime})$. 
In particular $\uExt((n^\pprime,d^\pprime),(n^\prime,d^\prime))$ is smooth and irreducible.
\end{corollary}
\begin{proof}
Denote by $\cF^\prime_{univ}$ and $\cF^\pprime_{univ}$ the universal sheaves on $\Coh_{n^\prime}^{d^\prime}\times C$ and $\Coh_{n^\pprime}^{d^\pprime}\times C$.
Denote by $p_{ij}$ the projection on the $i,j$-th factor of $\Coh_{n^\prime}^{d^\prime}\times \Coh_{n^\pprime}^{d^\pprime} \times C$.

We claim that for any substack of finite type of $\Coh_{n^\prime}^{d^\prime}\times \Coh_{n^\pprime}^{d^\pprime}$
the complex $\bR p_{12,*} (\cHom(p_{23}^* \cF^\pprime_{\univ},p_{13}^* \cF^\prime_{\univ}))$ is quasi-isomorphic to a complex $[\cE_0\to \cE_1]$ where $\cE_i$ are vector bundles. This holds since for any bounded family of coherent sheaves $\cF=\cHom(p_{23}^* \cF^\pprime,p_{13}^* \cF^\prime)$ on $C$ there exists $d\gg 0$ such that for any closed point $p\in C$ the sheaves $\cF(d\cdot p)$ have no higher cohomology. But then for $p,q\in C$ the complex $p_{12,*}(\cF(d \cdot p) \oplus \cF( d \cdot q) \to \cF(d (p+q)))$ is quasi-isomorphic to $\bR p_{12,*} (\cHom(p_{23}^* \cF^\pprime,p_{13}^* \cF^\prime))$.
 
By Proposition \ref{VectorBundleStack} we know that $[\cE_1/\cE_0]$ is isomorphic to $\uExt((n^\pprime,d^\pprime),(n^\prime,d^\prime))$. The relative dimension of the morphism is $\rk \cE_1-\rk \cE_2$, which by the Riemann-Roch formula is $n^\prime n^\pprime (g-1)+(n^\prime d^\pprime - n^\pprime d^\prime)$.
\end{proof}

The classes of vector bundle stacks are easy to compute in $\K0hat$:
\begin{lemma}\label{FormulaVectorBundleStack}
Suppose that $\cX,\cE_0,\cE_1$ are as in Proposition \ref{VectorBundleStack}, that $\cX$ is a local quotient stack, which defines a class in $\K0hat$, and that $\cE_0,\cE_1$ are of constant rank. Then we have 
$$[\cE_1/\cE_0] = [\cX]\bL^{\rk(\cE_1)-\rk(\cE_0)} \in \K0hat.$$
\end{lemma}
\begin{proof}
First, the lemma holds in the case that $X$ is a scheme: Stratify $X=\bigcup X_i$ such that over each $X_i$ the map $\phi|_{X_i}$ is of constant rank and both $\cE_1$ and $\cE_0$ are trivial. By the proposition we can, over each $X_i$, replace the complex $\cE_0 \map{d} \cE_1$ by $\ker(d) \map{0} \cE_1/\im(d)$. In this case we have
$$[\cE_1/\cE_0|_{X_i}]=[X_i] \times [\bA^{\rk(\cE_1/\im(d))}] \times [\Spec(k)/\bG_a^{\rk(\ker(d))}] = [X_i]\bL^{\rk(\cE_1)-\rk(\cE_0)}.$$

To prove the lemma for a local quotient stack it suffices to consider the case that $\cX=[X/\GL_n]$ is a global quotient. Let $p\colon X \to [X/\GL_n]$ denote the canonical projection. Note that the automorphism groups of all objects of $[\cE_1/\cE_0]$ are affine, so that by Kresch's result \cite[Proposition 3.5.9]{Kresch} the stack $[\cE_1/\cE_0]$ is again a local quotient stack.  Write $[\cE_1/\cE_0]=\cup [F_i/\GL_{n_i}]$ for some schemes $F_i$ and $n_i\in \bN$. This decomposition induces a decomposition $$[p^*\cE_1/p^*\cE_0]\cong X\times_{[X/\GL_n]} [\cE_1/\cE_0] = \bigcup X\times_{[X/\GL_n]}[F_i/\GL_{n_i}].$$
Moreover, $F_i\times_{[X/\GL_n]} X \to F_i$ is a $\GL_n$ torsor and $$F_i\times_{[X/\GL_n]} X  = [p^*\cE_1/p^*\cE_0] \times_{[\cE_1/\cE_0]} F_i \to X\times_{[X/\GL_n]}[F_i/\GL_{n_i}]$$ is a $\GL_{n_i}$-torsor.
Therefore:
\begin{align*}
[\cE_1/\cE_0] &= \sum_i \frac{[F_i]}{[\GL_{n_i}]}= \sum_i \frac{[F_i][\GL_n]}{[\GL_{n_i}][\GL_n]}= \sum_i \frac{[F_i\times_{[X/\GL_n]} X]}{[\GL_{n_i}][\GL_n]}\\
&= \frac{1}{[\GL_n]}\sum_i \Big[X\times_{[X/\GL_n]}[F_i/\GL_{n_i}]\Big]= \frac{[p^*\cE_1/p^*\cE_0]}{[\GL_n]}\\
& = \frac{[X]\bL^{\rk(\cE_1)-\rk(\cE_0)}}{[\GL_n]}= [\cX]\bL^{\rk(\cE_1)-\rk(\cE_0)}.
\end{align*}
This proves the lemma.
\end{proof}

\begin{example}\label{Bun2}
As simplest case of our problem, let us compute the motive of the stack of stable bundles of rank $2$ and degree $d$.

The Harder--Narasimhan stratum $\Bun_2^{d,\mu=l}$ of those bundles having a subbundle of degree $l>\frac{d}{2}$ is the stack classifying extensions $\cL \to \cE \to \cQ$ with $\deg(\cL)=l,\deg(\cQ)=d-l$. Thus we can apply Corollary \ref{ext-stack} and Lemma \ref{FormulaVectorBundleStack} to find:
$$ [\Bun_2^{d,\mu=l}]=\bL^{g-1+d-2l}[\Bun_1^{l}][\Bun_1^{d-l}]=\bL^{g-1+d-2l}\frac{P(1)^2}{(\bL-1)^2},$$
where the last equality uses \ref{sym}.
Thus for the class of the stack of semi-stable bundles we find:
\begin{align*}
[\Bun_2^{d,\sst}] &=  [\Bun_2^d] - \sum_{l>\frac{d}{2}} [\Bun_2^{d,\mu=l}]\\
                &=  \frac{P(1)P(\bL)}{(\bL-1)^2(\bL^2-1)} - \sum_{l>\frac{d}{2}} \bL^{g-1+d-2l}\frac{P(1)^2}{(\bL-1)^2}\\
&= \frac{(P(1)(P(\bL)-\bL^{g-1 +  (d\md 2)}P(1))}{(\bL-1)^2(\bL^2-1)}.
\end{align*}											

If one applies $P_{\textrm{Hodge}}$ to this class in case $d$ is odd, this formula gives a polynomial divided by $P_{\textrm{Hodge}}(\bL-1)$, corresponding to the fact that the stack is a $\bG_m$-gerbe over the coarse moduli space.

We'd like to point out that the same method would also allow us to compute the class of the coarse moduli space of stable bundles (see \cite{delBano}) in case the degree is even, by further discarding the strictly semi-stable bundles. 
\end{example}

\begin{remark}\label{Bunn}
For general $n$ one can write $[\Bun_n^{d,\sst}]=[\Bun_n^d]-\bigcup \textrm{HN-Strata}$. Corollary \ref{ext-stack} then gives a recursive formula for $[\Bun_n^{d,\sst}]$. Namely the Harder--Narasimhan strata are indexed by partitions $n=n_1+\dots+n_s$ and $d=\sum d_i$ with ${d_i}/{n_i}>{d_{i+1}}/{n_{i+1}}$ for all $i$. The class of such a stratum is given by $$\bL^{\sum_{i<j} n_in_j(g-1) + \sum_{i<j} d_jn_i-d_in_j} \prod_{i=1}^s [\Bun_{n_i}^{d_i,\sst}].$$ 

This recursive formula has been solved by Zagier \cite{Zagier} and Laumon and Rapoport in \cite{LaumonRapoport}, who formulated the result in terms of the Poincar\'e series in their article and used cohomology instead of cohomology with compact supports. Their argument shows:  
\begin{align*}
&[\Bun_n^{d,\sst}]= \\ 
& \sum_{s=1}^n (-1)^{s-1}\sum_{n=n_1+\dots+n_s\atop n_i>0} \prod_{i=0}^s [\Bun_{n_i}^0] \bL^{\sum_{i<j} n_in_j (g-1)} \prod_{i=1}^{s-1} \frac{\bL^{(n_i+n_{i+1})\frac{(n_1+\dots+n_i)d \md n}{n}}}{\bL^{n_i+n_{i+1}}-1}.
\end{align*}
Details on the computation using cohomology with compact supports can be found in \cite{deRonde}.
\end{remark}

%%%%%%%%%%%%%
\subsection{The class of the stack of modifications of bundles}
%%%%%%%%%%%%%%

The second basic ingredient for our computation is the class of the stack of modifications.
We write $\Hecke_{n,d}^l$ the stack classifying $(\cE_1\map{\phi} \cE_0)$ such that $\cE_i$ are vector bundles of rank $n$, degree $\deg(\cE_0)=d,\deg(\cE_1)=d-l$ and $\rank(\phi)=n$. This is usually called the stack of Hecke modifications of length $l$. %We will sometimes write $\Hecke_{n,d}^l:=\cM(n,n)_{\rk=n}^d$.

Also for any family of vector bundles $\cE$ of rank $n$ parameterized by a scheme of finite type (or a stack of finite type with affine stabilizer groups) $T$ we will write $\Hecke(\cE/T)^l$ for the stack classifying modifications $\cE^\prime \subset \cE$ with $\cE/\cE^\prime$ a torsion sheaf of length $l$. 

The argument of \cite{BGL} implicitly contains the next result. We give a slightly
different argument, since we need to work in $\K0hat$ and we will need a result over a general base.
\begin{proposition}\label{Mnn-rkn}The class of the stack of Hecke modifications is:
$$[\Hecke_{n,d}^l]= [\Bun_n^{d}] \times [(C\times \bP^{n-1})^{(l)}] \in \widehat{K}_0(\Var).$$
Similarly $[\Hecke(\cE/T)^l] = [T] \times [(C\times \bP^{n-1})^{(l)}].$
\end{proposition}
\begin{proof}
Since $\Hecke_{n,d}^l$ parameterizes pairs $\cE_1\subset \cE_0$ we have a canonical morphism:
\begin{align*}
 {\gr}\colon \Hecke_{n,d}^l &\to \Bun_n^{d} \times C^{(l)}\\
 (\cE_1\to \cE_0) &\mapsto (\cE_0, \supp(\cE_0/\cE_1))
\end{align*}
We first compute the fibers of $\gr$. For a point $P\in C$ and $\cE\in\Bun_n^d$ we denote $F_{\cE,lP}:={\gr}^{-1}(\cE,lP)$.

{\em Claim:} $[{\gr}^{-1}(\cE,lP)]=[\Sym^l \bP^{n-1}]$.

Although this is probably known,  for the sake of completeness we will give an inductive proof. For $n=1$ the map $\gr$ is an isomorphism, so the claim is clear.

In general, choose a trivialization $\cE|_{\cO_{C,P}}\cong \cO_{C,P}^{\oplus n}$ and chose a local parameter $t$ at $P$ in order to obtain an isomorphism $\widehat{\cO}_P\cong k[[t]]$. In particular the first summand of $\cO_{C,P}^{\oplus n}$ defines a subbundle $\cL \to \cE$ and we can stratify the space of modifications $\cE_1\to \cE$ according to the length $k$ of the image of $\cL$ in $\cE/\cE_1$:
$$\xymatrix{
\cL(-kP) \ar[r]\ar[d] & \cL \ar[r]\ar[d] & \cO/\cO(-kP)\ar[d] \\
\cE_1 \ar[r]\ar[d] & \cE \ar[r]^q\ar[d] & \cT_l\ar[d] \\
\cE_1/\cL(-kP) \ar[r] & \cE/\cL \ar[r]^{q^\pprime} & \cT_{l-k}.
}$$
The space of all such extensions is fibered over the space of modifications of $\cE^\pprime$ of length $l-k$.
The fibers are defined by the extensions of $\cT_{l-k}$ by $\cO/\cO(-kP)$ together with a choice of a map $q$.

The extensions are classified by $\uExt(\cT_{l-k},\cO/\cO(-kP))$ which is a vector bundle stack of rank $0$ over the space parameterizing the torsion sheaves $\cT_l$.
The choices of $q$ form a torsor under $\Hom(\cE^\pprime,\cO/\cO(-kp))$ which is a vector space of dimension $k\cdot (n-1)$. Thus we find $[F_{\cE,lP}] = \sum_{k=0}^l [F_{\cE^\pprime,(l-k)P}] \bL^{k(n-1)}$.

On the other hand, by Remark \ref{sym} we have
\begin{align*}
[\Sym^l(\bP^{n-1})] &= \coeff_{t^l} \frac{1}{(1-t)(1-\bL t)\cdots(1-\bL^{n-1}t)} \\
                  &= \sum_{k=0}^l \coeff_{t^{l-k}} \frac{1}{(1-t)(1-\bL t)\cdots(1-\bL^{n-2}t)}\bL^{k(n-1)}\\
                  &= \sum_{k=0}^l \bL^{k(n-1)} [\Sym^{l-k} \bP^{n-2}].
\end{align*}
This proves the claim.

An arbitrary point of $C^{(l)}$ is an effective divisor $D=\sum m_i P_i$ with $\sum m_i=l$ and $P_i\neq P_j$ for $i\neq j$. So the partitions $l=\sum m_i$ with $m_1\geq \dots \geq m_r$ define a stratification $C^{(l)}_{\un{m}}$ of $C^{(l)}$. In order to compute the class $[\Hecke_{n,d}^l]$ we stratify this stack accordingly. 

To deduce the result for families $(\cE,D=\sum m_i P_i)\in \Bun_n^{d}\times C^{(l)}_{\un{m}}(T)$ parameterized by some scheme $T$ we note that Zariski locally over $T$ we may choose local parameters at $P_i$. Also after replacing $T$ by a $\GL_n^r$-bundle over $T$ we may assume that there exist trivializations of $\cE$ at the points $P_i$. These local trivializations define an isomorphism of the fiber ${\gr}^{-1}(T)$ and $T\times \prod_{i} F_{\cO^{\oplus n},m_iP_i}$.
Thus we find:
$$[{\gr}^{-1}(T)]=[T]\times \prod_i [\Sym^{m_i}(\bP^{n-1})].$$
Similarly if $[T/\GL_n]$ is a any substack of finite type of $\Bun_{n}^d\times C^{(d)}_{\un{m}}$ we can deduce the same formula for this substack, because the map $\gr$ is representable and $\gr^{-1}([T/\GL_n])=[\gr^{-1}(T)/\GL_n]$.

To conclude, observe that the fiber of the projection $p: (C\times \bP^{n-1})^{(l)} \to C^{(l)}$ over a point $x\in C^{(l)}_{\un{m}}$ is isomorphic to $\prod_i \Sym^{m_i}(\bP^{n-1})$, since stabilizer in $S_l$ of a preimage of $x$ in $C^{(l)}$ is isomorphic to $\prod S_{m_i}$. Since $\bP^{n-1}$ is stratified by affine spaces and the permutation action is linear on the strata this implies (see \cite[Lemma 4.4]{Goettsche}):$$[p^{-1}(C^{(l)}_{\un{m}})]= \prod_i [\Sym^{m_i}(\bP^{n-1})] \cdot [C^{(l)}_{\un{m}}]\in \K0hat.$$ 
 Thus the sum over all strata $C^{(l)}_{\un{m}}$ can be written as:
$$ [\Hecke_{n,d}^l]= [\Bun_n^{d}]\times [(C\times \bP^{n-1})^{(l)}].$$
This proves the proposition.
\end{proof}

%\begin{remark} 
%We can give a more explicit formula for the class $[(C\times \bP^{n-1})^{(l)}]\in \K0hat$ occurring in the preceding proposition.  Since $[\bP^{n-1}]=\sum_{i=0}^{n-1} [\bA^i]$ we have that $Z(C\times \bP^{n-1},t)= \prod_i Z(C\times \bA^i,t)=\prod_i Z(C,\bL^i t)$. Therefore we have:
%$$[(C\times \bP^{n-1})^{(l)}] = \coeff_{t^l} \frac{\prod_i P(\bL^i t)}{\prod_i ((1-\bL^i t)(1-\bL^{i+1} t))}.$$
%\end{remark}

%%%%%%%%%%%%%%%%%%%%%
\section{Moduli stacks of chains: general results}\label{Extensions}\label{ChainsGeneral}
%%%%%%%%%%%%%%%%%%%%%

After recalling some basic definitions on chains (\cite{AG},\cite{AGS}), we will prove in this section that the Harder--Narasimhan strata in the moduli stacks of chains are vector bundle stacks over moduli stacks of chains of smaller rank. Moreover we construct another stratification of the stack of chains such that we can compute the classes of the strata in $\K0hat$. These are the key results needed to implement our strategy for the computation of the motives of the spaces of semistable chains.
 %We end the section with an application of these general results to the cohomology of the space of Higgs bundles with fixed determinant.

%%%%%%%%%%%%
\subsection{Stability of chains and basic properties of the moduli stack}
%%%%%%%%%%%%

A {\em (holomorphic) chain} on $C$ is a collection $((\cE_i)_{i=0,\dots,r},(\phi_i)_{i=1,\dots,r})$, where $\cE_i$ are vector bundles on $C$ and $\phi_i\colon\cE_i\to\cE_{i-1}$ are morphisms of $\cO_C$-modules. We will often abbreviate $((\cE_i)_{i=0,\dots,r},(\phi_i)_{i=1,\dots,r})$ as $\cE_\bullet$.

The rank of a chain is defined as $\rk(\cE_\bullet)=(\rk(\cE_i))_{i=0,\dots,r}$ and the degree is defined as $\deg(\cE_\bullet):=(\deg(\cE_i))_{i=0,\dots,r}$. Given $\cE_\bullet=(\cE_i,\phi_i)$ we will denote the compositions of the $\phi_i$ by $\phi_{ij}:=\phi_i\circ\dots\circ\phi_j$.

We will denote by $\cM(\un{n})_{\un{d}}$ the moduli stack of chains of rank $\un{n}$ and degree $\un{d}$. Write $\rcM(\un{n})_{\un{d}}\subset \ocM(\un{n})_{\un{d}}$ for the open substack such that $\phi_i\neq 0$ whenever $n_in_{i-1}\neq 0$. To show that the stack $\ocM(\un{n})_{\un{d}}$ is an algebraic stack, locally of finite type, we only have to observe that the forgetful map $\cM(\un{n})_{\un{d}}\to \prod_{i=0}^r\Bun_{n_i}^{d_i}$ is representable.  This holds, because the fibers parameterize morphisms of sheaves. 

Given $(\alpha_i)_{i=0,\dots,r}\in \bR^{r+1}$ the {\em $\un{\alpha}$-slope} of a chain $\cE_\bullet$ is defined as
\begin{align*}
\mu(\cE_\bullet)=\mu_\alpha(\cE_\bullet) &:= \sum_{i=0}^r \frac{\rk(\cE_i)}{|\rk(\cE_\bullet)|}(\mu(\cE_i) + \alpha_i),
\end{align*}
where $\mu(\cE_i):={\deg(\cE_i)}/{\rk(\cE_i)}$ is the slope of the vector bundle $\cE_i$ and the summand is read as $0$ if $\rk(\cE_i)=0$. Note that the $\alpha$-slope is a convex-combination of $\mu(\cE_i)+\alpha_i$. Since the slope only depends on the numerical invariants for fixed rank and degree $\un{n},\un{d}$, we also write $\mu(\un{n},\un{d})$ for the corresponding $\alpha$-slope.

A chain is called {\em $\alpha$-(semi)-stable} if for all proper subchains $\cE_\bullet^\prime \subset \cE_\bullet$ we have 
$$\mu(\cE^\prime_\bullet) (\leq) \mu(\cE_\bullet),$$
where we use the standard notation $(\leq)$ to abbreviate that the inequality $\leq$ defines semistability, whereas for stability, we require strict inequality.

\begin{remark}\label{HiggsStability}
Any chain $\cE_\bullet$ defines a Higgs bundle $\cE:=\bigoplus_i \cE_i\tensor \Omega^{-r-i}$ with Higgs field $\cE \to \cE\tensor \Omega$ given by the sum of the $\phi_i$. 
This Higgs-bundle is (semi)-stable if and only if the chain $\cE_\bullet$ is $\alpha$-(semi)-stable for the parameter $\alpha=(0,2g-2,\dots,r(2g-2))$.
\end{remark}

Note that --- as in the case of vector bundles --- given an extension $\cE_\bullet^\prime \to \cE_\bullet \to \cE_\bullet^\pprime$ of chains of rank $\un{n}^\prime$ and $\un{n}^\pprime$ we have
\begin{align*}
\mu(\cE_i) &= \frac{n_i^\prime}{n_i} \mu(\cE_i^\prime) + \frac{n_i^\pprime}{n_i}\mu(\cE_i^\pprime),\\
\mu(\cE_\bullet) %&= \mu(\cE_\bullet^\prime \oplus \cE_\bullet^\pprime) \\
								 &= \frac{|\un{n}^\prime|}{|\un{n}|}\mu(\cE_\bullet^\prime) + \frac{|\un{n}^\pprime|}{|\un{n}|}\mu(\cE_\bullet^\pprime).
\end{align*}
So the slope of an extension is a convex combination of the slope of the constituents. As for stability of vector bundles, this property immediately implies the following properties of stability of chains:
\begin{lemma}\label{stab}
 \begin{enumerate}
  \item A chain $\cE_\bullet$ is semistable if and only if for any quotient  $\cE_\bullet \tto \cE^\pprime_\bullet$ we have $\mu(\cE_\bullet)\leq \mu(\cE_\bullet^\pprime)$.
  \item If $\cE_\bullet,\cF_\bullet$ are semistable with $\mu(\cE_\bullet)>\mu(\cF_\bullet)$ then $\Hom(\cE_\bullet,\cF_\bullet)=0$.
  \item For every chain $\cE_\bullet$ there is a canonical Harder--Narasimhan flag of subchains $0\subset \cE_\bullet^{(1)} \subset \dots \cE^{(h)}_\bullet= \cE_\bullet$, such that $\mu(\cE^{(1)}_\bullet)>\dots >\mu(\cE_\bullet^{(h)})$ and the subquotients $\cE_\bullet^{(i)}/\cE_\bullet^{(i-1)}$ are semistable.
 \end{enumerate}
\end{lemma}

We will denote by $\cM(\un{n})_{\un{d}}^{\alpha-\textrm{ss}}\subset \cM(\un{n})_{\un{d}}$ the substack of $\alpha$-semistable chains. The same argument as for vector bundles shows that this is an open substack and that its complement is the disjoint union of the Harder--Narasimhan strata, i.e., the constructible substacks of those $\cE_\bullet$ such that the $\cE_\bullet^{(i)}$ are of some fixed rank and degree. %We will often drop the fixed degree $\un{d}$ from the notation.

Note that for any $c\in \bR$ and $\alpha=(\alpha_i)_{0=1,\dots,r}$ and $\alpha+c:=(\alpha_i+c)_{i=0,\dots,r}$ define the same stability condition. 
Also if we denote by $\overline{\alpha}:=(-\alpha_{r-i})_{0=1,\dots,r}$ then we have $$\mu_{\overline{\alpha}}(\cE_\bullet^\vee)=-\mu_\alpha(\cE_\bullet).$$
\begin{lemma}\label{dual}
Dualizing gives an isomorphism $$\cM(n_0,\dots,n_r)_{(d_0,\dots,d_r)}^{\alpha-\sst} \cong \cM(n_r,\dots,n_0)_{(-d_r,\dots,-d_0)}^{\overline{\alpha}-\sst}.$$

In particular for $\alpha_i:= i\cdot (2g-2)$ we have
$$\cM(n_0,\dots,n_r)_{(d_0,\dots,d_r)}^{\alpha-\sst} \cong \cM(n_r,\dots,n_0)_{(-d_r,\dots,-d_0)}^{\alpha-\sst}.$$
\end{lemma}
\begin{proof}
The first claim is immediate from the equivalent characterization of stability  given in Lemma \ref{stab} (1). The second follows, because in this case $\overline{\alpha}+r(2g-2)=\alpha$.
\end{proof}

Given $\un{n},\un{d}$ a parameter $\alpha$ is called {\em critical} (for $\un{n},\un{d}$), if there exist strictly semi\-stable chains of rank $\un{n}$ and degree $\un{d}$. Otherwise $\alpha$ is called {\em non-critical}.

Let us call $\alpha$ to be {\em good} if for any chain $\cE_\bullet^\prime$ occurring as subquotient $\cE_\bullet^i/\cE_\bullet^{i-1}$ in the Harder--Narasimhan filtration of a chain of rank $\un{n}$ and degree $\un{d}$ the following holds:
\begin{enumerate}
\item Whenever $\cE_i^\prime$ and $\cE_{i-1}^\prime$ are non-zero then $\phi_i$ is non-zero, and
\item the set $\{0\leq i \leq r \mid \cE_i^\prime \neq 0\}$ is an interval in $\bZ$.
\end{enumerate}
Note that any chain $\cE_\bullet^\prime$ violating one of the above conditions can be written as $\cE_\bullet^\prime=\cE_\bullet^1\oplus \cE_\bullet^2$ such that $\cE_\bullet^{1,2}$ are both non-trivial and for all $i$ either with $\cE_i^1$ or $\cE_i^2$ is $0$.  

Thus we see that a parameter $\alpha$ is good if the $\alpha_i$ are linearly independent over $\bQ$, because in that case $\mu(\cE_\bullet^1)\neq \mu(\cE_\bullet^2)$. In particular this implies that if $\alpha$ is not critical, then there is a good, non critical $\alpha^\prime$ such that $\cM(\un{n})^{\alpha-\sst}_{\un{d}}=\cM(\un{n})^{\alpha^\prime-\sst}_{\un{d}}$.

%%%%%%%%%%%
\subsection{Extensions of chains and the classes of Harder--Narasimhan strata}
%%%%%%%%%%%

Note that one can embed the category of chains into an abelian category, by allowing the $\cE_i$ to be coherent sheaves instead of vector bundles. In this category one can then do homological algebra (see, e.g., \cite{GothenKing}).

In this section it will further be useful to consider more generally chains $\cE_\bullet =(\cE_i,\phi_i)$ with $i\in \bZ$ such that only finitely many $\cE_i$ are non-zero. We will extend any chain $\cE_\bullet=((\cE_i)_{i=0,\dots r},(\phi_i)_{i=1,\dots,r})$ by putting $\cE_i:=0$ for all $i<0$ and all $i>r$. Similarly we will allow stability parameters $\alpha=(\alpha_i)_{i\in \bZ}$. 

\begin{notation}
Given chains $\cE^\prime_\bullet,\cE^\pprime_\bullet$ we denote by $\Hom(\cE_\bullet^\pprime,\cE^\prime_\bullet)$ the group of homomorphisms of chains, by $\Ext^1(\cE^\pprime_\bullet,\cE^\prime_\bullet)$ the set of isomorphism classes of extensions $\cE^\prime_\bullet \to \cE_\bullet \to \cE_\bullet^\pprime$ and by  $\underline{\Ext}(\cE^\pprime_\bullet,\cE^\prime_\bullet)$ the stack of such extensions. In particular we have $$\underline{\Ext}(\cE^\pprime_\bullet,\cE^\prime_\bullet) = [\Ext^1(\cE^\pprime_\bullet,\cE^\prime_\bullet)/\Hom(\cE_\bullet^\pprime,\cE^\prime_\bullet)].$$

For chains $\cE_\bullet^1,\dots,\cE_\bullet^h$ we denote by $\uExt(\cE_\bullet^h,\dots,\cE_\bullet^1)$ the stack of iterated extensions, i.e., chains $\cE_\bullet$ together with a filtration $0=\cF_\bullet^0\subset \cF_\bullet^1 \subset \dots \subset \cF_\bullet^h=\cE$ and isomorphisms $\cF_\bullet^i/\cF_\bullet^{i-1}\cong \cE_\bullet^i$. Similarly, fixing given ranks $\un{n}^i$ and degrees $\un{d}^i$, we denote by $\uExt(\un{n}^h,\dots,\un{n}^1)_{\un{d}^h,\dots,\un{d}^1}$ the stack of chains $\cE_\bullet$ together with a filtration $\cF_\bullet^i$ such that $(\rk(\cF_\bullet^i/\cF_\bullet^{i-1}),\deg(\cF_\bullet^i/\cF_\bullet^{i-1}))=(\un{n}^i,\un{d}^i)$. Also we will denote by $\uExt(\un{n}^h,\dots,\un{n}^1)_{\un{d}^h,\dots,\un{d}^1}^{\gr\,\alpha-\sst}$ the open substack of filtered chains such that the subquotients $\cF_\bullet^i/\cF_\bullet^{i-1}$ are $\alpha$-semistable.
\end{notation}

We need the following basic result, which can be found in \cite[Proposition 3.1 and 3.5]{AGS}:
\begin{proposition}\label{Ext-Rechnung}
Let $\cE^\pprime_\bullet$,$\cE_\bullet^\prime$ be chains. Then we have a long exact sequence:
\begin{align*}
 0 &\to \Hom(\cE_\bullet^\pprime,\cE^\prime_\bullet) \to \bigoplus_i \Hom(\cE^\pprime_i,\cE^\prime_i) \to \bigoplus_i \Hom(\cE^\pprime_{i},\cE^\prime_{i-1}) \\
 &\to \Ext^1(\cE^\pprime_\bullet,\cE^\prime_\bullet) \to \bigoplus_i \Ext^1(\cE^\pprime_i,\cE^\prime_i) \to \bigoplus_i \Ext^1(\cE^\pprime_i,\cE^\prime_{i-1}) \\
 &\to \Ext^2(\cE^\pprime_\bullet,\cE^\prime_\bullet)\to  0.
\end{align*}
If the $\phi_i^\pprime$ are injective for all $i$ or the $\phi_i^\prime$ are generically surjective for all $i$, then $\Ext^2(\cE^\pprime_\bullet,\cE^\prime_\bullet)=0$.
\end{proposition}

The above proposition is most useful, if the $\Ext^2$-term vanishes. We will need another criterion to show this. To this end, we recall that the long exact sequence computing $\Ext$-groups in the category of chains is obtained from the cohomology of the complex of sheaves $$\bigoplus_i \cHom(\cE_i^\pprime,\cE_i^\prime) \to \bigoplus_i \cHom(\cE_i^\pprime,\cE_{i-1}^\prime)$$ on $C$, where the differential is given by $(f_i) \mapsto (f_{i-1}\circ \phi^\pprime_i-\phi_i^\prime\circ f_i)$. In particular the last group in the sequence of Proposition \ref{Ext-Rechnung} is $\Ext^2(\cE_\bullet^\pprime,\cE_\bullet^\prime)$. Moreover, we can apply Serre-duality to this sequence and find that 
\begin{align*}
&{\bH}^i(C,\bigoplus_i \cHom(\cE_i^\pprime,\cE_i^\prime) \to \bigoplus_i \cHom(\cE_i^\pprime,\cE_{i-1}^\prime))^\vee \\
&\cong {\bH}^{2-i}(C,\bigoplus_i \cHom(\cE_{i-1}^\prime,\cE_i^\pprime\tensor \Omega_C) \to \bigoplus_i \cHom(\cE_i^\prime,\cE_{i}^\pprime \tensor \Omega_C)).
\end{align*}
The complex occurring in the second hypercohomology group is the one computing $\Ext^{2-i}(\cE_\bullet^\prime,\cE_{\bullet-1}^\pprime\tensor \Omega)$, where $\cE_{\bullet-1}^\pprime$ is the chain obtained by shifting the chain $\cE_\bullet^\pprime$ by placing $\cE_i^\pprime$ in degree $i-1$, so that the bundles of the resulting chains may be non-zero for $-1\leq i \leq r-1$. 
This implies:\footnote{S. Mozgovoy independently observed this lemma.}
\begin{lemma}\label{SerreChains}
Let $\cE^\pprime_\bullet$,$\cE_\bullet^\prime$ be chains, then we have 
$$ \Ext^2(\cE_\bullet^\pprime, \cE_\bullet^\prime)^\vee \cong \Hom(\cE_\bullet^\prime,\cE_{\bullet-1}^\pprime \tensor \Omega_C).$$
\end{lemma} 

In our applications we will mostly be interested in the case that $\alpha_i=(2g-2)i$. In this case the following lemma is the key to our the computation of the Harder--Narasimhan strata:
\begin{lemma}\label{HN-Ext-Rechnung}
Suppose for all $i$ we have $\alpha_i-\alpha_{i-1}\geq 2g-2$. 
\begin{enumerate}
\item Let $\cE^\prime_\bullet,\cE_\bullet^\pprime$ be chains of slope $\mu_{\min}(\cE_\bullet^\prime) > \mu_{\max}(\cE^\pprime_\bullet)$. Then $\Ext^2(\cE_\bullet^\pprime,\cE_\bullet^\prime)=0$. 
\item If $\cE_\bullet$ is an $\alpha$-stable chain, then $\Ext^2(\cE_\bullet,\cE_\bullet)=0$.
\end{enumerate}
\end{lemma}
\begin{proof}
Let us prove (1). The preceding lemma says that the statement is equivalent to the vanishing of
$$\Hom(\cE_\bullet^\prime, \cE_{\bullet-1}^\pprime \tensor \Omega_C).$$
%By assumption $\cE_\bullet^\prime$ is $\alpha$-semistable as is $\cE_\bullet^\pprime \tensor \Omega_C$.
The image of a morphism $f\colon \cE_\bullet^\prime \to \cE_{\bullet-1}^\pprime \tensor \Omega_C$, is a quotient of $\cE_\bullet^\prime$. Thus this image is of slope $>\mu_{\min}(\cE_\bullet^\prime)>\mu_{\max}(\cE_\bullet^\pprime)$. Also $f(\cE_{\bullet+1}^\prime)$ is a subchain of $\cE_\bullet^\pprime\tensor \Omega_C$, and we have
\begin{align*}
\mu(f(\cE_{\bullet+1}^\prime)) &= \mu(f(\cE_\bullet^\prime)) + \sum_i \frac{\rk(f(\cE_i))}{|\rk(f(\cE_\bullet^\prime))|} (\alpha_{i+1} -\alpha_i)\\
& \geq \mu(f(\cE_\bullet^\prime))+ 2g-2\\
& > \mu_{\max}(\cE_\bullet^\pprime) + 2g-2 = \mu_{\max}(\cE_\bullet^\pprime\tensor \Omega_C).
\end{align*}
This is a contradiction. %contradicts the semistability of $\cE_\bullet^\pprime\tensor \Omega_C$.

For (2) the same argument shows that any element of $\Hom(\cE_\bullet, \cE_{\bullet-1} \tensor \Omega_C)$ must be an isomorphism, which cannot exist, since $\cE_{-1}=0$.
\end{proof}
\begin{remark}
One should compare the above Lemma with \cite[Proposition 3.5]{AGS}, which gives a similar statement. 
Note however that statement (2) cannot be generalized to a pair of stable chains $\cE_\bullet^\prime,\cE_\bullet^\pprime$ of equal slope if $\alpha_i-\alpha_{i-1} = 2g-2$, e.g., the chains $\cE_\bullet^\pprime:=(\cO \to 0)$ and $\cE_\bullet^\prime=(0\to \Omega_C)$ are both stable, since they are of rank $1$, they have equal slope, but $\Ext^2(\cE_\bullet^\pprime,\cE_\bullet^\prime)=\Ext^1(\cO,\Omega)=H^0(\cO)^\vee$ is non-trivial.

In particular, in (iii) of loc.\ cit.\ the condition imposed by Lemma \ref{SerreChains} saying that $\Hom(\cE_\bullet^\prime, \cE_\bullet^\pprime\tensor \Omega_C)=0$ has to be added. This is for example satisfied if one of the chains is stable and all $\cE_i^\prime,\cE_i^\pprime$ are non-trivial.
The same applies to \cite[Proposition 3.6 (2)]{BGPG}.
\end{remark}

The above lemma allows us to describe the Harder--Narasimhan strata of chains:
\begin{proposition}\label{HN-Ext-Stack}
Let $\alpha$ be a stability parameter and $(\un{n}^i,\un{d}^i)_{i=1,\dots h}$ be ranks and degrees of chains. Suppose that $$\alpha_k-\alpha_{k-1}\geq 2g-2 \text{ for }k=1,\dots r \text{ and } \mu(\un{n}^i,\un{d}^i)>\mu(\un{n}^{i+1},\un{d}^{i+1}) \text{ for } i=1,\dots h-1.$$
Then the forgetful map: 
$$\gr:\uExt(\un{n}^h,\dots,\un{n}^1)_{\un{d}^h,\dots,\un{d}^1}^{\gr\,\alpha-\sst} \to \prod_{i=1}^h \cM(\un{n}^i)_{\un{d}^i}^{\alpha-\sst}$$
is smooth and its fibers are affine spaces of dimension $\chi=\sum_{i<j} \chi_{ij}$, where 
\begin{align*}
\chi_{ij}=& \sum_{k=0}^r n_k^i n_k^{j}(g-1) - n_k^i d_k^{j} + n_k^i d_k^{j}\\
&- \sum_{k=1}^r n_k^{j} n_{k-1}^i(g-1) - n_{k-1}^i d_k^{j} + n_k^{j} d_k^i.
\end{align*}
Moreover in $\K0hat$ we have 
$$ [\uExt(\un{n}^h,\dots,\un{n}^1)_{\un{d}^h,\dots,\un{d}^1}^{\gr\,\alpha-\sst}] = \bL^{\sum_{i<j}\chi_{ij}} \prod_{i=1}^h [\cM(\un{n}^i)_{\un{d}^i}^{\alpha-\sst}].$$ 
\end{proposition}
\begin{proof}
Let us denote the projections from the product $\uExt(\un{n}^{h-1},\dots,\un{n}^1)_{\un{d}^{h-1},\dots,\un{d}^1}^{\gr\,\alpha-\sst}\times \cM(\un{n}^h)_{\un{d}^h}^{\alpha-\sst} \times C$ onto the $i$-th factor by $p_i$ and similarly denote by $p_{ij}$ the projection onto the product of the $i$-th and $j$-th factor. Denote by $\cE^\prime_{\bullet,\univ},\cE^h_{\bullet,\univ}$ the universal chains on $\uExt(\un{n}^{h-1},\dots,\un{n}^1)_{\un{d}^{h-1},\dots,\un{d}^1}^{\gr\,\alpha-\sst} \times C$ and $\cM(\un{n}^2)_{\un{d}^2}^{\alpha-\sst}\times C$. Since $\cE_{\bullet,\univ}^\prime$ has a filtration with semistable subquotients of slope bigger than $\mu(\cE_{\bullet,\univ}^h)$, Lemma \ref{HN-Ext-Rechnung} implies that  
$$\bR p_{12,*} \Bigl(\bigoplus_i \cHom(p_{23}^*\cE_{i,\univ}^h,p_{13}^*\cE_{i,\univ}^\prime) \to \bigoplus_i \cHom(p_{23}^*\cE_{i,\univ}^h,p_{13}^*\cE_{i-1,\univ}^\prime)\Bigr)$$ is a complex with cohomology only in degrees $0,1$. As in the proof of Corollary \ref{ext-stack}, this complex can be represented by a complex of vector bundles $\cF_0\map{d_0}\cF_1\map{d_1} \cF_2$ and since its cohomology is only in degree $0,1$ it is quasi-isomorphic to $\cF_0 \to \ker(d_1)$. This is a complex of vector bundles of length $2$, to which we can apply Proposition \ref{VectorBundleStack}, showing that the vector bundle stack $[\ker(d_1)/\cF_0]$ is isomorphic to $\uExt(\un{n}^h,\dots,\un{n}^1)^{\gr \alpha-\sst}_{\un{d}^h,\dots,\un{d}^1}$. 
By the Riemann--Roch formula the dimension of the fibers of $\gr$ is $\sum_{j=1}^{h-1}\chi_{jh}$. 
The result now follows by induction on $h$.
\end{proof}

%%%%%%%%%%%%
\subsection{A stratification of the stack of chains with simple strata}\label{StratificationSection}
%%%%%%%%%%%%

Next, we want to define --- for any $\un{n},\un{d}$ --- a stratification of the stack of chains $\cM(\un{n})_{\un{d}}$ such that classes of the strata can be computed using the results of Section \ref{BasicTools}.
%compute the class of $\cM(\un{n})_{\un{d},\un{r},\un{k},\un{l}}$. 

Any chain $\cE_\bullet=(\cE_r \to \dots \to \cE_0)$ has a canonical subchain:
$$F^1\cE_\bullet= (\cE_r \to \phi_r(\cE_r)^{\sat} \to \dots \to \phi_{1r}(\cE_r)^{\sat} \subset \cE_\bullet),$$ in which all maps are generically surjective. We call this subchain the saturation of $\cE_r$ in $\cE_\bullet$ and write $(\rk(F^1\cE_\bullet),\deg(F^1\cE_\bullet))=:(\un{n}^\prime,\un{d}^\prime)$. Also we set $\un{l}:=(\deg(F^1\cE_i)-\deg(\phi_{i+1}(\cE_{i+1})))_{i=0,\dots r-1}$.

For fixed $\un{n}^\prime,\un{d}^\prime, \un{l}$, the stack of those chains such that the saturation of $\cE_r$ in $\cE_\bullet$ is of  rank $\un{n}^\prime$ and with degrees $\un{d}^\prime,\un{l}$ is a locally closed substack $\cM(\un{n})_{\un{d}}^{\un{n}^\prime,\un{d}^\prime,\un{l}^\prime}\subset \cM(\un{n})_{\un{d}}$. 

For any chain $\cE_\bullet$, the quotient $\cE_\bullet/F^1\cE_\bullet=(0 \to \cE_{r-1}/\phi_r(\cE_r)^{\sat} \to \dots \to \cE_0/\phi_{1r}(\cE_r)^{\sat})$ is a chain of shorter length. This has a corresponding subchain, the saturation of $\cE_{r-1}/\phi_r(\cE_r)^{\sat}$. Inductively this defines a filtration $$F^1\cE_\bullet \subset \dots \subset F^{r+1}\cE_\bullet=\cE_\bullet$$ such that the subquotients $F^i\cE_\bullet/F^{i-1}\cE_\bullet$ are chains of length $r-i+1$ in which all maps are generically surjective.
Proposition \ref{Ext-Rechnung} will therefore allow us to describe the substacks of chains such that the rank and degree of the chains in this filtration are constant.

Given $\un{n}$, $\un{d}$ and $\un{l}$, let us denote by $\cM(\un{n})_{\un{d}}^{\gensurj}$ the stack of chains of rank $\un{n}$ and degree $\un{d}$, such that all maps $\phi_i$ are generically surjective and denote by $\cM(\un{n})_{\un{d},\un{l}}^{\gensurj}$ the locally closed substack of $\cM(\un{n})_{\un{d}}^{\gensurj}$ defined by the condition $l_i=\deg(\cE_i)-\deg(\phi_i(\cE_{i+1}))$.  

\begin{lemma}\label{ChainStrata1}%\begin{enumerate}
The stack $\cM(\un{n})_{\un{d},\un{l}}^{\text{gen-surj}}$ is nonempty if and only if $n_r\geq n_{r-1} \geq \dots \geq n_0$, $l_i\geq 0$ and for all $i$ such that $n_i=n_{i+1}$ we have $l_i=d_{i}-d_{i+1}$. 

If these conditions are satisfied, the stack $\cM(\un{n})_{\un{d},\un{l}}^{\text{gen-surj}}$ is smooth and connected. 

Let  
$\chi_i:=\left\{\begin{array}{ll}
(n_{i+1}n_{i} - n_{i}^2)(g-1) +n_{i+1} d_{i} - n_{i}d_{i+1} & \textrm{ if } n_{i+1}\neq n_i\\
0 & \textrm{ if }n_{i+1}=n_i \end{array}\right.$.
Then we have:
$$[\cM(\un{n})_{\un{d},\un{l}}^{\gensurj}]= \Bun_{n_0}\prod_{i=0}^{r-1} [\Bun_{n_{i+1}-n_{i}}] [(C\times \bP^{n_{i}-1})^{(l_{i})}] \bL^{\sum \chi_i-n_{i+1}l_{i}}\in \K0hat.$$
\end{lemma}
\begin{proof}
For any chain in the substack write $\cK_i:=\ker(\phi_i), \cQ_i:=\im(\phi_i)$. Then we have $\rk(\cK_i)=n_i-n_{i-1}$ and $\deg(\cK_i))=d_i-(d_{i-1}-l_{i-1})$. 
Thus the stack classifies a collection of extensions $(\cK_i \to \cE_i \to \cQ_i)$ together with Hecke modifications $ \cQ_i\subset \cE_{i-1}$ of length $l_{i-1}$.

For fixed $\cK_i,\cQ_i$ the dimension of the stack of extensions of $\cQ_i$ by $\cK_i$ is 
\begin{align*}
&\rk(\cK_i)\rk(\cQ_i)(g-1)-(\rk(\cQ_i)\deg(\cK_i)-\rk(\cK_i)\deg(\cQ_i))\\
&=(n_i-n_{i-1})n_{i-1}(g-1)-(n_{i-1} (d_i-d_{i-1}+l_{i-1}) - (n_i-n_{i-1})(d_{i-1}-l_{i-1}))\\
&=(n_in_{i-1} - n_{i-1}^2)(g-1) +n_i d_{i-1} - n_i l_{i-1} - n_{i-1}d_i = \chi_i -n_il_{i-1}.
\end{align*}
Recall that for a family of bundles $\cE$ parameterized by some space $T$ we defined $\Hecke(\cE/T)^l$ to be the space of all modifications of length $l$ of $\cE$, i.e., the fibered product:
$$\xymatrix{
\Hecke^l(\cE/T) \ar[r]\ar[d] & \Hecke^l \ar[d]^{p}\\
T \ar[r] & \Bun_n^d.
}$$
Since $p$ is a smooth fibration with connected fibers, we see that if $T$ is connected $\Hecke(\cE/T)^l$ is connected as well. This proves the claimed connectedness. 

Using the formula for $[\Hecke(\cE/T)]$ (Proposition \ref{Mnn-rkn}) we find:
\begin{align*}
[(\cM(\un{n})_{\un{d},\un{l}}^{\gensurj}] &= \Bun_{n_0}\prod_{i=0}^{r-1} [\Bun_{n_{i+1}-n_{i}}] [(C\times \bP^{n_i-1})^{(l_i)}] \bL^{\sum \chi_i-n_{i+1}l_i},
\end{align*}
which is the claimed formula.
\end{proof}
\begin{corollary}\label{ChainStrata2}
The stack $\cM(\un{n})_{\un{d}}^{\gensurj}$ is non-empty if and only if $n_r\geq n_{r-1} \geq \dots \geq n_0$, and for all $i$ such that $n_i=n_{i-1}$ we have $d_i\leq d_{i-1}$.

If these conditions are satisfied, we have:
\begin{align*}
[\cM(\un{n})_{\un{d}}^{\gensurj}]&= \bL^{\sum \tilde{\chi}_{i}} [\Bun_{n_0}]\prod_{i\colon n_{i+1}\neq n_{i}} [\Bun_{n_{i+1}}]  \prod_{i\colon n_{i+1}=n_{i}} [(C\times \bP^{n_i-1})^{(d_{i}-d_{i+1})}],\\
\end{align*}
where 
$$\tilde{\chi}_{i}=\left\{\begin{array}{ll} n_{i+1}d_i-n_id_{i+1}+n_{i+1}n_i(1-g) & \textrm{ if }n_{i+1}\neq n_i\\ 0 &\textrm{ if }n_{i+1}=n_i.\end{array}\right.$$
\end{corollary}
\begin{proof}
Summing the formula obtained in Lemma \ref{ChainStrata1} all $l_i\geq 0$ we find:
\begin{align*}
[\cM(\un{n})_{\un{d}}^{\gensurj}]&= \bL^{\sum \chi_i} \Bun_{n_0}\prod_{i=0}^{r-1} [\Bun_{n_{i+1}-n_{i}}]  \prod_{i\colon n_{i+1}=n_{i}} [(C\times \bP^{n_i-1})^{(d_{i}-d_{i+1})}]\\ 
&\prod_{i \colon n_{i+1}-n_{i}>0} Z(C\times \bP^{n_i-1},\bL^{-n_{i+1}}).
\end{align*}
Moreover we know $Z(C\times \bP^{k-1},t) = \prod_{i=0}^{k-1} Z(C,\bL^it)$ and we have the formula $[\Bun_m]=\bL^{(m^2-1)(g-1)}[\uPic]\prod_{i=2}^m Z(C,\bL^{-i})$. Putting these two formulas together, we find $[\Bun_m]=\bL^{(2km-k^2)(g-1)}[\Bun_{m-k}]Z(C\times \bP^{k-1},\bL^{-m})$. To conclude we only need to substitute this expression for $m=n_{i+1},k=n_i$ in the above formula for $[\cM(\un{n})_{\un{d}}^{\gensurj}]$.
\end{proof}
\begin{remark}\label{FunnyFormula}
In case $n_r> \dots > n_1>n_0$ the above formula reduces to:
\begin{align*}
[\cM(\un{n})_{\un{d}}^{\gensurj}]&= \bL^{\sum_{i=1}^r n_id_{i-1}-n_{i-1}d_i+n_in_{i-1}(1-g)}\prod_{n=0}^r \Bun_{n_i}. 
\end{align*}
This formula looks quite surprising,  because in this case the class of $\cM(\un{n})_{\un{d}}^{\gensurj}$ is equal to the class of a vector bundle over $\prod_{i=0}^r\Bun_{n_i}^{d_i}$ of rank equal to the Euler characteristic of $\bigoplus_{i=1}^r \bR\Hom(\cE_i,\cE_{i-1})$, although the forgetful map $\cM(\un{n})_{\un{d}}^{\gensurj} \to \prod_{i=0}^r \Bun_{n_i}$ is far from being a bundle for any $\un{d}$.
\end{remark}
\begin{remark}\label{geninj}
By duality (Lemma \ref{dual}) the class of the stack of chains, such that all morphisms $\phi_i$ are injective, is also given by the expression given in Corollary \ref{ChainStrata2}.
\end{remark}
Lemma \ref{ChainStrata1} also allows us to define the stratification we are looking for. Given $\un{n},\un{d}$, we consider partitions $\un{n}=\sum_{i=0}^{r+1} \un{n}^i$, such that $n_k^i=0$ for $k>n-i$, so $n_r^0=n_r$, $n_{r-1}^1=n_{r-1}-n_{r-1}^0$, etc. Let $(\un{d}^i)_{i=1,\dots r}$ be such that $\sum_i \un{d}^i =\un{d}$ and such that $\un{n}^i,\un{d}^i$ satisfy the condition of Corollary \ref{ChainStrata2}. We call $\un{n}^i,\un{d}^i$ satisfying these conditions a {\em partition for the stratification by saturations}. Given such a partition let us denote $\cM(\un{n})_{\un{d}}^{(\un{n}^i,\un{d}^i)}\subset\cM(\un{n})_{\un{d}}$ the locally closed substack of chains $\cE_\bullet$ such that $(\un{n}^i,\un{d}^i)=(\rk(F^i\cE_\bullet/F^{i-1}\cE_\bullet),\deg(F^i\cE_\bullet/F^{i-1}\cE_\bullet))$.

\begin{proposition}\label{ChainStrata}
For any $\un{n},\un{d}$, the set of substacks $\cM(\un{n})_{\un{d}}^{(\un{n}^i,\un{d}^i)}\subset\cM(\un{n})_{\un{d}}$ indexed by partitions for the stratification by saturations defines a stratification of $\cM(\un{n})_{\un{d}}$.
Moreover in $\K0hat$ we have
\begin{align*}
[\cM(\un{n})_{\un{d}}^{(\un{n}^i,\un{d}^i)}]=& \bL^{\chi((\un{n}^i,\un{d}^i)_{i=0,\dots,r})}\prod_{i=0}^r [\cM(\un{n}^i)_{\un{d}^i}^{\text{gen-surj}}], \\
\text{where }\chi((\un{n}^i,\un{d}^i)_{i=0,\dots,r})=& \sum_{i<j} \big(\sum_{k=0}^{r-j} n_k^in_k^j(g-1) +n_k^id_k^j-n_k^jd_k^i 
\\ &- \sum_{k=1}^{r-j} n_{k-1}^in_k^j(g-1) + n_{k-1}^i d_k^j - n_{k}^j d_{k+1}^i\big).
\end{align*}
\end{proposition}
\begin{proof}
This holds, because by Proposition \ref{Ext-Rechnung} the canonical map $$\gr \colon \cM(\un{n})_{\un{d}}^{(\un{n}^i,\un{d}^i)} \to \prod_{i=0}^r \cM(\un{n}^i)_{\un{d}^i}^{\text{gen-surj}}$$ is a composition of vector bundle stacks and their dimension is given by the Riemann-Roch formula.
\end{proof}

%%%%%%%%%%%%%%%%%%%%%
\section{Application to Higgs bundles with fixed determinant}\label{Application}
%%%%%%%%%%%%%%%%%%%%%

In this section we give an application of Lemma \ref{ChainStrata1}, namely we prove that the action of the $n$-torsion points of $\Pic^0_C$ on the middle-dimensional cohomology of the space of rank $n$ Higgs bundles with fixed determinant is trivial. This answers a question of T.\ Hausel, which was motivated by the mirror symmetry conjecture of \cite{HT1}. The result will not be used in the rest of our article, but gives an instance where the results of the previous section can be applied without evaluating explicit formulas.
In order to simplify the exposition we will assume in this section that the characteristic of the ground field does not divide $n$. We denote our cohomology theory by $H^*(\cM)$, which means singular cohomology in case $k=\bC$ and  \'etale cohomology with $\bQ_\ell$-coefficients otherwise.

For a line bundle $L\in \Pic^d(C)$ we define $$\cM_{n}^L:=\left\langle (\cE,\phi\colon \cE\to \cE \tensor \Omega_C, \psi\colon \det(\cE) \map{\cong} L) \mid \tr(\phi)=0 \right\rangle$$ to be the stack of Higgs bundles of rank $n$ with fixed determinant $L$.

Tensoring with line bundles defines an action of the group of $n$-torsion points $\Gamma:=\Pic^0(C)[n]$ on this stack and this action respects the semistable locus $\cM_{n}^{L,\sst}$.

\begin{theorem}\label{TrivialAction}
Suppose $(n,\deg(L))=1$. Then the action of $\Gamma$ on $\cM_n^{L,\sst}$ induces the trivial action on $H^{\dim(\cM_n^L)}(\cM_n^{L,\sst})$.
\end{theorem}

\begin{proof}
As in the case of Higgs bundles without fixed determinant (Section \ref{recollection}), the cohomology of $H^*(\cM_n^{L,\sst})$ is a direct sum over cohomology groups for the fixed-point strata for the $\bG_m$-action on $\cM_n^{L,\sst}$. Moreover, these fixed-point strata are moduli spaces of semistable chains with fixed determinant. 
We denote these strata by $\cM(\un{n})_{\un{d}}^{L,\sst}$ and write $\cM(\un{n})_{\un{d}}^L$ for the corresponding stacks of all (not necessarily semistable) chains.

Since we assumed that $(n,\deg(L))=1$, the spaces of semistable chains with fixed determiant are projective, so their cohomology coincides with cohomology with compact supports. 

As in the previous section the spaces $\cM(\un{n})_{\un{d}}^L$ come equipped with a natural stratification defined by the saturation of the images of the bundles $\cE_i$. This stratification is obtained by pull-back from the stratification of the stack of chains with arbitrary determinant.

Let us fix a partition $\un{n}^i,\un{d}^i$ of $\un{n},\un{d}$ for the stratification by generic rank and fix moreover $\un{l}^i$ satisfying the numerical conditions of Lemma \ref{ChainStrata1}. Let us denote by $\cM(\un{n})_{\un{d}}^{\un{n}^i,\un{d}^i,\un{l}^i}\subset\cM(\un{n})_{\un{d}}^{\un{n}^i,\un{d}^i}$ the substack, such that $F^i\cE_\bullet/F^{i-1}\cE_\bullet \in \cM(\un{n}^i)_{\un{d}^i,\un{l}^i}^{\text{gen-surj}}$ for all $i$.
The action of $\Pic^0(C)$ on $\cM(\un{n})_{\un{d}}$ respects these substacks.

Also, by Lemma \ref{ChainStrata1} the stacks $\cM(\un{n}^i)_{\un{d}^i,\un{l}^i}^{\text{gen-surj}}$ are smooth and connected, and by Proposition \ref{Ext-Rechnung} the map 
$$\gr \colon \cM(\un{n})_{\un{d}}^{(\un{n}^i,\un{d}^i,\un{l}^i)} \to \prod_{i=0}^r \cM(\un{n}^i)_{\un{d}^i,\un{l}^i}^{\text{gen-surj}}$$
is a composition of vector bundle stacks. Thus the stacks $\cM(\un{n})_{\un{d}}^{(\un{n}^i,\un{d}^i,\un{l}^i)}$ are smooth and connected as well. 

We define $\cM(\un{n})_{\un{d}}^{(\un{n}^i,\un{d}^i,\un{l}^i),L}$ to be the analogous stack of chains with fixed determinant, i.e., the stack $\cM(\un{n})_{\un{d}}^{(\un{n}^i,\un{d}^i,\un{l}^i),L}$ is the fiber over $L$ of the determinant map $$\det\colon \cM(\un{n})_{\un{d}}^{(\un{n}^i,\un{d}^i,\un{l}_i)} \to \uPic.$$ 

The action of $\Gamma=\Pic^0(C)[n]$ respects these strata, so it is sufficient to show that $\Gamma$ acts trivially on the top cohomology with compact supports of each stratum, i.e., to show that the strata are connected.

By our description of $\cM(\un{n})_{\un{d}}^{(\un{n}^i,\un{d}^i,\un{l}_i)}$ the map $\det$ factorizes as:
\begin{align*}
\cM(\un{n})_{\un{d}}^{(\un{n}^i,\un{d}^i,\un{l}_i)} &\map{p_1} \prod_{i=0}^r \cM(\un{n}^i)_{\un{d}^i,\un{l}^i}^{\text{gen-surj}}\\
&\map{p_2} \prod_{i=0}^r \prod_{k=0}^{r-i} \Bun_{n_k^i-n_{k-1}^i}^{d_k^i-d_{k-1}^i-l_k^i} \times \prod C^{(l_k^i)}\\
&\map{p_3} \prod_{i=0}^r \prod_{k=0}^{r-i} \Pic^{d_k^i-d_{k-1}^i-l_k^i} \times \Pic^{l_k^i}\\
&\map{m} \Pic^d.
\end{align*}

Here, the map $p_1$ is given by mapping $\cE_\bullet$ to the subquotients $F^i\cE_\bullet/F^{i-1}\cE_\bullet$ of the filtration given by the saturation of the constituents $\cE_r,\dots,\cE_1$. This map is a composition of vector bundles stacks, so it has connected fibers. 

The map $p_2$ is the product of the natural maps $$\gr\colon \cM(\un{n}^i)_{\un{d}^i,\un{l}^i}^{\text{gen-surj}} \to \prod_{k=0}^{r-i} \Bun_{n_k^i-n_{k-1}^i}^{d_k^i-d_{k-1}^i-l_k^i} \times \prod C^{(l_k^i)}$$
sending a chain $\cE_\bullet$ to $(\cE_0,\ker(\phi_k)_{k=1\dots,r-i}, \supp(\cE_{k-1}/\phi_k(\cE_k))_{k=1,\dots,r-i})$, where we write $\supp(\cE_{k-1}/\phi_k(\cE_k))$ for the divisor defined by the torsion sheaf $\cE_{k-1}/\phi_k(\cE_k)$. This map is a composition of vector bundle stacks and Hecke modifications, so again this map has connected fibers.

The map $p_3$ is the product of the natural maps $\det\colon \Bun_m^e \to \Pic^e$ and $C^{(e)}\to \Pic^e$, which have connected fibers as well. 

Finally, the map $m$ is the map reconstructing the determinant of $\cE_\bullet$ out of the graded pieces defined by the $p_i$. It is of the form $m\colon(L_1,\dots,L_N)\mapsto \bigotimes_{i=1}^N L_i^{k_i}$, so this map has connected fibers if and only if $\gcd(k_1,\dots,k_N)=1$. We claim that this condition is satisfied, because of our coprimality assumption on rank and degree of the corresponding Higgs bundles. This is elementary, but slightly tedious:

Let us first consider a single factor: 
$$\det\colon\cM(\un{n}^i)_{\un{d}^i,\un{l}^i}^{\text{gen-surj}} \map{p_2\circ p_1} \prod_{i=0}^{r-i} \Pic^{d_k^i-d_{k-1}^i-l_k^i} \times C^{(l_i)} \map{m} \Pic.$$
For a chain $\cE_\bullet\in \cM(\un{n}^i)_{\un{d}^i,\un{l}^i}^{\text{gen-surj}}$ write $K_j:=\ker(\phi_j)$, $K_0:=\cE_0$ and let $D_j$ denote the divisor defined by the torsion sheaf $\cE_{j-1}/\phi_j(\cE_j)$. Then $$\det(\bigoplus \cE_j) = \bigotimes_{j=0}^r \det(K_j)^{r+1-j} \tensor \bigotimes_{j=1}^{r} \cO(-D_j)^{r+1-j}.$$

Suppose $k\geq 1$ divides all exponents occurring in this expression for $K_j\neq 0$ and $D_j\neq 0$, i.e., suppose  $k|r+1-j$ for all $j$ such that $K_j\neq 0$ or $D_j\neq 0$. Since $$\rk (\bigoplus\cE_j) = \sum_{j=0}^r (r+1-j) \rk(K_j),$$ this implies $k|\rk(\bigoplus\cE_j)$. 
Moreover, $$\deg(\bigoplus \cE_j) = \sum_{j=0}^r (r+1-j) \deg(K_j) + \sum_{j=1}^r (r+1-j) \deg(D_j),$$ thus $k|\deg(\bigoplus\cE_j)$.

Finally the degree of the corresponding Higgs bundle is 
\begin{align*}
\deg(\bigoplus \cE_j\tensor \Omega^{-(r-j)})&= \deg(\bigoplus \cE_j) - \sum_{j=0}^r (r-j)(2g-2)\bigg(\sum_{l=0}^j \rk(K_l)\bigg)\\
&= \deg(\bigoplus \cE_j) - \sum_{j=0}^r (2g-2)\rk(K_j) \frac{(r-j)(r-j+1)}{2}.
\end{align*}
Therefore we also have $k|\deg(\bigoplus \cE_i\tensor \Omega^{-(r-i)})$. 

Taking the product over all factors $\cM(\un{n}^i)_{\un{d}^i,\un{l}^i}^{\text{gen-surj}}$ this shows that, if $k$ divides all the exponents occurring in the map $m$, then $k$ has to divide $n$ and $\deg(L)$ so that $k=1$.
\end{proof}

%%%%%%%%%
\section{Moduli stacks of chains: recursion formulas and examples}\label{ChainsExamples}
%%%%%%%%%

In this section we will explain our strategy to compute the cohomology of moduli spaces of chains by giving several examples. Whenever the stack of chains of fixed invariants defines a class in $\K0hat$ our strategy immediately gives recursion formulas for the class of the space of $\alpha$-semi-stable chains whenever $\alpha$ satisfies $\alpha_{i+1}-\alpha_i\geq 2g-2$. We will say that $\alpha$ is  bigger or equal to $(0, 2g-2,\dots,r(2g-2))$, the stability parameter occurring in the study of Higgs bundles, if this condition is satisfied. We will begin with examples where this strategy works without further effort.

In general the stack of chains will have infinitely many strata of the same dimension, so that the sum over all strata does not converge. However, we know a priori that the stack of stable chains of fixed rank and degree is of finite type. In particular the convergence problem only stems from the fact that only finitely many strata will contain stable chains. We use this observation in some examples, in order to avoid this convergence problem.

In particular, we will explicitly work out the recursion formulas in the cases needed for our application to Higgs bundles of rank $4$ and odd degree. In this case the stability parameter for the corresponding chains is $\alpha=(0,2g-2,\dots,r(2g-2))$ and since there are no strictly semistable Higgs bundles this stability parameter will not be a critical value. Thus we may as well replace $\alpha$ by $\alpha_+$, such that $\alpha_+-\alpha>0$ is irrational and small. This will be helpful, because $\alpha$ might not be good for the spaces occurring in the Harder--Narasiman strata.

We use this only to simplify our formulas. We could as well use $\alpha$, but then we would have to include chains for which some of the maps are $0$, but this would increase the length of the formulas.

%%%%%%%%%%%%%
\subsection{Stacks of semistable chains of rank $(m,1,\dots,1)$}
%%%%%%%%%%%%

For all $m\in \bN, m>1$ and degrees $\un{d}=(d_0,d_1, \dots d_r)$ the stack $\rcM(m,1,\dots,1)_{\un{d}}$ is empty unless $d_r \leq d_{r-1} \leq \dots \leq d_{1}$ and that if this condition on $\un{d}$ holds we have (Remark \ref{geninj}):
\begin{align}\label{Mm111}
[\rcM(m,1,\dots,1)_{\un{d}}] &= \bL^{(d_0-md_1)+m(1-g)}[\Bun_{m}][\uPic]\prod_{i=2}^r [C^{(d_i-d_{i-1 })}].
\end{align}
Also the stack $\cM(m,1,\dots,1)_{\un{d}}$ is stratified by the substacks defined by the condition $\phi_i=0$ for $i\in I$ where $I\subset\{1,\dots, r\}$. To specify a subset $I$ is equivalent to the choice of an ordered partition $r+1 = \sum_{i=1}^l (r_i+1)$, where the $r_i$ are given by the length of the subchains with $\phi_i\neq 0$.  
Thus we find:
\begin{align*}
[\cM(m,1,\dots,1)_{\un{d}}] &= \sum_{k=1}^{r+1} \sum_{r_i\geq 0 \atop \sum_{i=1}^k r_i+1=r+1 } \prod_{j=1}^{k-1} [\rcM(\underbrace{1,\dots,1}_{r_j+1-\text{times}})_{\un{d}}] [\rcM(n,\underbrace{1,\dots,1}_{r_k-\text{times}})_{\un{d}}].
\end{align*}
Here the index $\un{d}$ on the right hand side of the formula refers to the corresponding subset of the degrees $d_i$.

\begin{example}[Rank $(m,1)$]\label{Mm1}
\begin{align*}
[\rcM(m,1)] &= \bL^{-m(g-1)+d_0-md_1} [\uPic][\Bun_m],\\ [\cM(m,1)] &=  (\bL^{-2(g-1)+d_0-2d_1}+1) [\uPic][\Bun_m].
%[\rcM(3,1)] &= \bL^{-3(g-1)+d_0-3d_1} [\uPic][\Bun_{3}]\\ [\cM(3,1)] &= (\bL^{-3(g-1)+d_0-3d_1}+1) [\uPic][\Bun_{3}]\\
\end{align*}
\end{example}

Next, we want to apply the general recursive procedure in order to compute the class of $\cM(m,1,\dots,1)^{\alpha-\sst}$. We have to study the Harder--Narasimhan strata, assuming that the stability parameter $\alpha$  satisfies $\alpha_{i+1}-\alpha_{i} \geq 2g-2$ and we will furthermore assume that $\alpha$ is good.

Let $\cE_\bullet\in\cM(m,1,\dots,1)_{\un{d}}$ be a chain and let $(\cE_\bullet^1 \subset \dots \subset \cE_\bullet^h=\cE_\bullet)$ be its Harder--Narasimhan flag. Since $\alpha$ is good and the subquotients $\cE_\bullet^i/ \cE_\bullet^{i-1}$ are semistable we find that 
$$\rk(\cE_\bullet^i/\cE_\bullet^{i-1})=\left\{\begin{array}{ll} (m_0^i,1,\dots,1,0,\dots,0) & \text{ if } m_0^i \neq 0 \\ (0,\dots,0,1,\dots,1,0,\dots,0) & \text { if } m_0^i=0 \end{array}\right. .$$
Removing the $0$ constituents of the subquotient we see that any such subquotient is a chain of rank $(m^\prime,1,\dots,1)$ of possibly shorter length. 

The type of the flag is therefore given by a partition $m_0^1+\dots +m_0^h=m_0$ with $m_0^i \geq 0$, given by the ranks of the subquotients $\cE_0^i/\cE_0^{i-1}$, a partition $r=r^1+\dots+ r^h$ with $r^i\geq 0$ giving the length of the sequences of $1$'s in the rank of the subquotients together with an index $l_i$, specifying the starting index, i.e., we set $l_i:= \min  \{ l\mid \cE_l^i/\cE_l^{i+1} \neq 0 \}$ and 
$$r_i:=\left\{\begin{array}{ll} \max\{l>0 \mid \cE_l^i/\cE_l^{i+1} \neq 0 \} - l_i+1 &\textrm{ if } \exists l>0 \textrm{ such that } \cE_l^i/\cE_l^{i+1} \neq 0\\ 0 &\textrm{ otherwise.}\end{array}\right.$$

These data satisfy $l_i=0$, if $n_0^i \neq 0$, and for all $1\leq j\leq
r$ there is exactly one index $i$, such that $l_i \leq j <
l_i+r_i$. Summing up, the 
Harder--Narasimhan  strata are indexed by partitions of $n_0$ and $r$ together with starting indices $\un{l}$ and a set of degrees $\un{d}^i$, such that $\un{d}=\sum \un{d}^i$.

In order to write down our formulas we use the notation $\cM(m,\un{1}_{k})$ to denote the stack of chains $\cE_k \to \dots \to \cE_0$ of rank $(m,1,\dots,1)$. For $m=0$ this is isomorphic to the stack of chains $\cE_{k-1} \to \dots \to \cE_0$ of rank $(1,\dots,1)$. 

Further, for a type $(m_0^i,r_i,l_i,\un{d}^i)$ of a Harder--Narasimhan stratum occurring in $\cM(m,1,\dots,1)_{\un{d}}$ we denote the corresponding stratum by $\cM(m,\un{1}_r)^{(m_0^i,r_i,l_i,\un{d}^i)}$.

\begin{proposition}
Fix a degree $\un{d}$ and let  $\alpha$ be a good stability parameter satisfying $\alpha_{i+1}-\alpha_{i} \geq 2g-2$. 
Let $(m_0^i,r_i,l_i,\un{d}^i)$ be the type of a Harder--Narasimhan stratum occurring in $\cM(m,1,\dots,1)_{\un{d}}$.  Then we have:
$$[\cM(m,\un{1}_r)^{(m_0^i,r_i,l_i,\un{d}^i)}]= \prod_{i<j} \bL^{\chi_{ij}} \cM^{\alpha-\sst}(m_0^i,\un{1}_{r_i}) \cM^{\alpha-\sst}(m_0^j,1_{r_j}),$$
where $\chi_{ij}$ is the Euler characteristic computed in Proposition \ref{HN-Ext-Stack}.

In particular this gives a recursive formula for
$$[\cM(m,\un{1}_r)^{\alpha-\sst}]=[\cM(m,\un{1}_r)] - \sum_{(m_0^i,r_i,l_i,\un{d}^i) \textrm{ HN-type}} [\cM(m,\un{1}_r)^{(m_0^i,r_i,l_i,\un{d}^i)}].$$
\end{proposition}

\begin{remark}
In the same way one can obtain a recursion formula for the space of semistable chains of rank $(m_0,1,\dots,1,m_r)$ for any $m_0,m_r\geq 1$ for any $r\geq 2$. 
\end{remark}

Let us evaluate this recursion formula for $m=2$ and $m=3$. After tensoring with a fixed line bundle we may assume that $d_1$ in both cases:
\begin{example}[Rank $(2,1)$]\label{Rank21}
Let $\alpha=(0,\sigma)$ be a good stability parameter with $\sigma\geq 2g-2$. 
\begin{enumerate}
\item $\cM(2,1)_{d_0,0}^{\alpha-\sst}$ is empty unless $\frac{\sigma}{2} \leq d_0 \leq 2\sigma$. If these inequalities hold then: 
\begin{align*}
[\cM(2,1)_{d_0,0}^{\alpha-\sst}] &= [\rcM(2,1)_{d_0,0}] - \sum_{l=\lceil\frac{2d_0-\sigma}{3}\rceil}^\infty [\uPic]^2 [C^{(l)}]\bL^{g-1+d_0-2l} - \sum_{l=0}^{\lceil \frac{2d_0-\sigma}{3}\rceil-1} [C^{(l)}]\bL^{l}\\
&= [\uPic]^2 \sum_{l=0}^{\lceil \frac{2d_0-\sigma}{3}\rceil-1} [C^{(l)}](\bL^{g-1+d_0-2l} - \bL^l).
\end{align*}
\item $\mathring{\cM}(2,1)_{d_0,0}^{\mu_{\min}>h}$ is empty unless $\min\{\sigma,d_0,\frac{d_0+\sigma}{3}\}>h$. If this holds then
\begin{align*}
[\rcM(2,1)_{d_0,0}^{\mu_{\min}>h}] &= [\rcM(2,1)_{d_0,0}] - \sum_{l=\lceil d_0-h\rceil}^\infty [\uPic]^2 [C^{(l)}]\bL^{g-1+d_0-2l} - \sum_{l=0}^{\lfloor 2h-\sigma \rfloor} [C^{(l)}]\bL^{l}\\
&=[\uPic]^2 \bigg(\sum_{l=0}^{\lceil d_0-h\rceil-1}[C^{(l)}] \bL^{g-1+d_0-2l}- \sum_{l=0}^{\lfloor 2h-\sigma\rfloor} [C^{(l)}]\bL^l \bigg).\\
\end{align*}
\end{enumerate}
\end{example}
\begin{proof}
To prove (1) first note that since $\alpha$ is good for any semistable triple $\cE_1\map{\phi_1} \cE_0$ the map $\phi_1$ is non trivial. Therefore $\cE_\bullet$ contains $0\to \cE_0$ and $\cE_1 \to \cE_1$ as subtriples. Thus there are no semistable triples unless $\frac{d_0}{2} \leq \frac{d_0+\sigma}{3} \Leftrightarrow d_0 \leq 2\sigma$ and $\frac{\sigma}{2}\leq \frac{d_0+\sigma}{3} \Leftrightarrow \frac{\sigma}{2}\leq d_0$. This proves the bounds on $d_0$.

Let us now list the Harder--Narasimhan strata according to the rank of the last step $\cE_\bullet^\prime := \cE^{h-1}_\bullet\subset \cE_\bullet$. Since $\cM(2,1)_{d_0,0}^{\alpha-\sst}\subset \rcM(2,1)$ we only need to determine the intersection of the HN-strata with $\rcM(2,1)$. We write $\un{d}^\prime:=\deg(\cE_\bullet^\prime)$:
%\begin{itemize}
\subsubsection*{Type $(1,1)$} The bounds on the degree are:
$$\mu(\cE^\prime_\bullet)>\mu(\cE_\bullet) \Leftrightarrow \frac{d_0^\prime+\sigma}{2}>\frac{d_0+\sigma}{3} \Leftrightarrow d_0^\prime > \frac{2d_0-\sigma}{3}.$$
In this case the Harder--Narasimhan stratum is one of the strata $\cM(2,1)_{l}$. Its class is $\bL^{(g-1)+d_0-2d_0^\prime} [\uPic]^2[C^{(d_0^\prime)}]$. Thus the sum over all these strata is 
$$\sum_{l=d^\prime_0>\frac{2d_0-\sigma}{3}} \bL^{(g-1)+d_0-2l} [\uPic]^2[C^{(l)}].$$
\subsubsection*{Type $(1,0)$} $\mu(\cE^\prime_\bullet)>\mu(\cE_\bullet)$ implies $d_0^\prime > {d_0+\sigma}/{3}$.
The quotient $\cE_1 \map{\phi^\prime} \cE_0/\cE_0^\prime$ has to be semistable. This implies $\phi^\pprime\neq 0$ and $0 \leq d_0-d_0^\prime < \sigma$. Here the second condition is implied by $\frac{d_0+\sigma}{3}< d_0^\prime$, since $d_0\leq 2\sigma$.
The dimension of $\uExt(\cE_\bullet/\cE_\bullet^\prime,\cE_\bullet^\prime)$ is $d_0-d_0^\prime$. 
Thus the class of a stratum is $\bL^{d_0-d_0^\prime}[\uPic]^2 [C^{(d_0-d_0^\prime)}]$ and putting $l=d_0^\prime-\lfloor\frac{d_0+\sigma}{3}\rfloor$ we find that the sum over all strata is:
$$\sum_{l=0}^{\lceil\frac{2d_0-\sigma}{3}\rceil-1} \bL^{l}[\uPic]^2[C^{(l)}].$$
%\end{itemize}
This proves (1).

For the second part we only have to adjust the inequalities.
First, any triple has the quotients $\cE_1 \to 0$ and $0\to \cE_0/\cE_1$. This implies that there are no triples with $\mu_{\min}>h$ unless $\sigma >h$ and $d_0>h$. Also we need $\mu(\cE_\bullet)>h$.
Moreover the Harder--Narasimhan strata of rank $(1,1)$ are in $\rcM(2,1)^{\mu_{\min}>h}$ unless $d_0-d_0^\prime \leq h$, i.e., $d_0^\prime \geq d_0-h$. 
For rank $(1,0)$ we only need to discard the strata with $\frac{d_0-d_0^\prime+\sigma}{2}\leq h$, i.e., $d_0-d_0^\prime \leq 2h-\sigma$. 
%Note that $2h-\sigma < \frac{2d_0-\sigma}{3} \Leftrightarrow h < \frac{d_0+\sigma}{3}=\mu(\cE_\bullet)$.
This proves (2).
\end{proof}
\begin{remark}
Using Example \ref{Hodge} we can read off the Hoge polynomial of $\cM(2,1)^{\sst}$ from the above formula and thereby obtain a rather short proof of the main result of \cite{MunozPairs}.
\end{remark}
\begin{example}[Rank $(3,1)$]\label{m31} Let $\alpha=(0,\sigma)$ be a good stability parameter satisfying $\sigma\geq 2g-2$. Then the space $\cM(3,1)^{\alpha-\sst}$ is empty unless $\sigma \leq d_0 \leq 3\sigma$. If these inequalities  hold, we have
\begin{align*}
[\cM(3,1)_{d_0,0}^{\alpha-\sst}] &= [\uPic][\Bun_3]\bL^{d_0-3(g-1)}\\
&- [\uPic]^2[\Bun_2] \left( \frac{\bL^{2d_0-2\lfloor\frac{3d_0-\sigma}{4}\rfloor}}{(\bL^2-1)} + \sum_{l=\lfloor \frac{d_0+\sigma}{4}\rfloor+1}^{\lceil d_0-\frac{\sigma}{2}\rceil-1} \bL^{-3l+2d_0-(g-1)} \right)\\
& - [\uPic][\Bun_2] \left( \sum_{l=0}^{\lceil \frac{d_0-\sigma}{2}\rceil-1} \bL^{2l}[C^{(l)}] + \sum_{l=\lfloor \frac{d_0-\sigma}{2}\rfloor +1}^\infty \bL^{-3l+d_0+2g-2}[C^{(l)}]\right) 
\\
& + [\uPic]^3 \left(\sum_{k=\lfloor \frac{d_0-\sigma}{2}\rfloor+1}^\infty [C^{(k)}]\bL^{-2k} \sum_{l=\lfloor \frac{3d_0-\sigma}{4}\rfloor +1}^{\infty} \bL^{2d_0+3g-3-2l} \right. \\
& + \sum_{k=\lfloor\frac{d_0-\sigma}{2}\rfloor+1}^\infty [C^{(k)}] \bL^{-2k} \sum_{l=\lfloor\frac{d_0+\sigma}{4}\rfloor+1}^{\lceil d_0-\frac{\sigma}{2}\rceil-1} \bL^{2g-2+2d_0-3l}\\ 
&+ \sum_{k=0}^{\lceil \frac{d_0-\sigma}{2}\rceil-1}[C^{(k)}] \bL^{-2k}  \sum_{l=\lfloor \frac{3d_0-\sigma}{4}\rfloor-k+1}^{d_0-\frac{\sigma}{2}}\bL^{2g-2+2d_0-3l} \\
&+ \left. \sum_{k=0}^{\lceil \frac{d_0-\sigma}{2}\rceil-1} [C^{(k)}]\bL^k\sum_{l=\lfloor \frac{d_0+\sigma}{4}\rfloor+1}^{\infty} \bL^{g-1+d_0-2l}\right).
\end{align*}
\end{example}
\begin{proof}
This is a bit tedious, but not difficult. 
As before, since $\alpha$ is good for any semistable triple $\cE_1\map{\phi_1} \cE_0$, the map $\phi_1$ is non trivial. Therefore $\cE_\bullet$ contains $0\to \cE_0$ and $\cE_1 \to \cE_1$ as subtriples. Thus there are no semistable triples unless $\frac{d_0}{3} \leq \frac{d_0+\sigma}{4} \Leftrightarrow d_0 \leq 3\sigma$ and $\frac{\sigma}{2}\leq \frac{d_0+\sigma}{4} \Leftrightarrow \sigma\leq d_0$. This proves the bounds on $d_0$.

We will group the HN-strata according to the rank of the last step of the HN-filtration $\cE^\prime_\bullet:=\cE_\bullet^{h-1}\subsetneq \cE_\bullet$. In particular $\cE_\bullet/\cE_\bullet^\prime$ must be semistable. We write $\deg(\cE_\bullet^\prime)=:\un{d}^\prime$. 
Finally, to add the contributions of the HN-strata we will --- as in the computation of $[\Bun_n^{d,\sst}]$ --- write each stratum as $$(\textrm{extensions of all chains})-(\textrm{contribution of the unstable locus}).$$ As in the case of vector bundles, there are cancellations between unstable contributions of different strata. Thus, it will be useful to use the same parameterizations in each occurrence.
%\begin{itemize}
\subsubsection*{Type $(2,1)$}
The bounds on the degree are:
$$\mu(\cE_\bullet/\cE_\bullet^\prime)<\mu(\cE_\bullet) \Leftrightarrow d_0-d_0^\prime<\frac{d_0+\sigma}{4} \Leftrightarrow d_0^\prime > \frac{3d_0-\sigma}{4}$$
and
$$\mu_{\min}(\cE_\bullet^\prime)>\mu(\cE_\bullet/\cE_\bullet^\prime)=d_0-d_0^\prime.$$

By Proposition \ref{HN-Ext-Stack} we have: $$\dim \uExt(\cE_\bullet/\cE_\bullet^\prime,\cE^\prime_\bullet) = 2(g-1)-d_0^\prime+2(d_0-d_0^\prime)=2d_0-3d_0^\prime+2g-2.$$
From Example \ref{Rank21} we know that the conditions on $\cM(2,1)_{d_0^\prime,0}^{\mu_{\min}>d_0-d_0^\prime}$ to be non-empty, are given by $\min\{\mu(\cE^\prime_\bullet),\mu(\cE_1\to 0),\mu(0\to\cE_0^\prime/\cE_1)\}>d_0-d_0^\prime=\mu(\cE_\bullet/\cE^\prime_\bullet)$. This condition is automatically satisfied for the first two terms, because their slope is $>\mu(\cE_\bullet)$ and it also holds for the last because  $\mu(0\to\cE_0^\prime/\cE_1)\leq \mu(\cE_\bullet/\cE^\prime_\bullet)<\mu(\cE_\bullet)$ implies that $\cE_1\to \cE_1$ is destabilizing, which we excluded. Thus the contribution of the strata is:
\begin{align*}
&\sum_{d_0^\prime=\lfloor \frac{3d_0-\sigma}{4}\rfloor+1}^\infty \bL^{2d_0-3d_0^\prime+2g-2} [\uPic] [\rcM(2,1)_{d_0^\prime,0}^{\mu_{\min}>d_0-d_0^\prime}]\\
&\stackrel{(\ref{Rank21})}{=} \sum_{k=\lfloor \frac{3d_0-\sigma}{4}\rfloor+1}^\infty \bL^{2d_0+2g-2-3k} [\uPic] [\rcM(2,1)_{k,0}]\\
&- \sum_{k=\lfloor \frac{3d_0-\sigma}{4}\rfloor+1}^\infty \sum_{l=2k-d_0}^\infty \bL^{2d_0+3g-3-2k-2l} [C^{(l)}] [\uPic]^3\\
&- \sum_{k=\lfloor \frac{3d_0-\sigma}{4}\rfloor+1}^\infty \sum_{l=0}^{2d_0-2k-\sigma} \bL^{2d_0+2g-2-3k+l}[C^{(l)}] [\uPic]^3 \\
&\stackrel{(\ref{Mm1})}{=} \frac{\bL^{2 d_0}[\uPic]^2[\Bun_2]}{(\bL^2-1)\bL^{2\lfloor\frac{3d_0-\sigma}{4}\rfloor}} - \sum_{k=\lfloor \frac{d_0-\sigma}{2}\rfloor+1}^\infty \sum_{l=\lfloor \frac{3d_0-\sigma}{4}\rfloor +1-k}^{\lfloor \frac{d_0-k}{2}\rfloor} \bL^{2d_0+3g-3-4k-2l} [C^{(k)}] [\uPic]^3 \\
&- \sum_{k=0}^{\lceil \frac{d_0-\sigma}{2}\rceil-1} \sum_{l=\lfloor\frac{3d_0-\sigma}{4}\rfloor+1-k}^{\lceil d_0-\frac{\sigma+k}{2}-k\rceil-1} \bL^{2d_0+2g-2-3l-2k}[C^{(k)}] [\uPic]^3.
\end{align*}
(In the second step we substituted $l \to k, k \to l+k$.)

\subsubsection*{Type $(2,0)$} The bounds of the degrees are: $$\mu(\cE_\bullet/\cE_\bullet^\prime)<\mu(\cE_\bullet^\prime) \Leftrightarrow \frac{d_0-d_0^\prime+\sigma}{2}<\frac{d_0^\prime}{2} \Leftrightarrow  d_0^\prime > \frac{d_0+\sigma}{2}.$$
Further $\cE_\bullet/\cE_\bullet^\prime$ must be semistable, so we need $0\leq d_0-d_0^\prime \leq \sigma \Leftrightarrow d_0-\sigma \leq d_0^\prime \leq d_0$. Note that $d_0-\sigma \leq \frac{d_0+\sigma}{2}\Leftrightarrow  d_0\leq 3\sigma$ is automatic, so that we end up with $$\frac{d_0+\sigma}{2}< d_0^\prime \leq d_0 \textrm{ and } \mu_{\min}(\cE_\bullet^\prime)>\frac{d_0-d_0^\prime+\sigma}{2}.$$
We have
\begin{align*} 
\dim \uExt(\cE_\bullet/\cE_\bullet^\prime,\cE_\bullet^\prime) &=2(g-1)+2(d_0-d_0^\prime)-d_0^\prime-(2(g-1)-d_0^\prime)
\\ &=2d_0-2d_0^\prime.
\end{align*}
Thus the sum over the strata contributes ($l:=d_0-d_0^\prime$):
\begin{align*}
&\sum_{l=0}^{\lceil \frac{d_0-\sigma}{2}\rceil-1} \bL^{2l} [\uPic][C^{(l)}] [\Bun_2^{d_0-l,\mu_{\min}>\frac{l+\sigma}{2}}]\\
&\stackrel{(\ref{Bun2})}= [\uPic][\Bun_2] \sum_{l=0}^{\lceil \frac{d_0-\sigma}{2}\rceil-1} \bL^{2l}[C^{(l)}]
- [\uPic]^3 \sum_{l=0}^{\lceil \frac{d_0-\sigma}{2}\rceil-1} \sum_{k=\lfloor d_0-\frac{3l+\sigma}{2}\rfloor+1}^\infty [C^{(l)}]\bL^{d_0+l-2k+g-1}.
\end{align*}

\subsubsection*{Type $(1,1)$} The bounds on the degree are: $$\mu(\cE_\bullet/\cE_\bullet^\prime)<\mu(\cE_\bullet^\prime) \Leftrightarrow \frac{d_0-d_0^\prime}{2}< \frac{d_0^\prime+\sigma}{2} \Leftrightarrow d_0^\prime > \frac{d_0-\sigma}{2}.$$
Furthermore we need $\mu_{\min}(\cE_\bullet^\prime)>\mu(\cE_\bullet/\cE_\bullet^\prime) \Leftrightarrow \min\{\frac{d_0^\prime+\sigma}{2},\sigma\}>\frac{d_0-d_0^\prime}{2}$. The condition $\sigma > \frac{d_0-d_0^\prime}{2}$ is automatic, because $d_0-2\sigma-\frac{d_0-\sigma}{2}=\frac{d_0-3\sigma}{2}\leq 0$.
Thus the only condition is $$d_0^\prime > \frac{d_0-\sigma}{2}.$$

The space of extensions is of dimension $2g-2+d_0-3d_0^\prime$. Thus the contribution of the strata is
\begin{align*}
&\sum_{l=\lfloor \frac{d_0-\sigma}{2}\rfloor +1}^\infty \bL^{2g-2+d_0-3l}[\uPic][C^{(l)}][\Bun_{2}^{d_0-l,\sst}]\\
&\stackrel{(\ref{Bun2})}=\bL^{2g-2+d_0}[\uPic][\Bun_2] \sum_{l=\lfloor \frac{d_0-\sigma}{2}\rfloor +1}^\infty \bL^{-3l}[C^{(l)}] \\
&- \sum_{l=\lfloor \frac{d_0-\sigma}{2}\rfloor +1}^\infty\sum_{k=\lfloor\frac{d_0-l}{2}\rfloor+1}^\infty [\uPic]^3[C^{(l)}] \bL^{3g-3+2d_0-4l-2k}.\\
\end{align*}

\subsubsection*{Type $(1,0)$} The bounds on the degree are: $$\mu(\cE_\bullet/\cE_\bullet^\prime)<\mu(\cE_\bullet^\prime) \Leftrightarrow d_0^\prime > \frac{d_0-d_0^\prime+\sigma}{3} \Leftrightarrow d_0^\prime > \frac{d_0+\sigma}{4}.$$
Also $\cE_\bullet/\cE_\bullet^\prime$ must be semistable. By Example \ref{Rank21} the space of these semistable triples is non-empty only if  $\frac{\sigma}{2}\leq d_0-d_0^\prime\leq 2\sigma$ or equivalently $d_0-2\sigma \leq d_0^\prime \leq d_0-\frac{\sigma}{2}$. However $d_0-2\sigma> \frac{d_0+\sigma}{4} \Leftrightarrow d_0> 3\sigma$, which we already discarded. Thus we find the condition $$\frac{d_0+\sigma}{4}< d_0^\prime \leq d_0-\frac{\sigma}{2}.$$

Finally $\dim \uExt(\cE_\bullet/\cE_\bullet^\prime,\cE_\bullet^\prime)= 2(g-1)+d_0-3d_0^\prime-(g-1-d_0^\prime)=g-1+d_0-2d_0^\prime$.
Thus the contribution of the strata is
\begin{align*}
& \bL^{(g-1)+d_0}[\uPic] \sum_{l=\lfloor \frac{d_0+\sigma}{4}\rfloor+1}^{\lceil d_0- \frac{\sigma}{2}\rceil-1} \bL^{-2l} [\cM(2,1)_{d_0-l,0}^{\alpha-\sst}]\\
&= \bL^{(g-1)+d_0}[\uPic] \sum_{l=\lfloor \frac{d_0+\sigma}{4}\rfloor+1}^{\lceil d_0- \frac{\sigma}{2}\rceil-1} \bL^{-2l} [\cM(2,1)_{d_0-l,0}]\\
&- \sum_{l=\lfloor \frac{d_0+\sigma}{4}\rfloor+1}^{\lceil d_0-\frac{\sigma}{2}\rceil-1} [\uPic]^3 \sum_{k=\lceil\frac{2d_0-2l-\sigma}{3}\rceil}^\infty \bL^{2g-2+2d_0-3l-2k} [C^{(k)}] \\ 
&-\sum_{l=\lfloor \frac{d_0+\sigma}{4}\rfloor+1}^{\lceil d_0- \frac{\sigma}{2}\rceil-1} [\uPic]^3 \sum_{k=0}^{\lceil \frac{2d_0-2l-\sigma}{3}\rceil-1}  \bL^{g-1+d_0-2l+k} [C^{(k)}]\\
&= \bL^{2d_0-(g-1)}[\uPic]^2[\Bun_2] \sum_{l=\lfloor \frac{d_0+\sigma}{4}\rfloor+1}^{\lceil d_0-\frac{\sigma}{2}\rceil-1} \bL^{-3l}\\ &- \sum_{k=0}^\infty \sum_{l=\lfloor \max\{\frac{d_0+\sigma}{4},d_0-\frac{3k+\sigma}{2}\}\rfloor+1}^{d_0-\frac{\sigma}{2}} [\uPic]^3 \bL^{2g-2+2d_0-3l-2k} [C^{(k)}] \\
&- \sum_{k=0}^{\lceil \frac{d_0-\sigma}{2}\rceil-1} \sum_{l\lfloor\frac{d_0+\sigma}{4}\rfloor+1}^{\lceil d_0- \frac{3k+\sigma}{2}\rceil-1}[\uPic]^3 \bL^{g-1+d_0-2l+k} [C^{(k)}].
\end{align*}
%\end{itemize}
Adding up the above terms we find the claimed formula. 
\end{proof}

\begin{example}[Chains of rank $(2,1,1)$]\label{m211}
Assume for simplicity that the stability parameter is of the form $\alpha=(0,\sigma,2\sigma)$, that $\alpha$ is good (e.g., $\sigma$ is irrational) with $\sigma\geq 2g-2$ and $\un{d}$ is such that $\mu(\cE_\bullet) \not\in\bZ$. 
Then $\cM(2,1,1)_{\un{d}}^{\alpha-\sst}$ is empty unless: $$d_2\leq d_1, d_0<2\mu(\cE_\bullet), d_0+d_1< 3\mu(\cE_\bullet)-\sigma \textrm{ and } d_2+2d_1< 3\mu(\cE_\bullet)-3\sigma.$$
If these inequalities hold then
$$[\cM(2,1,1)_{\un{d}}^{\alpha-\sst}] = [\uPic]^2[C^{(d_1-d_2)}]\sum_{l=0}^{\lceil d_0-d_1-\mu(\cE_\bullet)\rceil-1} (\bL^{g-1+d_0-2d_1-2l} - \bL^{l}) [C^{(l)}].$$
\begin{proof}
%(1) Follows immediately from Proposition \ref{ChainStrata}. %$\cM(2,1,1)_{\un{d},\phi_i\neq 0}\cong \cM(2,1)_{d_0,d_1} \times C^{(d_1-d_2)}$ and Lemma \ref{M21}. 
Any chain $\cE_\bullet\in \rcM(2,1,1)$ has subchains $(0\to 0 \to \cE_0)$ which is of slope $\frac{d_0}{2}$, $(0\to\cE_1\to \cE_0)$ of slope $\frac{d_0+d_1+\sigma}{3}$ and $(\cE_2\to \cE_1\to \cE_1)$ of slope $\frac{d_2+2d_1+3\sigma}{3}$. This proves the claimed inequalities.

If the inequalities hold, we already excluded destabilizing subchains of rank $(2,1,0)$, $(2,0,0)$ and $(1,1,0)$.

If $\cE^\prime_\bullet=\cE^{h-1}_\bullet$ is of rank $(1,1,1)$ we have 
$$\mu(\cE_\bullet/\cE_\bullet^\prime)=d_0-d_0^\prime<\mu(\cE_\bullet) \Leftrightarrow d_0^\prime-d_1 > d_0-d_1-\mu(\cE_\bullet) \textrm{ and } \mu_{\min}(\cE_\bullet^\prime)>\mu(\cE_\bullet/\cE_\bullet^\prime).$$ The last inequality is automatic, since the quotients of $\cE^\prime_\bullet$ are $(\cE_2 \to 0 \to 0)$ and $(\cE_2\to \cE_1\to 0)$ which are of slope $>\mu(\cE_\bullet)$.
A Harder--Narasimhan stratum of this type contributes $$\bL^{g-1+d_0-2d_1-2(d_0^\prime-d_1)}[\uPic]^2[C^{(d_1-d_2)}][C^{(d_0^\prime-d_1)}].$$

For $\cE_\bullet^{h-1}$ to be of rank $(1,0,0)$ we need $d_0^\prime > \mu(\cE_\bullet)$ and the quotient $\cE_\bullet/\cE_\bullet^\prime$ has to be semistable. This holds if and only if $d_1\leq d_0-d_0^\prime$, because the quotients $(\cE_2 \to 0 \to 0)$ and $(\cE_2\to \cE_1\to 0)$ are of slope $\geq \mu(\cE_\bullet)>\mu(\cE_\bullet/\cE_\bullet^\prime)$. Thus we find $$\mu(\cE_\bullet) < d_0^\prime \leq d_0-d_1.$$
Finally $\dim \uExt(\cE_\bullet/\cE_\bullet^\prime,\cE_\bullet^\prime)= d_0-d_1-d_0^\prime$. So the class of such a stratum is $$[\uPic]^2[C^{(d_1-d_2)}][C^{(d_0-d_0^\prime-d_1)}]\bL^{d_0-d_1-d_0^\prime}.$$
Thus using \ref{Mm111} we find:
\begin{align*}
[\cM(2,1,1)_{\un{d}}^{\alpha-\sst}] &= [\rcM(2,1,1)_{\un{d}}] - \sum [\textrm{HN-Strata}] \\
&= \bL^{g-1+d_0-2d_1} [\uPic]^2[C^{(d_1-d_2)}]\sum_{l=0}^\infty [C^{(l)}]\bL^{-2l}\\
&- \bL^{g-1+d_0-2d_1} [\uPic]^2[C^{(d_1-d_2)}]\sum_{l>d_0-d_1-\mu(\cE_\bullet)} [C^{(l)}]\bL^{-2l} \\
&- [\uPic]^2 [C^{(d_1-d_2)}] \sum_{l>\mu(\cE_\bullet)}^{d_0-d_1} \bL^{d_0-d_1-l}[C^{(d_0-d_1-l)}]\\
&= [\uPic]^2[C^{(d_1-d_2)}]\sum_{l=0}^{\lceil d_0-d_1-\mu(\cE_\bullet)\rceil-1} (\bL^{g-1+d_0-2d_1-2l} - \bL^{l}) [C^{(l)}].
\end{align*}
\end{proof}
\end{example}

\begin{example}[Rank (1,1,1,1)]\label{M1111} $\cM(1,1,1,1)^{\un{\alpha}-\sst}$ for $\alpha=(0,\sigma,2\sigma,3\sigma)$ is non empty if and only if: 
\begin{enumerate}
\item $d_0\geq d_1 \geq d_2 \geq d_3$
\item $d_0+d_1+d_2\leq 3 d_3 +6\sigma$
\item $d_0+d_1 \leq d_2+d_3 + 4\sigma$
\item $3d_0 \leq d_1+d_2+d_3 + 6\sigma$
\end{enumerate}
If these conditions hold, and $\alpha$ is non-critical, then $$[\cM(1,1,1,1)^{\sigma-\sst}]=[\uPic][C^{(d_0-d_1)}][C^{(d_1-d_2)}][C^{(d_2-d_3)}].$$

\begin{proof}
The first condition is the necessary and sufficient condition for $\cM(1,1,1,1)$ to be non-empty.  The possible subchains of $\cE_\bullet$ are $ 0 \to \cE_2 \to \cE_1 \to \cE_0$, $0\to 0\to \cE_1 \to \cE_0$ and $0 \to \dots \to \cE_0$. These give the other conditions.

The second part follows immediately from the first, since in this case all chains in $\rcM(1,1,1,1)$ are semistable and this stack classifies chains of line bundles with non-zero maps between them.
\end{proof}
\end{example}

%%%%%%%%%%%%%%%%%%%%%%%%%%%%%%%%%%%%%%%%%%%%%%%%%%%%%%%%%%%%%%%%%%%%%%%%%%%%%%%%%%%%%%%%%
\subsection{Chains of rank $(n,\dots,n)$ for $d_{i-1}-d_{i}<\alpha_i-\alpha_{i-1}$}
%%%%%%%%%%%%%%%%%%%%%%%%%%%%%%%%%%%%%%%%%%%%%%%%%%%%%%%%%%%%%%%%%%%%%%%%%%%%%%%%%%%%%%%%%
In this section we give another case where we can obtain an inductive formula for any rank. For our application to Higgs bundles of rank $4$ we only need the case of chains of rank $(2,2)$, but this has a natural generalization for chains of rank $(n,\dots,n)$, which does not require an extra effort. We therefore formulate the result in the more general case.

The following proposition improves  \cite[Proposition 6.4]{BGPG} and also extends the result to chains. 

\begin{proposition}\label{alphanice}
Fix $n,r\in \bN$ and write $\un{n}=(n,\dots,n)$. Fix a degree $\un{d}=(d_0,\dots,d_r)$ with $d_r\leq d_{r-1} \leq \dots \leq d_0$ and $\un{\alpha}$ a stability parameter. Suppose that for all $i>0$ we have $d_{i-1}-d_i < \alpha_{i}-\alpha_{i-1}$.
\begin{enumerate}
\item For any $\alpha$ semistable chain of rank $\un{n}$ and degree $\un{d}$ all maps $\phi_i$ are injective, i.e., $\cM(\un{n})_{\un{d}}^{\un{\alpha}-\sst}\subset \cM(\un{n})_{\un{d}}^{\gensurj}$.
\item Suppose $\cE_\bullet \in \cM(\un{n})_{\un{d}}^{\gensurj}$ is a chain with HN-flag $\cE^1_\bullet \subset \dots \subset \cE_\bullet^h =\cE_\bullet$. Then for any $j$ we have $\rk(\cE^j_{\bullet}/\cE^{j-1}_\bullet)=(m_j,\dots,m_j)$ for some $m_j\in \bN$. 
\end{enumerate}
\end{proposition}
\begin{proof}
To show (1) suppose $\cE_\bullet$ was a semistable chain with $\rk(\Ker(\phi_i))=m <n$ for some $i$. Then $\cK_\bullet =(\cdots 0 \to \ker(\phi_i) \to 0 \cdots)$ and $\cE^\prime_\bullet:=(\cE_r \to \cdots \to \cE_i \to \cE_i/\ker(\phi_i) \to \cE_{i-2} \to \cdots \cE_0)$ are subchains of $\cE_\bullet$. Denote $\deg(\ker(\phi_i)):=k$.
Thus we have 
\begin{align*}
\mu(\cK_\bullet)= \frac{k}{m}+\alpha_i &\leq \mu(\cE_\bullet) = \frac{\sum_{j} d_j + n\sum \alpha_j}{(r+1)n} \\
\Leftrightarrow \quad\quad\quad\quad (r+1)nk &\leq m\sum_{j} d_j + mn\sum \alpha_j -m(r+1)n\alpha_i
\end{align*}
and
\begin{align*}
\mu(\cE_\bullet^\prime)= \frac{\sum_{j\neq i-1} d_j +d_i-k + n\sum \alpha_j -m\alpha_{i-1}}{(r+1)n-m}&\leq \frac{\sum_{j} d_j + n\sum \alpha_j}{(r+1)n} \\
%\Leftrightarrow
%(r+1)n(d_i-k) + (r+1)n(n-m)\alpha_{i-1} + m(\sum_{j} d_j + n \sum \alpha_j)&\leq ((r+1)n-m)d_{i-1} + ((r+1)n-m)\alpha_{i-1}\\
\Leftrightarrow 
(r+1)n(d_i-d_{i-1}-m\alpha_{i-1}) + m\bigg(\sum_j d_j +n\sum_j\alpha_j\bigg) &\leq (r+1)nk.
\end{align*}
This implies
\begin{align*}
(r+1)n(d_i-d_{i-1}-m\alpha_{i-1}) &\leq -m(r+1)n\alpha_i \\
\Leftrightarrow \quad\quad\quad\quad\quad\quad\quad\quad \alpha_i-\alpha_{i-1} & \leq \frac{d_{i-1}-d_i}{m}.
\end{align*}
This contradicts our assumption.

The proof of (2) is by induction. Suppose $\cE_\bullet^\prime=\cE_\bullet^{i}\subset \cE_\bullet$ was a destabilizing subchain such that not all $\cE_k^i$ have equal rank. We will denote by $\cE_\bullet^\pprime := \cE_\bullet/\cE_\bullet^\prime$ the quotient chain. 

By assumption all maps $\phi_i^\prime$ are injective, so that $\rk(\cE_\bullet^\prime)=\un{n}^\prime$ with $n_r^\prime \leq \dots \leq n_0^\prime\leq n$. Let $i$ be the minimal integer such that $n_i^\prime < n_{i-1}^\prime$. 

Then $\cK_\bullet:=(0\to \dots \to \ker(\phi_i^\pprime) \to 0 \to \dots)$ is a subchain of $\cE_\bullet^\pprime$ and $\cQ_\bullet := (\cdots 0\to \cE_{i-1}^\prime/\cE_{i}^\prime  \to 0 \cdots)$ is a quotient of $\cE_\bullet^\prime$. 
Thus we have:
$$\mu(\cQ_\bullet) \geq \mu_{\min}(\cE_\bullet^\prime) > \mu_{\max}(\cE_\bullet^\pprime) \geq \mu(\cK_\bullet),$$
i.e.,
\begin{align*}
\frac{d_{i-1}^\prime-d_i^\prime}{n_{i-1}^\prime-n_i^\prime} + \alpha_{i-1} &> \mu(\cK_\bullet) \geq \frac{d_i^\pprime-d_{i-1}^\pprime}{n_{i}^\pprime -n_{i-1}^\pprime} + \alpha_i\\
\Rightarrow\quad
\frac{d_{i-1}^\prime-d_i^\prime}{n_{i-1}^\prime-n_i^\prime} + \alpha_{i-1} &> \frac{d_i-d_{i-1}+d_{i-1}^\prime-d_{i}^\prime}{n_{i-1}^\prime -n_{i}^\prime} + \alpha_i\\
\Leftrightarrow\quad\quad\quad\quad\;
\frac{d_{i-1}-d_i}{n_{i-1}^\prime-n_i^\prime} &> \alpha_i-\alpha_{i-1},
\end{align*}
which again contradicts our assumption.
\end{proof}
This proposition allows us to deduce the following recursion formula for the motive of $\cM(\un{n})_{\un{d}}^{\alpha-\sst}$:
\begin{corollary}
Let $\un{n}=(n,\dots,n)$ be constant. If $\un{\alpha},\un{d}$ satisfy $d_{i-1}-d_i < \alpha_{i}-\alpha_{i-1}$ for all $i>0$ then we have
\begin{align*} 
[\cM(\un{n})_{\un{d}}^{\alpha-\sst}] &= [\Bun_n^{d_0}] \prod_{i=1}^r [\Sym^{d_{i-1}-d_i}(C\times \bP^{n-1})]\\
                       &-\bigg(\sum_{\un{m},\un{e},\un{k}} \bL^{\sum_{k<j} \chi_{kj}} \prod_i [\cM(\un{m}_j)_{\un{e_j}}^{\alpha-\sst}]\bigg),
\end{align*}
where the sum runs over all partitions $n=\sum_{j=1}^l m_j$, $\un{d}=\sum_j \un{e}^{(j)}$ such that for all $i,j$ we have $e_i^{(j)} \leq e_{i-1}^{(j)}$ and for $\mu^{(j)}:=\frac{\sum_i e_i^{(j)}}{rm_j}$ we have $\mu_1 > \dots > \mu_l$.
We have written $\chi_{kj}= m_jm_k(g-1)+ \sum_{i=0}^r (m_ke^{(j)}_i-m_je^{(k)}_i)- \sum_{i=1}^r (m_ke^{j}_i -m_je^{(k)}_{i-1})$.
\end{corollary}
\begin{proof}
From Proposition \ref{alphanice} we know that under our assumption on $\alpha$ all semi\-stable chains are contained in the substack of chains such that the $\phi_i$ have full rank and moreover for any such chain all subquotients of the HN-filtration also satisfy this condition. 

Thus we have
$$\cM(\un{n})_{\un{d}}^{\alpha-\sst}= \cM(\un{n})_{\un{d}}^{\text{gen-surj}} - \bigcup \textrm{ Harder--Narasimhan strata}.$$

The Harder--Narasimhan strata are given by rank and degrees as claimed. Since in all occurring subquotients the morphisms $\phi$ are injective we can apply Proposition \ref{Ext-Rechnung} to conclude that these strata are (iterated) vector bundle stacks over $\prod_i \cM(\un{m}^{(j)})_{\un{e}^{(j)}}^{\alpha-\sst}$. The dimension of the fibers are: 
\begin{align*}
\dim \uExt(\cE^j_\bullet/\cE^{j-1}_\bullet,\cE^k_\bullet/\cE^{k-1}_\bullet) &= m_jm_k(g-1)+ \sum_{i=0}^r (m_ke^{(j)}_i-m_je^{(k)}_i)\\&- \sum_{i=1}^r (m_ke^{j}_i -m_je^{(k)}_{i-1}).
\end{align*}
This proves the corollary.
\end{proof}
For the computation of rank $4$ Higgs-bundles we will only need the following special case:
\begin{corollary}
Let $\un{d}=(d_0,0)$ with $d_0<\alpha_1-\alpha_0$ and $d_0$ odd. Then:
$$[\cM(2,2)^{\alpha-\sst}_{(d_0,0)}] = [\Bun_2^{d_0}][(C\times \bP^1)^{(d_0)}] - \frac{\bL^{(g-1)+\frac{d_0-1}{2}}}{1-\bL^{-1}} [\uPic]^2\frac{[(2C)^{(d_0)}]}{2}.$$
\end{corollary}
\begin{proof}
A subtriple $\cE_\bullet^\prime\subset \cE_\bullet$ is destabilizing, if and only if it is of rank $(1,1)$ and its degree satisfies:
\begin{enumerate}
\item $\mu(\cE_\bullet^\prime)> \mu(\cE_\bullet) \Leftrightarrow d_0^\prime+d_1^\prime > \frac{d_0}{2}$.
\item $\cE_\bullet^\prime$ is semistable, so $0 \leq d_0^\prime-d_1^\prime \leq \alpha_1-\alpha_0$
\item $\cE_\bullet/\cE_\bullet^\prime$ is semistable so  $d_0-\alpha_1+\alpha_0 \leq  d_0^\prime -d_1^\prime \leq d_0.$
\end{enumerate}
Since $d_0 < \alpha_1-\alpha_0$ we find $0 \leq d_0^\prime-d_1^\prime \leq d_0$. The dimension of $\uExt(\cE_\bullet/\cE_\bullet^\prime,\cE_\bullet^\prime)$ is $(g-1)+d_0-d_0^\prime-d_1^\prime$.
Thus the sum over all HN-strata of rank $(1,1)$ is:
\begin{align*}
[\cM(2,2)_{(1,1)}] &= \sum_{k>\frac{d_0}{2}} \sum_{l=0}^{d_0} \bL^{(g-1)+d_0-k} [\uPic]^2[C^{(l)}][C^{(d_0-l))}]\\
 &= \sum_{k=\frac{d_0+1}{2}} \bL^{(g-1)+d_0-k}[\uPic]^2 \sum_{e=0}^{\frac{d_0-1}{2}} [C^{(e)}][C^{(d_0-e)}] \\
 &= \frac{\bL^{(g-1)+\frac{d_0-1}{2}}}{1-\bL^{-1}} [\uPic]^2\frac{[(2C)^{(d_0)}]}{2}.
\end{align*}
Therefore we find 
$$[\cM(2,2)^{\alpha-\sst}] = [\Bun_2^{d_0}][(C\times \bP^1)^{(d_0)}] - \frac{\bL^{(g-1)+\frac{d_0-1}{2}}}{1-\bL^{-1}} [\uPic]^2\frac{[(2C)^{(d_0)}]}{2},$$
which is the claimed formula.
\end{proof}

\subsection{Chains of rank $(2,2)$ for $d_0-d_1>\alpha_1-\alpha_0$}\label{M22Section}
To complete the computation for the class of the space of Higgs bundles of rank $4$ we are left to compute the class of the space of semistable chains of rank $(2,2)$ in case $d_0-d_1>\alpha_1-\alpha_0$ and the class of the stack of semistable chains of rank $(1,2,1)$. In both situations the stack of all chains of fixed degree does not have a class in $\K0hat$.  In this section we will show how to handle this problem for chains of rank $(2,2)$. The case of rank $(1,2,1)$, which is a bit simpler, will be done in the next section. 

For this section we fix $\alpha=(0,\sigma)$. We will only need to consider chains such that $d_0+d_1$ is odd. Since we may dualize and tensor with line bundles we may therefore assume $d_1=0$ and $d_0>\sigma$ odd. In that case we may also assume that $d_0<2\sigma$ since otherwise for every chain $\cE_\bullet$ the subchain $\cE_0$ is destabilizing.
Also, in order to simplify one formula (type $(1,1)\subset(2,2)$ below) we assume that $\lfloor \sigma \rfloor$ is even, which will be satisfied in our application to Higgs bundles.

First let us explicitly compute the stratification by generic rank given in Proposition \ref{ChainStrata} in the case of rank $(2,2)$. Let us write $\cM(2,2)^{1,k,l}_{(d_0,0)}$ for the space of chains $\cE_1\map{\phi}\cE_0$ such that $\rk(\phi)=1$, $\deg(\ker(\phi))=k$, $\deg(\phi(\cE_1)^\sat)=l$. 
\begin{lemma}\label{M22class}
$$[\cM(2,2)^{1,k,l}_{(d_0,0)}]= [\uPic]^3 [C^{(l+k)}] \bL^{2(g-1)+d_0-2k-2l}.$$ 
In particular for fixed $d_0$ this class only depends on $k+l$. It is non-zero if $k+l\geq 0$.
\end{lemma}

This implies that the sum over all possible $k,l$ does not converge, so that $\cM(2,2)$ does not define a class in $\K0hat$. However, if either $k$ or $l$ are large, then we see that all triples in $\cM(2,2)_{1,k,l}$ will be unstable, since either $\ker(\phi)\to 0$ or $\cE_1 \to \im(\phi)^\sat$ will be a destabilizing subchain.
More precisely $\ker(\phi) \to 0$ is destabilizing if $k+\sigma > \mu(\cE_\bullet)$ and $\cE_1 \to \im(\phi)^\sat$ is destabilizing if $\mu(\cE_\bullet) > d_0-l$.

We will therefore define the following open substack of $\cM(2,2)$:
$$\cM(2,2)^{\fin}_{\un{d}} := \cM(2,2)_{\un{d}}^{\text{gen-surj}} \cup \bigcup_{(k,l) : 0\leq k+l \atop  k+\sigma < \mu(\cE_\bullet) < d_0-l} \cM(2,2)_{\un{d}}^{1,k,l}.$$
From Lemma \ref{M22class} we see that in this stack the class of $\cM(2,2)^{1,k,l}$ occurs at most for $0\leq k+l < d_0- \sigma$. For fixed value of $k+l=m$ in this range there are $\lfloor \frac{d_0}{4}-\frac{\sigma}{2}\rfloor +\lfloor \frac{3d_0}{4}-\frac{\sigma}{2}\rfloor -m+1$ such strata. 
Thus this stack does have a well-defined class in $\K0hat$ which by Lemma \ref{M22class} and Proposition \ref{Mnn-rkn} is:
\begin{align*}
&[\cM(2,2)^{\fin}_{\un{d}}] = \\&[\Bun_2^{d_0}][(C\times \bP^1)^{(d_0)}] + [\uPic]^3\bL^{2(g-1)+d_0} \sum_{m=0}^{\lceil d-\sigma\rceil-1}(d-\lfloor\sigma\rfloor-m) \bL^{-2m}[C^{(m)}].
\end{align*}
Finally denote by $\cM(2,2)^{\out}_{\un{d}} := \cM(2,2)_{\un{d}} - \cM(2,2)^{\fin}_{\un{d}}$. This is the substack of triples such that either $\ker(\phi) \to 0$ or $\cE_0 \to \im(\phi)^\sat$ is a destabilizing subtriple.

Next we compute the Harder--Narasimhan strata such that the HN-flag does not contain triples of rank $(0,1)$ or $(1,2)$, since only these can intersect $\cM(2,2)^{\fin}_{\un{d}}$. Also destabilizing subtriples of rank $(2,0)$ cannot occur by our assumption $d_0<2\sigma$. We will denote the HN-flags by $\dots \subset\cE^\pprime_\bullet \subset \cE^\prime_\bullet \subset \cE_\bullet$ and $d_i^\prime=\deg(\cE_i^\prime)$ etc.

As before, we group the strata according to the rank of the HN-flag. For each rank we will first compute the bounds on the degrees given by the characterizing property of the HN-flag. Then we compute the dimension of the $\uExt$-space from Proposition \ref{HN-Ext-Stack}. Finally we compute the intersection of the stratum with $\cM(2,2)_{\un{d}}^{\out}$.

%%%%%%%%%%%%%%%%%%%%%%%%%%%%%%%%%%%%%%%
%\begin{itemize}
%%%%%%%%%%%%%%%%%%%%%%%%%%%%%%%%%%%%%%
\subsubsection*{Type $(2,1)\subset(2,2)$} The bounds on the degree $\un{d}^\prime$ are:
$$\mu(\cE_\bullet/\cE_\bullet^\prime)<\mu(\cE_\bullet) \Leftrightarrow -d_1^\prime +\sigma <\frac{d_0}{4}+\frac{\sigma}{2} \Leftrightarrow d_1^\prime > \frac{\sigma}{2}-\frac{d_0}{4}.$$
In order to have $\cM(2,1)_{(d_0,d_1^\prime)}^{\alpha-\sst}\neq \emptyset$ we need (Example \ref{Rank21})  $d_1^\prime<\frac{d_0}{2}-\frac{\sigma}{4}$. 
Thus we found the bounds:
$$  \frac{\sigma}{2}-\frac{d_0}{4} < d_1^\prime < \frac{d_0}{2}-\frac{\sigma}{4}.$$
%In particular this implies $\frac{\sigma}{2}-\frac{d_0}{4} < \frac{d_0}{2} -\frac{\sigma}{4}$, i.e., $\sigma<d_0$.
By Proposition \ref{HN-Ext-Stack} we have $\dim\uExt(\cE_\bullet/\cE_\bullet^\prime,\cE_\bullet^\prime)= d_0-(g-1).$ 

We claim that strata of this type are contained in $\cM(2,2)^{\fin}_{\un{d}}$: Since $\cE^\prime_\bullet$ is semistable the morphism $\phi^\prime$ is not zero, so that $\ker(\phi)\to 0$ injects into $\cE_\bullet/\cE_\bullet^\prime$. Thus $\ker(\phi)$ cannot be destabilizing. Also if $\rank(\phi)=1$ then $0\to \cE_0/\im(\phi)$  is a quotient of $\cE_\bullet^\prime$, so it is of slope $\geq \mu(\cE_\bullet^\prime)>\mu(\cE_\bullet)$, so it cannot be a destabilizing quotient. Thus these strata contribute:
$$\sum_{d_1^\prime=\lfloor \frac{\sigma}{2}-\frac{d_0}{4}\rfloor+1}^{\lceil\frac{d_0}{2}-\frac{\sigma}{4}\rceil-1} \bL^{d_0-(g-1)} [\uPic][\cM(2,1)_{d_0,d_1^\prime}^{\alpha-\sst}].$$

%%%%%%%%%%%%%%%%%%%%%%%%%%%%%%%%%%%%%%%%
\subsubsection*{Type $(1,1)\subset (2,1) \subset(2,2)$} The bounds on the degrees are: $$\mu(\cE_\bullet^\pprime)>\mu(\cE_\bullet^\prime/\cE_\bullet^\pprime) \Leftrightarrow \frac{d_0^\pprime+d_1^\prime+\sigma}{2}>d_0-d_0^\pprime \Leftrightarrow d_0^\pprime-d_1^\prime > \frac{2d_0-4d_1^\prime-\sigma}{3}$$ and
$$\mu(\cE_\bullet^\prime/\cE_\bullet^\pprime) > \mu(\cE_\bullet/\cE_\bullet^\prime) \Leftrightarrow d_0-d_0^\pprime > -d_1^\prime+\sigma \Leftrightarrow d_0^\pprime-d_1^\prime< d_0-\sigma.$$
%\item[$\bullet$] 
%\end{itemize}
Thus we find $$\frac{2d_0-4d_1^\prime-\sigma}{3} < d_0^\pprime-d_1^\prime < d_0-\sigma,$$
and in particular $d_1^\prime > \frac{\sigma}{2} - \frac{d_0}{4}$.

Semistability of $\cE_\bullet^\pprime$ implies furthermore $0 \leq d_0^\pprime - d_1^\prime \leq \sigma$.
%Note that $d_1^\prime>\frac{2d_0-d_1^\prime-\sigma}{3}$ is equivalent to $d_1> \frac{d_0}{2}-\frac{\sigma}{4}$.
%In particular this cannot happen if $d_0 < \sigma$. 
Since $d_0<2\sigma$ the right hand inequality follows from the first set of inequalities. 
%Thus we find:
%$$\max(\frac{2d_0-d_1^\prime-\sigma_+}{3},d_1^\prime) \leq d_0^\prime < d_0+d_1^\prime-\sigma.$$
Also $\frac{2d_0-4d_1^\prime-\sigma}{3} < 0 \Leftrightarrow \frac{d_0}{2}-\frac{\sigma}{4} < d_1^\prime$.

We have $\dim \uExt(\cE_\bullet^\prime/\cE_\bullet^\pprime,\cE_\bullet^\pprime) = (g-1)+d_0-2d_0^\pprime$ and as in the previous stratum $\dim\uExt(\cE_\bullet/\cE_\bullet^\prime,\cE_\bullet^\prime)= d_0-(g-1).$
Thus the contribution of the strata is:
\begin{align*}
&\biggl(\sum_{d_1^\prime = \lfloor \frac{\sigma}{2}-\frac{d_0}{4}\rfloor+1}^{\lfloor \frac{d_0}{2}-\frac{\sigma}{4}\rfloor} \sum_{l = \lfloor \frac{2d_0-4d_1^\prime-\sigma}{3}\rfloor +1}^{\lfloor d_0-\sigma\rfloor} \bL^{2d_0 - 2(l+d_1^\prime)} [C^{(l)}]\\
&+ \sum_{d_1^\prime =\lfloor \frac{d_0}{2}-\frac{\sigma}{4}\rfloor +1}^\infty \sum_{l=0}^{\lfloor d_0-\sigma\rfloor}  \bL^{2d_0 - 2(l+d_1^\prime)} [C^{(l)}]\biggr) [\uPic]^3.\\
%&= \biggl(\sum_{d_1^\prime > \frac{\sigma}{2}-\frac{d_0}{4}}^{\frac{d_0}{2}-\frac{\sigma}{4}} \sum_{l > \frac{2d_0-4d_1^\prime-\sigma}{3}}^{d_0-\sigma} \bL^{2d_0 - 2(l+d_1^\prime)} [C^{(l)}] + \frac{\bL^{\lceil \frac{d_0}{2}-\frac{\sigma}{4} \rceil+2d_0}}{1-\bL^{-2}} \sum_{l=0}^{d_0-\sigma} \bL^{-2l}[C^{(l)}]\biggr) [\uPic]^3.
\end{align*}

Finally we compute the intersection of these strata with $\cM(2,2)^{\out}_{\un{d}}$. First, note that $\ker(\phi)\to 0$ would inject into $\cE_\bullet/\cE_\bullet^\prime$ so this cannot be destabilizing. 

If $\rank(\phi)=1$ then $\im(\phi)^\sat=\cE_0^\pprime$. Thus $\cE_1 \to \cE_0^\pprime$ is a destabilizing chain if $d_0^\pprime > d_0-\mu(\cE_\bullet) = \frac{3}{4}d_0-\frac{\sigma}{2}$. Equivalently $d_0^\pprime-d_1^\prime>\frac{\sigma}{2}-\frac{d_0}{4}+d_0-\sigma -d_1^\prime$.
Since $d_1^\prime > \frac{\sigma}{2}-\frac{d_0}{4}$ this implies $d_0^\pprime-d_1^\prime > \frac{2d_0-4d_1^\prime-\sigma}{3}$. Thus we find the conditions:  
$$\max\Big\{\frac{3d_0}{4}-\frac{\sigma}{2} -d_1^\prime,0\Big\} \leq d_0^\pprime-d_1^\prime  < d_0-\sigma.$$

The class of the intersection of the stratum with $\cM(2,2)^{\out}_{\un{d}}$ is the stack of extensions in $\uExt(\cE_\bullet/\cE_\bullet^\prime,\cE_\bullet^\prime/\cE_\bullet^\pprime,\cE_\bullet^\pprime)$ that satisfy $\phi(\cE_1)\subset \cE_0^\pprime$, which is equivalent to $\uExt(\cE_\bullet/\cE_\bullet^\prime,\cE_\bullet^\pprime)\times \uExt(\cE_0/\cE_0^\pprime,\cE_0^\pprime)$. The dimension of the first $\Ext$-stack is $(g-1)-2d_1^\prime - ((g-1)-d_1^\prime-d_0^\pprime)=d_0^\pprime -d_1^\prime$, the dimension of the second is $g-1+d_0-2d_0^\pprime$. 

Thus the intersection of the HN-strata with $\cM(2,2)^{\out}_{\un{d}}$ is 
%$\bL^{g-1+d_0-d_0^\prime -3d_1^\prime} [\uPic]^3 [C^{(d_0^\prime-d_1^\prime)}].$

%In total we find that the class of the intersection of the strata with $\cM(2,2)_{out}$ is (we set $l=d_0^\prime-d_1^\prime$)
\begin{align*}
&\sum_{d_1^\prime =\lfloor \frac{\sigma}{2}-\frac{d_0}{4}\rfloor+1}^{\lfloor \frac{3d_0}{4}-\frac{\sigma}{2}\rfloor} \sum_{l=\lfloor \frac{3d_0}{4}-\frac{\sigma}{2}\rfloor+1-d_1}^{\lfloor d_0-\sigma\rfloor} \bL^{-l-2d_1^\prime+d_0+(g-1)} [C^{(l)}] [\uPic]^3\\
&+\sum_{d_1^\prime =\lfloor \frac{3d_0}{4}-\frac{\sigma}{2}\rfloor+1}^\infty \sum_{l=0}^{\lfloor d_0-\sigma\rfloor}\bL^{-l-2d_1^\prime+d_0+(g-1)}[C^{(l)}] [\uPic]^3.
\end{align*}

%\begin{itemize}
%%%%%%%%%%%%%%%%%%%%%%%%%%%%%%%%%%%%%%%%%%%%%%%%%%%%%%
\subsubsection*{Type $(1,0)\subset (1,1)\subset (2,1) \subset (2,2)$} We claim that this cannot occur, because we need 
$\mu(\cE^{3}_\bullet)>\mu(\cE^\pprime_\bullet/\cE^{3}_\bullet) \Leftrightarrow  d_1^\prime+\sigma<d_0^\pprime$ and 
$\mu(\cE^\prime_\bullet/\cE_\bullet^\pprime)>\mu(\cE_\bullet/\cE^\prime_\bullet) \Leftrightarrow d_0-d_0^\pprime>-d_1^\prime+\sigma \Leftrightarrow d_0^\pprime < d_0+d_1^\prime-\sigma$.
This implies $d_1^\prime +\sigma < d_0^\prime < d_0 +d_1^\prime -\sigma \Rightarrow 2\sigma < d_0$, contradicting our assumption on $d_0<2\sigma$.

%%%%%%%%%%%%%%%%%%%%%%%%%%%%%%%%%%%%%%%%%%
\subsubsection*{Type $(1,0) \subset (2,1) \subset (2,2)$}
The bounds on the degrees are: 
$\mu(\cE^\pprime_\bullet)>\mu(\cE_\bullet^\prime/\cE^\pprime_\bullet) \Leftrightarrow d_0^\pprime > \frac{d_0-d_0^\pprime+d_1^\prime+\sigma}{2} \Leftrightarrow d_0^\pprime> \frac{d_0+d_1^\prime+\sigma}{3}$
and
%\begin{enumerate}
%\item 
$\mu(\cE^\prime_\bullet/\cE_\bullet^\pprime)>\mu(\cE_\bullet/\cE^\prime_\bullet) \Leftrightarrow 
\frac{d_1^\prime+d_0-d_0^\pprime+\sigma}{2} > -d_1^\prime+\sigma \Leftrightarrow d_0^\pprime < 3d_1^\prime +d_0-\sigma$.
Thus we find $$\frac{d_0+d_1^\prime+\sigma}{3}< d_0^\pprime < 3d_1^\prime +d_0-\sigma, \text{ in particular }\frac{\sigma}{2}-\frac{d_0}{4}<d_1^\prime.$$

The quotient $\cE^\prime_\bullet/\cE^\pprime_\bullet$ has to be semistable, i.e., $0\leq d_0-d_0^\pprime -d_1^\prime \leq \sigma$ equivalently $d_0-d_1^\prime-\sigma \leq d_0^\pprime \leq d_0-d_1^\prime$.
In particular we need $\frac{d_0+d_1^\prime+\sigma}{3}< d_0-d_1^\prime \Leftrightarrow d_1^\prime < \frac{d_0}{2}-\frac{\sigma}{4}$. Thus we find $ \frac{\sigma}{2}-\frac{d_0}{4} < d_1^\prime < \frac{d_0}{2}-\frac{\sigma}{4}$.

Also $d_0-d_1^\prime-\sigma> \frac{d_0+d_1^\prime+\sigma}{3} \Rightarrow  2d_0-4\sigma > 4d_1^\prime > 2\sigma-d_0 \Rightarrow d_0 > 2\sigma$ which we excluded. Finally $3d_1^\prime +d_0-\sigma>d_0-d_1^\prime \Leftrightarrow d_1^\prime > \frac{\sigma}{4}$. Thus we find 
$$\frac{d_0+d_1^\prime+\sigma}{3} < d_0^\pprime \left\{\begin{array}{ll} < 3d_1^\prime +d_0-\sigma & \text{if } d_1\leq \frac{\sigma}{4}\\ \leq d_0-d_1^\prime & \text{if } d_1^\prime > \frac{\sigma}{4}.\end{array}\right.$$  

We have $\dim\uExt(\cE_\bullet^\prime/\cE_\bullet^\pprime,\cE_\bullet^\pprime)= d_0-d_0^\pprime-d_1^\prime$ and $\dim\uExt(\cE_\bullet/\cE_\bullet^\prime,\cE_\bullet^\prime)= d_0-(g-1)$.
%The class of such a stratum is given by:
%$$ \bL^{d_0-(g-1)}\bL^{d_0-d_0^\pprime-d_1^\prime}[\uPic]^3 [C^{(d_0-d_0^\pprime-d_1^\prime)}].)$$
Strata of this type are contained in $\cM(2,2)^{\fin}$, because $\ker(\phi)\to 0$ injects into $\cE_\bullet/\cE^\prime_\bullet$, so this cannot be destabilizing. Also if $\rank(\phi)=1$ then $\cE_0^\pprime$ has to inject into $\cE_0/\im(\phi)$ because $\cE_1^\prime \hookrightarrow \cE_0/\cE_0^\pprime$. So this cannot be a destabilizing quotient.

Thus the contribution of these strata is ($l=d_0-d_0^\pprime-d_1^\prime$):
\begin{align*}
& \biggl(\sum_{d_1^\prime =\lfloor \frac{\sigma}{2}-\frac{d_0}{4}\rfloor +1}^{\lfloor \frac{\sigma}{4}\rfloor} \sum_{l=\lfloor\sigma\rfloor+1-4d_1^\prime}^{\lceil\frac{2d_0-4d_1^\prime-\sigma}{3}\rceil-1} \bL^{l} [C^{(l)}]\\
&+ \sum_{d_1^\prime =\lfloor \frac{\sigma}{4}\rfloor+1}^{\lfloor\frac{d_0}{2}-\frac{\sigma}{4}\rfloor} \sum_{l=0}^{\lceil\frac{2d_0-4d_1^\prime-\sigma}{3}\rceil-1} \bL^{l} [C^{(l)}]\biggr) \bL^{d_0-(g-1)}[\uPic]^3.
\end{align*}

%%%%%%%%%%%%%%%%%%%%%%%%%%%%%%%
\subsubsection*{Type $(1,1) \subset (2,2)$} The bounds on the degree are:
%\begin{enumerate}
%\item 
$\mu(\cE_\bullet^\prime)> \mu(\cE_\bullet) \Leftrightarrow d_0^\prime+d_1^\prime > \frac{d_0}{2}$.
Moreover $\cE_\bullet^\prime$ and $\cE_\bullet/\cE^\prime_\bullet$ are semistable, so $0 \leq d_0^\prime-d_1^\prime \leq \sigma$ and $0 \leq d_0-d_0^\prime+d_1^\prime \leq \sigma \Leftrightarrow d_0-\sigma \leq  d_0^\prime -d_1^\prime \leq d_0.$
Since $d_0> \sigma$ we find the bounds $$d_0-\sigma \leq d_0^\prime-d_1^\prime \leq \sigma,d_0^\prime+d_1^\prime > \frac{d_0}{2}.$$

We have $\dim\uExt(\cE_\bullet/\cE_\bullet^\prime, \cE_\bullet^\prime)= (g-1)+d_0-d_0^\prime-d_1^\prime$.
Strata of this type automatically satisfy $\rank(\phi)=2$, so they are contained in $\cM(2,2)^{\fin}$. Thus the sum over all these strata is:
\begin{align*}
& \sum_{d_1^\prime+d_0^\prime=\lfloor \frac{d_0}{2}\rfloor+1}^\infty \bL^{-(d_0^\prime+d_1^\prime)+(g-1)+d_0}\sum_{d_0^\prime-d_1^\prime= d_0-\lfloor \sigma\rfloor\atop d_0^\prime-d_1^\prime \equiv d_0^\prime+d_1^\prime \md 2}^{\lfloor \sigma\rfloor} [C^{(d_0^\prime-d_1^\prime)}][C^{(d_0-(d_0^\prime-d_1^\prime))}][\uPic]^2 \\
&=\sum_{k>\frac{d_0}{2}} \bL^{-k+d_0+(g-1)}[\uPic]^2 \sum_{l=0}^{\frac{2\lfloor \sigma\rfloor-d_0-1}{2}} [C^{(d_0+1-\lfloor \sigma\rfloor+2l)}][C^{\lfloor \sigma\rfloor-1 -2l}].%\\
%&=\frac{\bL^{\lceil\frac{d_0}{2}\rceil+d_0+(g-1)}}{1-\bL^{-1}}[\uPic]^2 \sum_{l=d_0-\sigma}^{\sigma-\frac{d_0}{2}} [C^{(l)}][C^{(d_0-l)}].
\end{align*}
Here the equality uses that $\lfloor \sigma \rfloor$ is even and that $d_0$ is odd, so that one can simplify the congruence condition in the summation.

%%%%%%%%%%%%%%%%%%%%%%%%%%%%%%%%%%%%%%%%%%
\subsubsection*{Type $(1,0) \subset (1,1) \subset (2,2)$}
The bounds on the degrees are: %\begin{enumerate}
%\item  $\mu(\cE^\pprime_\bullet)\mu(\cE_\bullet) \Leftrightarrow d_0^\prime > \frac{d_0}{4}+\frac{\sigma}{2}$
$\mu(\cE^\pprime_\bullet)>\mu(\cE^\prime_\bullet/\cE^\pprime_\bullet) \Leftrightarrow d_0^\prime > d_1^\prime +\sigma$ and $\mu(\cE^\prime_\bullet/\cE^\pprime_\bullet)>\mu(\cE_\bullet/\cE^\prime_\bullet) \Leftrightarrow d_1^\prime+\sigma > \frac{d_0-d_0^\prime-d_1^\prime+\sigma}{2} \Leftrightarrow d_1^\prime > \frac{d_0-d_0^\prime-\sigma}{3}$.
So we find $$ \frac{d_0-d_0^\prime-\sigma}{3} < d_1^\prime < d_0^\prime -\sigma, \text{ in particular } \frac{d_0}{4}+\frac{\sigma}{2} < d_0^\prime.$$
Also $\cE_\bullet/\cE_\bullet^\prime$ is semistable, thus: $0 \leq d_0-d_0^\prime+d_1^\prime \leq \sigma \Leftrightarrow d_0^\prime-d_0\leq d_1^\prime \leq d_0^\prime -d_0 +\sigma$. Since $d_0<2\sigma$ the right hand inequality is automatically satisfied. Also $d_0^\prime-d_0 > \frac{d_0-d_0^\prime-\sigma}{3} \Leftrightarrow d_0^\prime > d_0-\frac{\sigma}{4}$.
 % and the right hand inequality can only be satisfied if $d_0>\sigma$.
%\end{enumerate}
Thus we find: 
$$ d_0^\prime - \sigma >  d_1^\prime  \left\{\begin{array}{ll} > \frac{d_0-d_0^\prime-\sigma}{3}&  d_0^\prime \leq d_0-\frac{\sigma}{4}\\ \geq d_0^\prime -d_0& \text{if } d_0^\prime > d_0-\frac{\sigma}{4}.\end{array} \right.$$

We have $\dim(\uExt( \cE_\bullet/\cE_\bullet^\pprime,\cE_\bullet^\pprime))= d_0-(g-1)$ and $\dim (\uExt(\cE/\cE^\prime,\cE^\prime/\cE^\pprime))= g-1-2d_1^\prime$.
Thus 
%the class of any such stratum is $\bL^{d_0-2d_1^\prime}[\uPic]^3[C^{(d_0-d_0^\prime+d_1^\prime)}].$
the sum over the strata is (set $l=d_0-d_0^\prime+d_1^\prime$):
\begin{align*}
&\biggl(\sum_{d_0^\prime=\lfloor\frac{d_0}{4}+\frac{\sigma}{2}\rfloor+1}^{\lfloor d_0-\frac{\sigma}{4}\rfloor} \sum_{l=\lfloor \frac{d_0-d_0^\prime-\sigma}{3}\rfloor+1+d_0-d_0^\prime}^{\lceil d_0-\sigma\rceil} \bL^{-2d_0^\prime-2l} [C^{(l)}]\\
&+ \sum_{d_0^\prime =\lfloor d_0-\frac{\sigma}{4}\rfloor +1} \sum_{l=0}^{\lfloor d_0-\sigma\rfloor} \bL^{-2d_0^\prime-2l} [C^{(l)}]\biggr) \bL^{3d_0}[\uPic]^3.
\end{align*}
Finally we compute the intersection with $\cM(2,2)^{\out}_{\un{d}}$. (We already did the dual case above). 
Since $\cE_\bullet/\cE_\bullet^\prime$ is semistable, the morphism of this chain is non trivial. Therefore if  $\rank(\phi)=1$ then $\cE_0^\prime\neq \im(\phi)^\sat$, so that $\cE^\prime_0$ injects into $\cE_0/\im(\phi)^\sat$. Thus $\cE_0/\im(\phi)^\sat$ cannot be a destabilizing quotient. 

Also if $\phi^\prime=0$ then $(\cE_1^\prime \to 0)=(\ker(\phi) \to 0)$ is destabilizing if and only if $d_1^\prime +\sigma > \frac{d_0}{4} +\frac{\sigma}{2} \Leftrightarrow d_1^\prime > \frac{d_0}{4} -\frac{\sigma}{2}$. 
Since $\frac{d_0-d_0^\prime-\sigma}{3} < \frac{d_0}{4}-\frac{\sigma}{2} \Leftrightarrow \frac{d_0}{4}+\frac{\sigma}{2} < d_0^\prime$, we find the conditions for $d_1^\prime$  
$$ d_0^\prime - \sigma >  d_1^\prime  \left\{\begin{array}{ll} > \frac{d_0}{4}-\frac{\sigma}{2}& \text{if } d_0^\prime \leq \frac{d_0}{4}+\frac{\sigma}{2}+(d_0-\sigma)\\ \geq d_0^\prime -d_0& \text{if } d_0^\prime > \frac{d_0}{4}+\frac{\sigma}{2}+(d_0-\sigma).\end{array} \right.$$

If these conditions are satisfied, chains in the intersection of the HN-stratum with $\cM(2,2)^{\out}_{\un{d}}$ are given as an extension of $\cE_\bullet/\cE_\bullet^\prime$ by $(0\to \cE_0^\prime)$ together with an extension of $\cE_1/\cE_1^\prime$ by $\cE_1^\prime$. We have $\dim(\uExt(\cE_1/\cE_1^\prime, \cE_1^\prime))=g-1-2d_1^\prime$ and $\dim\uExt(\cE_\bullet/\cE_\bullet^\prime,(0\to\cE_0^\prime))= d_0-2d_0^\pprime - (-d_1^\prime-d_0^\pprime)=d_0-d_0^\prime+d_1^\prime$. In total we find $g-1-d_0^\prime-d_1^\prime+d_0$.

%So the class of the intersiection of a stratum with $\cM(2,2)_{out}$ is $\bL^{g-1-d_0^\pprime-3d_1^\prime} [\uPic]^3 [C^{(d_0-d_0^\pprime+d_1^\prime)}].$
Thus the sum over the intersections of the HN-strata with $\cM(2,2)_{\un{d}}^{\out}$ is (set $l=d_0-d_0^\pprime+d_1^\prime$):
\begin{align*}
& \biggl( \sum_{d_0^\prime=\lfloor \frac{d_0}{4}+\frac{\sigma}{2}\rfloor +1}^{\lfloor \frac{d_0}{4}+\frac{\sigma}{2}+(d_0-\sigma)\rfloor} \sum_{l=\lfloor \frac{5d_0}{4}-\frac{\sigma}{2}\rfloor +1-d_0^\prime}^{\lfloor d_0-\sigma\rfloor} \bL^{(g-1)+2d_0-2d_0^\prime-l}[C^{(l)}]\\
&+ \sum_{d_0^\prime=\lfloor \frac{5d_0}{4}-\frac{\sigma}{2}\rfloor+1} \sum_{l=0}^{\lfloor d_0-\sigma\rfloor} \bL^{(g-1)+2d_0-2d_0^\prime-l}[C^{(l)}]\biggr) [\uPic]^3.
\end{align*}

%%%%%%%%%%%%%%%%%%%%%%%%%%%
\subsubsection*{Type $(1,0)\subset (2,2)$} The bounds on the degree are $\mu(\cE_\bullet^\prime)>\mu(\cE_\bullet) \Leftrightarrow d_0^\prime > \frac{d_0+2\sigma}{4}$. 

The quotient $\cE_\bullet/\cE_\bullet^\prime$ has to be semistable of degree $(d_0-d_0^\prime,0)$. This can only happen if $\frac{\sigma}{4} \leq d_0-d_0^\prime \leq \sigma \Leftrightarrow d_0-\frac{\sigma}{4} \geq d_0^\prime \geq d_0-\sigma$. The right hand inequality is automatic because $d_0<2\sigma$.
%Thus $\mu(\ker(\phi^\pprime)) \geq d_0^\prime-d_0 + \sigma < \frac{d_0-d_0^\prime+2\sigma}{3} \Leftrightarrow 4d_0^\prime< 4d_0-\sigma$. 
So we find $\frac{d_0+2 \sigma}{4} < d_0^\prime \leq d_0-\frac{\sigma}{4}$. %, which implies $\sigma < d_0 (<2\sigma)$. 

We have $\dim\uExt(\cE_\bullet/\cE_\bullet^\prime,\cE_\bullet^\prime)= d_0-(g-1)$. 
Using the isomorphisms $$\cM(1,2)_{d_0-d_0^\prime,0}^{\alpha-\sst}\cong \cM(2,1)^{\alpha-\sst}_{0,d_0^\prime-d_0}=\cM(2,1)^{\alpha-\sst}_{2d_0-2d_0^\prime,0},$$ we find that the sum over the  strata is:
\begin{equation*}
\bL^{d_0-(g-1)}[\uPic] \sum_{d_0^\prime=\lfloor \frac{d_0+2\sigma}{4}\rfloor +1}^{\lfloor d_0-\frac{\sigma}{4}\rfloor}[\cM(2,1)^{\alpha-\sst}_{(2d_0-2d_0^\prime,0)}].
\end{equation*}

Strata of this type are contained in $\cM(2,2)^{\fin}_{\un{d}}$, because $(\ker(\phi)\to 0)$ injects into $\cE_\bullet/\cE^\prime_\bullet$, so this cannot be destabilizing. Also if $\rank(\phi)=1$, then $\cE_0^\prime$ injects into $\cE_0/\im(\phi)$ and so $0\to \cE_0/\im(\phi)$ cannot be a destabilizing quotient.
%\end{itemize}

Now we can sum:
\begin{proposition}\label{m22}
Assume that $\alpha=(\sigma,0)$ is good and that $\lfloor\sigma\rfloor$ is even. Then for $d_0>\sigma\geq 2g-2$ we have:
\begin{align*}
[\cM(2,2)_{d_0,0}^{\alpha-\sst}] &= [\cM(2,2)^{\fin}] - \bL^{2d_0-3(g-1)}[\uPic]^2[\Bun_2]\bigg(\sum_{k=\lfloor\sigma-\frac{d_0}{2}\rfloor+1}^{\lceil d_0-\frac{\sigma}{2}\rceil-1}\bL^{-k}\bigg)\\
%% Here use that the (2,1) and (1,0) contributions can be joined using k -> k+(d_0/2)
&- \bL^{d_0+g-1}[\uPic]^2\sum_{k=\lfloor\frac{d_0}{2}\rfloor+1}^{\infty}\bL^{-k} \sum_{l=0}^{\lceil\sigma-\frac{d_0}{2}\rceil-1} [C^{(d_0-\lfloor\sigma\rfloor+2l)}][C^{(\lfloor\sigma\rfloor-2l)}]\\
&- \bL^{2d_0}[\uPic]^3\sum_{k=\lfloor d-\frac{\sigma}{2}\rfloor+1}^\infty \bL^{-k} \sum_{l=0}^{\lceil d_0-\sigma\rceil-1} \bL^{-2l}[C^{(l)}]\\
&+ \bL^{2d_0}[\uPic]^3\sum_{k=\lfloor\sigma-\frac{d_0}{2}\rfloor+1}^{\lceil d_0-\frac{\sigma}{2}\rceil-1}\bL^{-k}\sum_{l=\lfloor d-\sigma\rfloor+1}^\infty \bL^{-2l}[C^{(l)}] \\
&+\bL^{d_0+(g-1)}[\uPic]^3\bigg(\sum_{k=\lfloor\sigma-\frac{d_0}{2}\rfloor+1}^{\lceil\frac{3d_0}{2}-\sigma\rceil-1}\bL^{-k}\sum_{l=\lfloor\frac{3d_0}{4}-\frac{\sigma}{2}-\frac{k}{2}\rfloor+1}^{\lceil d_0-\sigma\rceil-1}\bL^{-l}[C^{(l)}]\\
&+\sum_{k=\lfloor\frac{3d_0}{2}-\sigma\rfloor+1}^\infty \bL^{-k} \sum_{l=0}^{\lceil d_0-\sigma\rceil-1} \bL^{-l}[C^{(l)}]\bigg)\\
&+\bL^{d_0-(g-1)}[\uPic]^3\bigg( \sum_{k=\lfloor\sigma-\frac{d_0}{2}\rfloor+1}^{\lceil d-\frac{\sigma}{2}\rceil-1} \sum_{l=0}^{\lceil \frac{2d-2k-\sigma}{3}\rceil-1} \bL^l [C^{(l)}] \\
&- \sum_{k=\lfloor\frac{\sigma}{2}-\frac{d_0}{4}\rfloor+1}^{\lceil\frac{\sigma}{4}\rceil-1}\sum_{l=\lfloor\sigma-4k\rfloor+1}^{\lceil\frac{2d-4k-\sigma}{3}\rceil-1}\bL^l [C^{(l)}]
- \sum_{k=\lfloor\frac{\sigma}{4}\rfloor+1}^{\lceil\frac{d_0}{2}-\frac{\sigma}{4}\rceil-1} \sum_{l=0}^{ \lceil\frac{2d-4k-\sigma}{3}\rceil-1} \bL^l[C^{(l)}]\bigg).
\end{align*}
\end{proposition}

\subsection{Stacks of chains of rank $(1,2,1)$}
To compute the class of the moduli space of semistable chains rank $(1,2,1)$ and degree $\un{d}$, we can again tensor with a line bundle in order to reduce to the case that $d_2=0$. Moreover, for our application we are only interested in stability parameters of the form $\alpha=(0,\sigma,2\sigma)$, so for simplicity we will only consider such $\alpha$. Our computation will show that such $\alpha$ are good, if $\sigma$ is irrational.
Again we consider the stratification of $\cM(1,2,1)_{\un{d}}$ by saturations of the $\cE_i$ as defined in Section \ref{StratificationSection}. Since $\alpha$ is good all semistable chains $\cE_\bullet$ satisfy $\phi_i\neq 0$. Our descriptions of these strata (Proposition \ref{ChainStrata}) implies:
\begin{lemma}
\begin{enumerate}
\item We have a decomposition 
$$\rcM(1,2,1)_{d_0,d_1,0}= \bigcup_{0\leq l \leq d_0} \rcM(1,2,1)_{\un{d}}^{(1,1,1),(d_0,l,0)} \cup \bigcup_{l\geq 0} \rcM(1,2,1)_{\un{d}}^{(1,1,0),(0,l,0)}.$$  
\item $[\rcM(1,2,1)_{\un{d}}^{(1,1,1),(d_0,l,0)}] = [\uPic]^2 [C^{(l)}][C^{(d_0-l)}] \bL^{d_0-l}$.
\item $[\rcM(1,2,1)_{\un{d}}^{(1,1,0),(0,l,0)}] = [\uPic]^2 [C^{(l)}][C^{(d_0-d_1+l)}] \bL^{(g-1)+d_1-2l}$. This stratum exists for $\max\{0,d_1-d_0\}\leq l$.\hfill$\square$
\end{enumerate}
\end{lemma}
This lemma shows that --- as for $\cM(2,2)$ (Section  \ref{M22Section}) --- the stack $\cM(1,2,1)$ does not define an element in $\K0hat$, because all of the strata with $\phi_1\circ\phi_2=0$ are of the same dimension. However almost all of these are unstable: Set $l=\deg(\im(\cE_2)^\sat)$. 

Then $\cE_2\to \im(\phi_1)^\sat\to 0$ is destabilizing if $\frac{l+3\sigma}{2}>\frac{d_0+d_1}{4}+\sigma \Leftrightarrow l>\frac{d_0+d_1}{2}-\sigma$. 

Also $\cE_2 \to \im(\phi_1)^{\sat} \to \cE_0$ is destabilizing if $\frac{d_0+d_1}{4}+\sigma> d_1-l+\sigma \Leftrightarrow l > \frac{3d_1-d_0}{4}$.

Denote by 
\begin{align*}
\cM(1,2,1)_{\un{d}}^{\fin}&:= \bigcup_{{0\leq l \leq \min\{d_0,\lceil\frac{3d_1-d_0}{4}\rceil-1\}}} \cM(1,2,1)_{\un{d}}^{(1,1,1),(d_0,l,0)} \\
&\cup \bigcup_{0\leq l \leq \lceil\frac{d_0+d_1}{2}-\sigma\rceil-1} \cM(1,2,1)_{\un{d}}^{(0,1,1),(0,l,0)}.
\end{align*}
This is an open substack of $\cM(1,2,1)_{\un{d}}$.
 
\begin{lemma}\label{m121} Let $\un{d}=(d_0,d_1,0)$, $\un{\alpha}=(0,\sigma,2\sigma)$ and suppose that $\sigma\geq 2g-2$ is irrational. Then $\cM(1,2,1)^{\alpha-\sst}_{\un{d}}$ can be non empty only if 
$$d_0+d_1 \leq 4 \sigma,\;\; 3d_0-d_1\leq 4\sigma \textrm{ and } 0\leq d_0 \leq 3d_1\leq 5d_0.$$ 
If these inequalities hold, $\cM(1,2,1)_{\un{d}}^{\alpha-\sst}$ is non-empty and we have %for $d_0<\frac{3}{2}\sigma$
\begin{align*}
[\cM(1,2,1)_{\un{d}}^{\alpha-\sst}]= & [\cM(1,2,1)^{\fin}_{\un{d}}] -  [\uPic]^2\bL^{d_0-(g-1)}\sum_{l=\lfloor \frac{d_1-d_0}{2}+\sigma\rfloor+1}^{\min\{d_0,d_1\}} [C^{(d_0-l)}][C^{(d_1-l)}]\\
&- [\uPic]^2\sum_{l=0}^{\lceil \frac{3d_1-d_0}{4}\rceil -1}  [C^{(l)}][C^{(d_0-l)}] \bL^l.
\end{align*}
\end{lemma}
\begin{proof}
To obtain the necessary conditions  we first list the ranks of canonical subchains:
%\begin{enumerate}
\subsubsection*{Type $(1,2,0)$} There are no semistable chains if: $\mu(\cE_\bullet/\cE_\bullet^\prime)<\mu(\cE_\bullet^\prime) \Leftrightarrow 2\sigma <\frac{d_0+d_1+2\sigma}{3} \Leftrightarrow d_0+d_1> 4\sigma$.
\subsubsection*{Type $(1,0,0)$} There are no semistable chains if: $\mu(\cE_\bullet^\prime)>\mu(\cE_\bullet/\cE_\bullet^\prime) \Leftrightarrow d_0 > \frac{d_1+4\sigma}{3} \Leftrightarrow 3d_0-d_1> 4\sigma$.
\subsubsection*{Type $(1,1,1)$} We always have a subchain $\cE_2 \to \cE_2 \to \cE_0$, so we find the necessary condition $0\leq 3d_1-d_0$.
\subsubsection*{Type $(0,1,0)$} Dually to the previous type we always have a subchain $0\to \ker(\phi_1) \to 0$, so we need $d_1-d_0 \leq (d_0+d_1)/4 \Leftrightarrow 3d_1 \leq 5d_0$.
%\end{enumerate}

%In particular we find $d_0\leq 3d_1\leq 5d_0$ so that $0\leq d_0$. 
Thus we may assume 
\begin{equation}\label{ineq}
d_0+d_1\leq 4\sigma, \;\; 3d_0-d_1\leq 4\sigma \textrm{ and } d_0 \leq 3d_1 \leq 5d_0.
\end{equation}

Let us first exclude strata that do not intersect $\cM(1,2,1)^{\fin}_{\un{d}}$. We list them according to the ranks of $\cE_\bullet^{h-1}$:
%\begin{enumerate}
\subsubsection*{Type $(1,1,1)$} These strata do not intersect $\cM(1,2,1)^{\fin}_{\un{d}}$: If $\cE^\prime_\bullet\subset \cE_\bullet$ is of rank $(1,1,1)$ and $\phi_1^\prime \neq 0$ this holds by definition. If $\phi_1^\prime=0$ then either $\cE_2\to \cE^1_1 \to 0$ or $0 \to 0 \to \cE_0$ is a destabilizing subchain. Thus again the chain does not lie in $\cM(1,2,1)_{\un{d}}^{\fin}$.
\subsubsection*{Type $(0,1,1)$} By definition these strata do not intersect $\cM(1,2,1)_{\un{d}}^{\fin}$.
%\end{enumerate}

Finally we list the HN-strata intersecting $\cM(1,2,1)_{\un{d}}^{\fin}$. We denote the Harder--Narasimhan flag by $\cE_\bullet^\pprime \subset \cE_\bullet^\prime \subset \cE_\bullet$ and list the strata by specifying the rank of the subchains:
%\begin{enumerate}
\subsubsection*{Type $(1,1,0)$}
The bound on the degree is 
$$\mu(\cE_\bullet/\cE_\bullet^\prime)< \mu(\cE_\bullet^\prime) \Leftrightarrow \frac{d_1-d_1^\prime+3\sigma}{2} < \frac{d_0+d_1^\prime+\sigma}{2} \Leftrightarrow d_1^\prime>\frac{d_1-d_0}{2}+\sigma.$$
Also $\cE_\bullet/\cE_\bullet^\prime$ has to be semistable, i.e., $0\leq d_1-d_1^\prime \leq \sigma \Leftrightarrow d_1-\sigma \leq d_1^\prime \leq d_1$ and $\cE_\bullet^\prime$ has to be semistable, i.e., $0\leq d_0-d_1^\prime \leq \sigma \Leftrightarrow d_0-\sigma \leq d_1^\prime \leq d_0$. The lower bounds in these inequalities are automatically satisfied because $d_1-\sigma > \frac{d_1-d_0}{2}+\sigma \Leftrightarrow d_1+d_0 > 4\sigma$ and $d_0-\sigma >\frac{d_1-d_0}{2}+\sigma  \Leftrightarrow 3d_0-d_1>4\sigma$, which we excluded (\ref{ineq}).
Thus the conditions on $d_1^\prime$ are:
$$ \frac{d_1-d_0}{2}+\sigma < d_1^\prime \leq \min\{d_0,d_1\}.$$
We have $\dim(\uExt(\cE_\bullet/\cE_\bullet^\prime,\cE_\bullet^\prime))= d_0 -(g-1)$.

Finally, we claim that HN-strata of this form are contained in $\cM(1,2,1)_{\un{d}}^{\fin}$: 
If $\phi_2\circ\phi_1=0$ the subchain $\cE_2\to \im(\phi_2)^{\sat} \to 0$ is a subchain of $\cE_\bullet/\cE_\bullet^\prime$,  so this cannot be destabilizing. Also the subchain $\cE_2 \to \im(\phi_2)^{\sat} \to \cE_0$ cannot be destabilizing because $0 \to 0 \to \cE_0$ and $\cE_2 \to \im(\phi_1)^{\sat} \to 0$ both have slope $<\mu(\cE_\bullet)$. 
%then $d_1-d_1^\prime > \frac{3d_1-d_0}{4} \Leftrightarrow d_1^\prime < \frac{d_0+d_1}{4}$. Thus $2d_1-2d_0+4\sigma < d_0+d_1 \Leftrightarrow 4\sigma < 3d_0 -d_1 \leq 4\sigma$, a contradiction.

Thus the strata contribute %($l:=\min{d_1,d_0}-d_1^\prime$)
$$ \sum_{d_1^\prime = \lfloor \frac{d_1-d_0}{2}+\sigma  \rfloor +1}^{\min\{d_0,d_1\}} [\uPic]^2[C^{(d_0-d_1^\prime)}][C^{(d_1-d_1^\prime)}]\bL^{d_0-(g-1)}.$$

\subsubsection*{Type $(0,1,0)\subset (1,1,0)$} The bounds on the degree are 
$$d_1^\prime +\sigma > d_0 > \frac{d_1-d_1^\prime+3\sigma}{2}, \textrm{ i.e., } d_1^\prime > \max\{d_0-\sigma, d_1-2d_0+3\sigma\}.$$ 
Also $\cE_\bullet/\cE_\bullet^\prime$ has to be semistable, so as before $d_1 \geq d_1^\prime \geq d_1-\sigma$.
Now $d_1-2d_0+3\sigma \geq d_0-\sigma$ is automatic because $3d_0-d_1\leq 4\sigma$ by (\ref{ineq}). And  $d_1-2d_0+3\sigma \geq d_1-\sigma \Leftrightarrow d_0\leq 2\sigma$ but this is automatic unless $d_1<d_0$ and in that case we already know $d_1-\sigma < d_0-\sigma \leq d_1-2d_0+3\sigma$.
Thus the bounds on the degree are $$d_1-2d_0+3\sigma < d_1^\prime \leq d_1,$$ and this implies $d_0>\frac{3}{2}\sigma$.

For an extension of $\cE_\bullet^\pprime,\cE_\bullet^\prime/\cE_\bullet^\pprime,\cE_\bullet/\cE_\bullet^\prime$ to lie in $\cM(1,2,1)^{\fin}_{\un{d}}$ we need the extension of $\cE_\bullet/\cE_\bullet^\prime$ by  $\cE_\bullet^\prime/\cE_\bullet^\pprime$ to be non-trivial, since otherwise the last map of the chain $\cE_\bullet$ would be $0$.
If this holds, the extension is contained in $\cM(1,2,1)^{\fin}_{\un{d}}$: First, $\cE_2\to \im(\phi_2)^\sat$ has to inject into $\cE_\bullet/\cE^\prime_\bullet$, so it cannot be destabilizing. Second, if $\phi_1\circ\phi_2\neq 0$, then $\cE_2 \to \im(\phi_2)^\sat \to \cE_0$ has to inject into $\cE_\bullet/\cE_\bullet^\pprime$, so again this cannot be destabilizing.

Therefore the strata occur only for $d_0>\frac{3}{2}\sigma$ and in this case their contribution can be calculated as for the $(0,1,0)$-strata to be ($l=d_1-d_1^\prime$):
$$\sum_{l=0}^{\lceil 2d_0-3\sigma\rceil-1}[\uPic]^2 [C^{(l)}][C^{(d_0-l)}] \bL^l.$$

\subsubsection*{Type $(0,1,0)$} The bound on the degree is $$\mu(\cE_\bullet^\prime) > \mu(\cE_\bullet) \Leftrightarrow d_1^\prime > \frac{d_0+d_1}{4}.$$

Also $\cE_\bullet/\cE_\bullet^\prime$ has to be semistable, so that $0\leq d_1-d_1^\prime \leq d_0$, $d_0\leq \frac{d_1-d_1^\prime+3\sigma}{2}$ and $2\sigma \geq \frac{d_0+d_1-d_1^\prime+\sigma}{2}$, i.e., 
$$ \max\{d_1-d_0, d_0+d_1-3\sigma\} \leq d_1^\prime \leq \min\{d_1,d_1-2d_0+3\sigma\}.$$  
We have $d_1-d_0 \leq d_0+d_1-3\sigma \Leftrightarrow d_0 \geq \frac{3}{2}\sigma$ and $d_1\geq d_1-2d_0+3\sigma \Leftrightarrow d_0 \geq \frac{3}{2}\sigma$.
Moreover $\frac{d_0+d_1}{4} > d_0+d_1-3\sigma \Leftrightarrow d_0+d_1 \leq 4 \sigma$ so this is automatic and similarly we already excluded the strata with $d_1^\prime\leq d_1-d_0$.
So we find  
$$ \lfloor \frac{d_0+d_1}{4}\rfloor +1 \leq d_1^\prime \leq \min\{d_1,d_1-2d_0+3\sigma\}.$$
Finally $\dim(\uExt(\cE_\bullet/\cE_\bullet^\prime,\cE_\bullet^\prime))=d_1-d_1^\prime$.

We claim that again, any such HN-stratum is contained in $\cM(1,2,1)_{\un{d}}^{\fin}$: Since $\cE_\bullet/\cE_\bullet^\prime$ is semistable we have $\phi_2\circ \phi_1\neq 0$, and the subchain $\cE_2 \to \im(\phi_2)^\sat \to \cE_0$ must inject into $\cE_\bullet/\cE_\bullet^\prime$, so it cannot be destabilizing. 

%The class of a stratum is $[\uPic]^2[C^{(d_1-d_1^\prime)}][C^{(d_0-d_1+d_1^\prime)}]\bL^{d_1-d_1^\prime}$.

Thus the contribution of these strata is:
$$\sum_{l=\max\{0,2d_0-3\sigma\}}^{\lceil \frac{3d_1-d_0}{4}\rceil -1} [\uPic]^2 [C^{(l)}][C^{(d_0-l)}] \bL^l.$$
%\end{enumerate}
Adding the above contributions we find the claimed formula. The statement that $\cM(1,2,1)^{\alpha-\sst}_{\un{d}}$ is non-empty in this case follows from this formula, since the dimension of each HN-stratum is strictly smaller than $(g-1)^2+2d_0$, which is the dimension of the largest stratum occurring in $\cM(1,2,1)^{\fin}_{\un{d}}$.
\end{proof}

\section{Application: Higgs bundles of rank $4$ and odd degree}\label{HiggsRk4}

From the results of the previous section it is now very easy to deduce the class of $M_4^d$, the moduli space of stable Higgs bundles of rank $4$ and odd degree $d$. This is the aim of this section. In particular the expression we find  gives an explicit formula for the Hodge- and Poincar\'e polynomials of $M_4^d$. 

Let us first note that the moduli spaces $M_4^d$ with $d$ odd are all isomorphic, since by tensoring with a line bundle of fixed degree we can reduce to the case that $d=\pm 1$ and dualization gives an isomorphism $M_4^{-1}\cong M_4^1$. So in the following we will assume that $d=1$. 

We already know from Corollary \ref{HT}, that 
\begin{equation}\label{HTFormula}
[M_{4}^1]=\bL^{16(g-1)+1} \sum_i F_i,
\end{equation}
where $F_i$ are moduli spaces of $\alpha$-semistable chains of some length $r$ rank $\un{n}$ and degree $\un{d}$ with $\sum n_i=4$, $\sum_{i=0}^r d_i-(l-i)n_i(2g-2)=1$ and $\alpha=(0,2g-2,\dots,r(2g-2))$ by Remark \ref{HiggsStability}. 
We used the notation $\cM(\un{n})_{\un{d}}^{\alpha-\sst}$ for the moduli stack of semistable chains and we will write $M(\un{n})_{\un{d}}$ for the corresponding coarse moduli space. 

Since semistability implies stability for Higgs bundles if rank and degree are coprime, the same holds for the moduli spaces of chains occurring as fixed point strata in $M_4^1$. In particular the stability parameter $\alpha$ is not critical so that we may replace $\alpha$ by a good stability parameter $\alpha^\prime$ defining the same moduli space.

Furthermore, since stable Higgs bundles only admit $\bG_m$ as automorphims we know the stack of stable Higgs bundles is a $\bG_m$ gerbe over its coarse moduli space. This gerbe is trivial, because we rank and degree are coprime (see e.g.\ \cite[Lemma 3.10 and Corollary 3.12]{Lectures}). In particular the restriction of the gerbe to the fixed point strata $F_i$ is still trivial. Therefore, as in Example \ref{Bun2} we find that $$[M(\un{n})_{\un{d}}^{\alpha-\sst}]=[\cM(\un{n})_{\un{d}}^{\alpha-\sst}][\bL-1]$$
for all stacks of chains occurring in $M_4^1$.

For all partitions $\un{n}$ with $\sum n_i=4$ we computed the class of the moduli stack $\cM(\un{n})_{\un{d}}^{\alpha-\sst}$ in the previous section and found conditions on $\un{d}$ such that these spaces are non-empty. Let us list the possible range of $\un{d}$ for the different partitions $\un{n}=(n_0,\dots,n_r)$:
%\begin{enumerate}
\subsubsection*{Type $(4)$} Here $d_0=1$.
\subsubsection*{Type $(3,1)$} Here $d_0+d_1-3(2g-2)=1$. By Example \ref{m31} there are no semistable chains of rank $(3,1)$ and degree $(d_0,d_1)$ unless $\sigma\leq d_0-3d_1\leq 3(2g-2)$. Thus the strata $\cM(3,1)_{(e,0)}^{\alpha-\sst}$ contribute for $2g-2 < e < 3(2g-2)$ and $e\equiv (1-\sigma) \md 4$.
\subsubsection*{Type $(1,3)$} In this case $d_0+d_1-(2g-2)=1$. Dualizing and tensoring with line bundles we find $\cM(1,3)_{(d_0,d_1)}^{\alpha-\sst}\cong \cM(3,1)_{(e=3d_0-d_1,0)}^{\alpha-\sst}$. Thus Example \ref{m31} shows that there are no semistable triples of degree $\un{d}$ unless $2g-2 \leq 3d_0-d_1\leq 3(2g-2)$. Thus we find strata $\cM(1,3)_{d_0,d_1)}^{\alpha-\sst}\cong \cM(3,1)_{(e=3d_0-d_1,0)}^{\alpha-\sst}$ for $e\equiv -1-(2g-2) \md 4$ and $2g-2 < e < 3(2g-2)$.
\subsubsection*{Type $(2,2)$} Here $d_0+d_1-4(g-1)=1$. We know (Section \ref{M22Section}) that there are no $\alpha$-semistable chains of this degree unless $0\leq d_0-d_1<2(2g-2)$. As in \ref{M22Section} we use that $\cM(2,2)_{d_0,d_1}^{\alpha-\sst}\cong \cM(2,2)_{(e=d_0-d_1,0)}^{\alpha-\sst}$, and $e=d_0-d_1$ is odd. We thus find that the strata of rank $(2,2)$ occurring are isomorphic to $[\cM(2,2)_{e,0}^{\alpha-\sst}]$ with $0<e<2(2g-2)$ and $e$ odd.
\subsubsection*{Type $(2,1,1)$} We have $d_0+d_1+d_2-5(2g-2)=1$. Write $\overline{d}_0:=d_0-2d_2,\overline{d}_1:=d_1-d_2$, so that we need $\overline{d}_0+\overline{d}_1+4d_2-5(2g-2)=1$, i.e., $\overline{d}_0+\overline{d}_1\equiv 1+(2g-2) \md 4$.

In Example \ref{m211} we have seen there are no semistable chains of rank $(2,1,1)$ unless $0\leq \overline{d}_1$, $\overline{d}_0-\overline{d}_1 \leq 3(2g-2), \overline{d}_0+\overline{d}_1 \leq 5(2g-2), 3(2g-2)\leq 3\overline{d}_0-5\overline{d}_1$,
%Writing $d:=\overline{d}_0+\overline{d}_1$ we find:
%$0\leq \overline{d}_1$, $d-2\overline{d}_1 \leq 3(2g-2), d\leq 5(2g-2), 3(2g-2)\leq 3d-8\overline{d}_1$, 
i.e., $\overline{d}_0+\overline{d}_1\equiv 1+(2g-2) \md 4$ and
\begin{align*}
0 &\leq \overline{d_1} \leq 3g-3,\\
2g-2+\frac{5}{3}\overline{d}_1 &\leq \overline{d}_0 \leq \min(3(2g-2)+\overline{d}_1,5(2g-2)-\overline{d}_1).
\end{align*}

\subsubsection*{Type $(1,1,2)$} Here $d_0+d_1+d_2=1 + 3(2g-2)$. This case is dual to the previous one.
Writing $e_0:=(-d_2+2d_0), e_1:=(-d_1+d_0)$ we find $-(e_0+e_1)+4d_0=1+6g-6$, i.e., we need $e_0+e_1\equiv -1+2g-2 \md 4$ and in this case the bounds on $e_i$ are the same as the bounds on $\overline{d}_i$ of the previous case.

Together with the previous chains we therefore find:
\begin{align*}
\sum_{l=0}^{2g-2} \sum_{k>\frac{5}{3}l+2g-2\atop k+l \equiv 1 \md 2}^{l+6g-6} [\cM(2,1,1)_{k,l,0}^{\alpha-\sst}]
+\sum_{l=2g-1}^{3g-3} \sum_{k>\frac{5}{3}l+2g-2 \atop k+l \equiv 1 \md 2}^{10g-10-l} [\cM(2,1,1)_{k,l,0}^{\alpha-\sst}].
\end{align*}

\subsubsection*{Type $(1,2,1)$} Here $d_0+d_1+d_2-4(2g-2)=1$. Put $\overline{d}_0:=d_0-d_2, \overline{d}_1:=d_1-2d_2$. Then the conditions on the degrees are  $\overline{d}_0+\overline{d}_1+4d_1-4(2g-2)=1$, i.e., we need $\overline{d}_0+\overline{d}_1\equiv 1 \md 4$.
By Lemma \ref{m121} we know that $\cM(1,2,1)^{\alpha-\sst}_{(\overline{d}_0,\overline{d_1},0)}$ is non-empty only if:
$$3\overline{d}_0-\overline{d}_1 \leq 4(2g-2),\; \overline{d}_0+\overline{d}_1 \leq 4\sigma,\; \overline{d}_0 \leq 3\overline{d}_1 \leq 5\overline{d}_0.$$
We put $\overline{d}:=\overline{d}_0+\overline{d}_1$. Then the above inqualities read:
$$ 3\overline{d}-4(2g-2) \leq 4\overline{d}_1,\; \overline{d}\leq 4(2g-2),\; \overline{d} \leq 4\overline{d}_1,\; \overline{d}_1\leq \frac{5}{8}d.$$
So for $0\leq \overline{d} \leq 2(2g-2)$ we have $\overline{d}\leq 4d_1 \leq 4\overline{d}$ and for
$2(2g-2)< \overline{d} < 4(2g-2)$ we have $3d-4(2g-2) \leq 4 \overline{d}_1 \leq 4\overline{d}$ and $\overline{d}=4k+1$, i.e., the strata contribute (here $\overline{d}=4k+1, l=\overline{d}_1$):
$$ \sum_{k=0}^{g-2} \sum_{l=k+1}^{\lfloor\frac{5}{8}(4k+1)\rfloor} \cM(1,2,1)^{\alpha-\sst}_{(4k+1-l,l,0)} +  \sum_{k=g-1}^{2g-3} \sum_{l=3k+1-(2g-2)}^{\lfloor\frac{5}{8}(4k+1)\rfloor} \cM(1,2,1)^{\alpha-\sst}_{(4k+1-l,l,0)}.$$

%Let us abbreviate $K=2g-2$. Then the inequalities for the degrees are:
%$$d_0\leq 2 K,\;\; d_2 >0,\;\; 2d_0+d_2 > 3 K,\;\; 2d_2+d_0 \leq 3K.$$ %checked more than once Jochen 2.5.2010
%Let us rewrite this as
%\begin{enumerate}
%\item For $3g-3 < d_0 \leq 4g-4$ we have $1\leq d_2 < 3g-3 - \frac{d_0}{2}+\frac{1}{8}$.
%\item for $2g-2 < d_0 \leq 3g-3$ we have $\frac{3}{4}+6g-6-2d_0 < d_2 <3g-3 - \frac{d_0}{2}+\frac{1}{8}.$
%\end{enumerate}
\subsubsection*{Type $(1,1,1,1)$} Here $d_0+d_1+d_2+d_3-6(2g-2)=1$. We write $k:=d_2-d_3, l:=d_1-d_2, m:=d_0-d_1$, so that $4d_3+3k+2l+m-6(2g-2)=1$, i.e., we need $3k+2l+m\equiv 1 \md 4$ and $k,l,m\geq 0$.
For semistable chains to exist we need furthermore (Example \ref{M1111}):
\begin{align*}
3k+2l+m& \leq 6(2g-2),\\
k+2l+m & \leq 4(2g-2),\\
k+2l+3m& \leq 6(2g-2).
\end{align*}
These inequalities are equivalent to:
\begin{align*}
0&\leq m \leq 2(2g-2),\\
0&\leq 2l \leq \min\{6(2g-2)-3m,4(2g-2)-m\},\\
0&\leq 3k \leq \min\left\{\begin{array}{l} 6(2g-2)-m-2l, \\ 12(2g-2)-3m-3(2l), \\18(2g-2)-9m-3(2l) \end{array}\right\}.
\end{align*}
Moreover we have $6(2g-2)-3m\leq 4(2g-2)-m \Leftrightarrow (2g-2)\leq m$ and $12(2g-2)-3m-3(2l)\geq 18(2g-2)-9m-3(2l) \Leftrightarrow m \geq (2g-2)$. 

Finally $6(2g-2)-m-2l\leq 12(2g-2)-3m-3(2l) \Leftrightarrow 2l\leq 3(2g-2)-m$ and $6(2g-2)-m-2l\leq 18(2g-2)-9m-3(2l) \Leftrightarrow 2l \leq 6(2g-2)-4m$.
Thus Example \ref{M1111} shows that the sum over these strata is:
\begin{align*}
[\uPic] \Bigg(\sum_{m=0}^{2g-2} \bigg(\sum_{l=0}^{\lceil 3(g-1)-\frac{m}{2}\rceil-1} \sum_{k=0\atop 3k+2l+m\equiv 1 \md 4}^{\lfloor 4(g-1)-\frac{m+2l}{3}\rfloor}[C^{(k)}][C^{(l)}][C^{(m)}]\\
+ \sum_{l=\lceil 3g-3-\frac{m}{2}\rceil}^{\lfloor 4g-4-\frac{m}{2}\rfloor} \sum_{k=0\atop 3k+2l+m\equiv 1 \md 4}^{(8g-8)-m-2l}[C^{(k)}][C^{(l)}][C^{(m)}]\bigg) \\
+ \sum_{m=2g-1}^{4g-4} \bigg(\sum_{l=0}^{\lceil 6g-6-2m\rceil-1} \sum_{k=0\atop 3k+2l+m\equiv 1 \md 4}^{\lfloor 4(g-1)-\frac{m+2l}{3}\rfloor}[C^{(k)}][C^{(l)}][C^{(m)}]\\
+ \sum_{l=6g-6-2m}^{\lceil6g-6-\frac{3m}{2}\rceil-1} \sum_{k=0\atop 3k+2l+m\equiv 1 \md 4}^{(12g-12)-3m-2l}[C^{(k)}][C^{(l)}][C^{(m)}]\bigg)\Bigg).
\end{align*}

%\end{enumerate}
Inserting the above inequalities together into Formula \ref{HTFormula} we find:
\begin{theorem}\label{Higgs4}
The class of the moduli space of stable rank $4$ Higgs bundles of odd degree is
\begin{align*}
&[M_{4}^1]= \\
& \bL^{16(g-1)}(\bL-1) \Bigg( [\Bun_4^{1,\sst}] + \sum_{k=g-1}^{3g-4} [\cM(3,1)_{(2k+1,0)}^{\alpha-\sst}]+ \sum_{k=0}^{2g-3} [\cM(2,2)^{\alpha-\sst}_{(2k+1,0)}] \\
&					+ \sum_{k=0}^{g-2} \sum_{l=k+1}^{\lfloor\frac{5}{8}(4k+1)\rfloor} [\cM(1,2,1)^{\alpha-\sst}_{(4k+1-l,l,0)}] +  \sum_{k=g-1}^{2g-3} \sum_{l=3k+1-(2g-2)}^{\lfloor\frac{5}{8}(4k+1)\rfloor} [\cM(1,2,1)^{\alpha-\sst}_{(4k+1-l,l,0)}]\\
&+\sum_{l=0}^{2g-2} \sum_{k=\lfloor\frac{5}{3}l+2g-2\rfloor+1\atop k+l \equiv 1 \md 2}^{l+6g-6} [\cM(2,1,1)_{k,l,0}^{\alpha-\sst}]+\sum_{l=2g-1}^{3g-3} \sum_{k=\lfloor\frac{5}{3}l+2g-2\rfloor+1 \atop k+l \equiv 1 \md 2}^{10g-10-l} [\cM(2,1,1)_{k,l,0}^{\alpha-\sst}]\\
&+[\uPic] \bigg(\sum_{m=0}^{2g-2} \Big(\sum_{l=0}^{\lceil 3(g-1)-\frac{m}{2}\rceil-1} \sum_{k=0\atop 3k+2l+m\equiv 1 \md 4}^{\lfloor 4(g-1)-\frac{m+2l}{3}\rfloor}[C^{(k)}][C^{(l)}][C^{(m)}]\\
&+ \sum_{l=\lceil 3g-3-\frac{m}{2}\rceil}^{\lfloor 4g-4-\frac{m}{2}\rfloor} \sum_{k=0\atop 3k+2l+m\equiv 1 \md 4}^{(8g-8)-m-2l}[C^{(k)}][C^{(l)}][C^{(m)}]\Big) \\
&+ \sum_{m=2g-1}^{4g-4} \Big(\sum_{l=0}^{6g-6-2m} \sum_{k=0\atop 3k+2l+m\equiv 1 \md 4}^{\lfloor 4(g-1)-\frac{m+2l}{3}\rfloor}[C^{(k)}][C^{(l)}][C^{(m)}]\\
&+ \sum_{l=6g-6-2m+1}^{\lfloor 6g-6-\frac{3m}{2}\rfloor} \sum_{k=0\atop 3k+2l+m\equiv 1 \md 4}^{(12g-12)-3m-2l}[C^{(k)}][C^{(l)}][C^{(m)}]\Big)\bigg)
	\Bigg).
\end{align*}
The classes $[\Bun_4^{1,\sst}],[\cM(3,1)_{(2k+1,0)}^{\alpha-\sst}],[\cM(2,2)^{\alpha-\sst}_{(2k+1,0)}],[\cM(1,2,1)^{\alpha-\sst}_{(4k+1-l,l,0)}] $ and $[\cM(2,1,1)_{k,l,0}^{\alpha-\sst}]$ are given by the formulas in Remark \ref{Bunn}, Example \ref{m31}, Proposition \ref{m22}, Example \ref{m211} and Lemma \ref{m121} for $\alpha=(0,\dots,r(2g-2))$.
\end{theorem}

Since the cohomology of $M_4^1$ is known to have a pure Hodge structure (see e.g. \cite{HauselMirror}), one can immediately read off the Poincar\'e- and Hodge-polynomials of $M_4^1$ from the above theorem using the formulas collected in Section \ref{sym}. For genus $\leq 21$ we evaluated the above formula using {\tt Maple} and found that the result coincides with conjectured result for the Poincar\'e polynomial from \cite{HRV}.

\section{Appendix: Higgs bundles of rank $2$ and $3$}

For completeness we give the formulas for the classes in $\K0hat$ of the spaces of Higgs bundles of rank $n=2,3$. 
For $n=2$ this is contained in Hitchin's original article, where the result is formulated in terms of the Poincar\'e polynomial. For $n=3$ the formula for the Poincar\'e polynomial is due to Gothen \cite{Gothen}.

\begin{theorem}
Let $M_{n}^d$ denote the moduli space of semistable Higgs-bundles of rank $n$ and degree $d$ on $C$.
We have
\begin{align*}
&[M_{2}^1]=\\
&\bL^{4(g-1)+1}P(1) (\frac{P(\bL)-\bL^g P(1)}{(\bL-1)(\bL^2-1)} + \sum_{k=1}^{g-1} [C^{(2k-1)}])= \\
&\bL^{4(g-1)+1}P(1) \left(\frac{P(\bL)}{(\bL-1)(\bL^2-1)} +\frac{P(1)t^{2g-1}}{(1-t^2)(\bL-1)} +\left.\1halb\bigg(Z(C,t)-Z(C,-t)\bigg)\right)\right|_{t=1}
\end{align*}
and
\begin{align*}
&[M_3^1]=\\
&\bL^{9(g-1)-1}P(1)\Bigg(
\frac{P(\bL)P(\bL^2)}{(\bL-1)(\bL^2-1)^2(\bL^3-1)} 
- \frac{\bL^{2(g-1)}(\bL^2+\bL)P(1)P(\bL)}{(\bL-1)^2(\bL^2-1)(\bL^3-1)}\\
&+ \frac{\bL^{3(g-1)+2}P(1)^2}{(\bL-1)^2(\bL^2-1)^2}\\
        &+ \frac{P(1)}{\bL-1} \Big(\sum_{k=0}^{3(g-1)} \sum_{l=0}^{\lfloor{\frac{2k}{3}}\rfloor} [C^{(l)}]\big(\bL^{2(g-1)+k-2l} - \bL^l\big) - \sum_{k=0}^{g-1}\sum_{l=0}^{2k} [C^{(l)}]\big(\bL^{2(g-1)+3k-2l} -\bL^l\big) \Big)\\
        &+ \sum_{k=0}^{2(g-1)}\sum_{l=0 \atop l\not\equiv k \md 3}^{k} [C^{(l)}][C^{(k)}] + \sum_{k=2g-1}^{3(g-1)}\sum_{l=0 \atop l\not\equiv k \md 3}^{6(g-1)-2k} [C^{(l)}][C^{(k)}]\Bigg).
\end{align*}
\end{theorem}
\begin{proof}
We know from Corollary \ref{HT} that 
\begin{equation*}\label{HTFormulaB}
[M_{n}^1]=\bL^{n^2(g-1)+1} \sum_i F_i,
\end{equation*}
where the $F_i$ are the $\alpha$-semistable chains of some length $r$ rank $\un{n}$ and degree $\un{d}$ with $\sum n_i=n$, $\sum_{i=0}^r d_i-(l-i)n_i(2g-2)=1$ and $\alpha=(0,2g-2,\dots,r(2g-2))$ by Remark \ref{HiggsStability}. 

For $n=2$, the fixed point strata for the $\bG_m$ action on $M_{2}^1$ are $[\overline{\Bun}_{2}^{1,\sst}]$ (Example \ref{Bun2}) and spaces of $\alpha$-semistable chains of rank $(1,1)$ and degree $d_0+d_1=1+2g-2$, which exist for $0\leq d_0-d_1 \leq 2g-2$ (see Example \ref{M1111}). Using that $$\sum_{t=0}^{g-1} [C^{2k+1}]t^{2k+1}= \1halb(Z(C,t)-Z(C,-t))-\sum_{k=g}^\infty [C^{(2k+1)}]$$ and for $N>2g-2$ we have $[C^{(N)}]=[\uPic](\bL^{N+1-g}-1)$ we obtain the claimed formula.

% Maple compared the formula for the Poincaré polynomial with Tamas' formula for g=3,13,20,30,47

%Since the $P(t)$ is a polynomial of degree $2g$ the rational function occurring in the residue is holomorphic at $\infty$ so one can compute the residue explicitly by summation over the residues at the other poles, which are of order $1$ and $2$.

The fixed point strata for $n=3$ are of rank $(3),(2,1),(1,2),(1,1,1)$. The class of $\Bun_3^{\sst}$ has been computed in Example \ref{Bunn}. For rank $(2,1)$ semistable chains of degree $d_0+d_1=1+4g-4$ can occur and from Example \ref{Rank21} we know that these are non-empty only if $g-1\leq d_0-2d_1 \leq 4g-4$, i.e.,
$g-1 \leq 1+4(g-1)-3d_1\leq 4(g-1)$. % $d_0-2d_1=g+3l$ for $l=0,\dots,g-2$.

Similarly, for rank $(1,2)$ we have semistable chains of degree $d_0+d_1=1+2g-2$ satisfying $g-1 \leq -(d_1-2d_0) \leq 4g-4$, i.e. $g-1\leq -1-(2g-2)+3d_0 \leq 4g-4$. 

The sum over these strata is (Example \ref{Rank21}):
\begin{align*}
&\sum_{k=1\atop k\not\equiv 0\md 3}^{3(g-1)} (\bL-1)[\cM(2,1)_{g-1+k,0}^{\alpha-\sst}] 
= \frac{P(1)^2}{\bL-1} \sum_{k=1\atop k\not\equiv 0\md 3}^{3(g-1)} \sum_{l=0}^{\lceil \frac{2k}{3}\rceil-1} [C^{(l)}](\bL^{2(g-1)+k-2l} - \bL^l)\\
=&\frac{P(1)^2}{\bL-1} \bigg(\sum_{k=0}^{3(g-1)} \sum_{l=0}^{\lfloor{\frac{2k}{3}}\rfloor} [C^{(l)}](\bL^{2(g-1)+k-2l} - \bL^l) - \sum_{k=0}^{g-1}\sum_{l=0}^{2k} [C^{(l)}](\bL^{2(g-1)+3k-2l} -\bL^l) \bigg).
\end{align*}

For rank $(1,1,1)$ we find $d_0+d_1+d_2=1+6g-6$. Write $l=d_1-d_2$, $k=d_0-d_1$. For semistable chains of rank $(1,1,1)$ and degree $(l+k,l,0)$ to exist we need $0\leq l,k$ and $2l+k\leq 6g-6, l+2k\leq 6g-6$.

Thus the fixed point strata contribute $\rcM(1,1,1)_{l+k,l,0}$ for $2l+k\equiv 1 \md 3$ with $0\leq l,k$ and $2l+k\leq 6g-6, l+2k\leq 6g-6$:

\begin{align*}
\sum_{k=0}^{2(g-1)}\sum_{l=0 \atop l\equiv 1-2k \md 3}^{2g-2} [C^{(l)}][C^{(k)}][\uPic] + \sum_{k=2g-1}^{3(g-1)}\sum_{l=0 \atop l\equiv 1-2k \md 3}^{6(g-1)-2k} [C^{(l)}][C^{(k)}][\uPic] \\
+ \sum_{l=2g-1}^{3(g-1)}\sum_{k=0 \atop 2k\equiv 1-l \md 3}^{6(g-1)-2l} [C^{(k)}][C^{(l)}][\uPic] \\
= [\uPic]\bigg(\sum_{k=0}^{2(g-1)}\sum_{l=0 \atop l\not\equiv k \md 3}^{k} [C^{(l)}][C^{(k)}] + \sum_{k=2g-1}^{3(g-1)}\sum_{l=0 \atop l\not\equiv k \md 3}^{6(g-1)-2k} [C^{(l)}][C^{(k)}] \bigg).
\end{align*}
This proves the claimed formula.
\end{proof}

%%%%%%%%%%%%%%%%%%%%%%%%%%%%%%%%%%%%%%%%%%%%%%%%%%%%%%%%%%%%%%%%%%%%%%%%%%%%%%%%%%%%%%%%%%%%%%%%%%%%%%%
%%% Literaturverzeichnis %%%
%%%%%%%%%%%%%%%%%%%%%%%%%%%%%%%%%%%%%%%%%%%%%%%%%%%%%%%%%%%%%%%%%%%%%%%%%%%%%%%%%%%%%%%%%%%%%%%%%%%%%%%

\end{document}